\pdfoutput=1
\RequirePackage{ifpdf}
\ifpdf % We~are running pdfTeX in pdf mode
\documentclass[pdftex]{sigma}
\else
\documentclass{sigma}
\fi

\numberwithin{equation}{section}

\newtheorem{Theorem}{Theorem}[section]
\newtheorem*{Theorem*}{Theorem}
\newtheorem{Corollary}[Theorem]{Corollary}
\newtheorem{Lemma}[Theorem]{Lemma}
\newtheorem{Proposition}[Theorem]{Proposition}

\theoremstyle{definition}
\newtheorem{Definition}[Theorem]{Definition}

\newtheorem{Remark}[Theorem]{Remark}

\usepackage{mathrsfs}
\usepackage{physics}
\usepackage{tikz-cd}
\usepackage{bm}
\usetikzlibrary{backgrounds}
\usetikzlibrary{calc}
\usetikzlibrary{hobby}
\usetikzlibrary{decorations.markings}
\usetikzlibrary{arrows.meta}
\usetikzlibrary{patterns}
\usepackage{array}
\usepackage{afterpage}

\DeclareMathOperator{\ram}{ram}
\DeclareMathOperator{\slope}{slope}
\DeclareMathOperator{\Hom}{Hom}
\DeclareMathOperator{\End}{End}

\DeclareMathOperator{\irr}{irr}

\DeclareMathOperator{\Levels}{Levels}
\DeclareMathOperator{\GL}{GL}
\DeclareMathOperator{\SL}{SL}

\newcommand{\wh}{\widehat}
\newcommand{\IC}{\mathbb C}
\newcommand{\IP}{\mathbb P}

\newcommand{\Ical}{\mathcal{I}}

\newcommand{\Ccal}{\mathcal{C}}

\newcommand{\cir}[1]{\bigl\langle #1 \bigr\rangle}

\renewcommand{\proofname}{Proof}

\begin{document}

\allowdisplaybreaks

\newcommand{\arXivNumber}{2409.12864}

\renewcommand{\PaperNumber}{027}

\FirstPageHeading

\ShortArticleName{Basic Representations of Genus Zero Nonabelian Hodge Spaces}

\ArticleName{Basic Representations of Genus Zero Nonabelian \\ Hodge Spaces}

\Author{Jean DOU\c{C}OT}

\AuthorNameForHeading{J.~Dou\c{c}ot}

\Address{``Simion Stoilow'' Institute of Mathematics of the Romanian Academy, \\ Calea Grivi\c{t}ei~21,
010702-Bucharest, Sector 1, Romania}
\Email{\mail{jeandoucot@gmail.com}}
\URLaddress{\url{https://www.normalesup.org/~doucot/}}

\ArticleDates{Received May 28, 2025, in final form February 13, 2026; Published online March 20, 2026}

\Abstract{In some previous work, we defined an invariant of genus zero nonabelian Hodge spaces taking the form of a diagram. Here, enriching the diagram by fission data to obtain a~refined invariant, the enriched tree, including a partition of the core diagram into~$k$~subsets, we show that this invariant contains sufficient information to reconstruct $k+1$ different classes of admissible deformations of wild Riemann surfaces, that are all representations of one single nonabelian Hodge space, so that the isomonodromy systems defined by these representations are expected to be isomorphic. This partially generalises to the case of arbitrary singularity data the picture of the simply-laced case featuring a diagram with a complete $k$-partite core. We illustrate this framework by discussing different Lax representations for Painlev\'e equations.}

\Keywords{irregular connections; Lax representations; Painlev\'e equations; Fourier transform; nonabelian Hodge diagrams}

\Classification{34M40; 44A10; 32C38; 32G34}

\section{Introduction}

\subsection{General motivation: classification of (wild) nonabelian Hodge spaces}

The main aim of this work is to formulate a convenient combinatorial way of identifying many conjectural isomorphisms between different genus zero nonabelian Hodge spaces, which should in turn extend to isomorphisms between the related isomonodromy systems, generalising known (multi)dualities between different Lax representations of Painlev\'e-type equations to the case of arbitrary irregular types.

The term `nonabelian Hodge space', introduced in \cite[Definition~7]{boalch2018wild}, refers to a class of manifolds that can be seen as moduli spaces for different types of objects. On one hand, they are moduli spaces of meromorphic connections with possibly irregular singularities on curves (the de Rham side). Via the Riemann--Hilbert--Birkhoff correspondence, they are also moduli spaces of generalised monodromy data, known as wild character varieties \cite{boalch2014geometry, boalch2015twisted} (the Betti side), and, via the wild nonabelian Hodge correspondence \cite{biquard2004wild}, as moduli spaces of irregular Higgs bundles (the Dolbeault side). Furthermore, they carry rich (symplectic, hyperk\"ahler) geometric structures, and give rise to (isospectral and isomonodromic) integrable systems, including many classical integrable systems of interest in mathematical physics such as all Painlev\'e equations (see, e.g., \cite{boalch2018wild} for an overview).

A nonabelian Hodge space $\mathcal M$ depends on the choice of a tuple $(\Sigma, \mathbf a, \bm\Theta,\bm{\mathcal C})$, that we call a~\textit{wild Riemann surface with boundary data}. Here the triple $(\Sigma, \mathbf a, \bm\Theta)$ is a \textit{wild Riemann surface}~\cite{boalch2014geometry, boalch2021topology}, meaning that $\Sigma$ is a compact Riemann surface, $\mathbf a=\{a_1,\dots, a_m\}\subset \Sigma$ is a finite set of points on $\Sigma$ corresponding to the singularities of the connections giving rise to points of $\mathcal M_{\rm dR}$, the global \textit{irregular class} $\bm\Theta$ encodes the type of singularities that these connections have at each singular point, and $\bm{\mathcal C}$ is a collection of conjugacy classes, encoding the \emph{formal monodromies} of the connections.

Remarkably, it turns out that there exist some isomorphisms between moduli spaces defined by different wild Riemann surfaces with boundary data, corresponding to connections with different ranks, number of singularities and pole orders. It is expected that when this happens the moduli spaces are really isomorphic as nonabelian Hodge spaces, i.e., all the structures of the nonabelian Hodge package match in the natural way. In other words, we would like to view those different wild Riemann surfaces with boundary data as different \textit{realizations}, or \textit{representations}, of one single abstract nonabelian Hodge space $\mathcal M$, which raises the question of their classification. Notice that here we will only be considering representations coming from the case of vector bundles, however wild nonabelian Hodge spaces arising from vector bundles may also be isomorphic to nonabelian Hodge spaces arising from principal $G$-bundles, for~$G$ a~complex reductive group different from $\GL_n(\mathbb C)$, see, for example,~\cite{boalch2013symmetric}.

Varying the wild Riemann surface, one can consider moduli spaces of irregular connections in families, giving rise to isomonodromic deformations \cite{boalch2014geometry}. If the connections are written down explicitly in coordinates, this gives systems of nonlinear PDEs which include, among the simplest nontrivial examples, all Painlev\'e equations \cite{jimbo1981monodromyI} (see also \cite{bertola2023hamiltonian,gaiur2023isomonodromic,marchal2025hamiltonian} for works investigating more systematically the question of obtaining explicit Hamiltonians for isomonodromic deformations). In this picture, the wild Riemann surface $(\Sigma, \mathbf a, \bm\Theta)$ giving rise to a given isomonodromy system corresponds to what is called a \textit{Lax representation}, \textit{Lax pair}, or \textit{linearisation} for the system in the literature on integrable systems, albeit in a more abstract sense than is often the case (typically the term Lax pair is used to refer to an explicit parametrisation in coordinates of the coefficients of the connections). The existence of nontrivial isomorphisms between different nonabelian Hodge spaces is related to the existence of different Lax representations for many Painlev\'e-type equations, and the question of the classification of nonabelian Hodge spaces extends to the problem of finding all possible Lax representations for a given isomonodromy system.

\subsection{Diagrams and isomorphisms between moduli spaces of connections}

Interestingly, it turns out that in many genus zero cases, i.e., when the underlying Riemann surface is $\mathbb P^1(\mathbb C)$, there is a clean combinatorial way to pass between different representations.

Possibly the simplest instance where such ``dualities'' between different moduli spaces of connections on $\mathbb P^1$ occur is Harnad duality \cite{adams1988isospectral,alameddine2024explicit,harnad1994dual, yamakawa2011middle}, which relates dual connections on trivial vector bundles over the affine line of the form
$
{\rm d}-\bigl(A+P(Y-z)^{-1}Q\bigr){\rm d}z$ and ${\rm d}+\bigl(Y+Q(A-z)^{-1}P\bigr){\rm d}z$,
where $P\colon W\to V$ and $Q\colon V\to W$ are linear maps, where $V$, $W$ are finite dimensional complex vector spaces, and $A\in \End(V)$, $Y\in \End(W)$. The connections on both sides have an irregular singularity at infinity, with a pole of order two, as well as regular singularities at finite distance, whose positions are given by the eigenvalues of~$A$ and~$Y$. Harnad showed that the isomonodromic deformations of these two connections, governed by the JMMS equations~\cite{jimbo1980density}, are equivalent. The duality induces an isomorphism between the de Rham moduli spaces corresponding to both sides.

Harnad duality can be conveniently described in a graphical way as follows. The linear maps $A$, $Y$, $P$, $Q$ can be viewed as representations of the quiver drawn below, and the duality corresponds to exchanging the roles of the two vertices:
$$
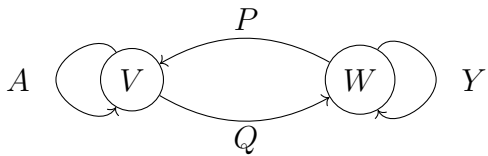
\begin{tikzpicture}
\tikzstyle{vertex}=[circle,fill=black,minimum size=6pt,inner sep=0pt]
\node[circle, draw] (A) at (0,0){$V$};
\node[circle, draw] (B) at (3,0){$W$};
\draw[->] (A) to[bend right] (B);
\draw[->] (B) to[bend right] (A);
\draw[->] (A) to[in=90, out=120] (-1,0) to[in=-120, out=-90] (A);
\draw[->] (B) to[in=90, out=60] (4,0) to[in=-60, out=-90] (B);
\draw (-1.5, 0) node {$A$};
\draw (4.5, 0) node {$Y$};
\draw (1.5,0.8) node {$P$};
\draw (1.5,-0.8) node {$Q$};
\end{tikzpicture}
$$
When $A$ and $Y$ are semisimple, the vector spaces $V$ and $W$ can further be split into the corresponding eigenspaces, giving rise to a complete bipartite graph, and Harnad duality amounts to exchanging the two subsets of vertices.

This duality was enlarged by Boalch \cite{boalch2008irregular,boalch2012simply, boalch2016global} in the \textit{simply-laced} case, which is the case of connections on $\mathbb P^1$ having one irregular singularity, located at infinity, with an (unramified) pole of order at most~3, possibly together with regular singularities at finite distance: now there is a~``multiduality'' between more than two moduli spaces of connections.

This multiduality can also be conveniently encoded in a graphical way, in terms of a \textit{supernova quiver}: it consists of a core, which is a complete $k$-partite graph for some integer~$k$ (i.e., the core is a disjoint union of $k$ subsets of vertices, and each vertex in one subset is linked by one edge to every vertex belonging to a different subset) to which are then glued \textit{legs}, each leg being a linear (i.e., type A) quiver.

The main statement is then as follows \cite{boalch2012simply}: such a supernova quiver, together with some extra data, canonically determines $k+1$ moduli spaces of connections of simply-laced type, with in general different ranks and number of singularities, all isomorphic to each other. In each representation, the vertices of the core diagram are in one-to-one correspondence with the exponential factors of the connections. There is one \emph{generic} representation, of maximal rank, with only one (irregular) singularity, located at infinity, and every $i\in \{1,\dots, k\}$ defines a~\emph{nongeneric} representation, in which the vertices of the $i$-th subset of the core diagram correspond to regular singularities at finite distance, while all other core vertices correspond to exponential factors at infinity. In other words, one can read from the diagram $k+1$ representations of the same de Rham moduli space (see the example on the figure below). More concretely, the generic representation corresponds to connections of the form
$
{\rm d}-(Az+B){\rm d}z$,
%\label{eq:simply_laced_connection}
with $A$, $B$ constant matrices, and $A$ semisimple, with $k$ distinct eigenvalues in one-to-one correspondence with the subsets of the partition of the core diagram.
Furthermore, the equations governing the isomonodromic deformations of each representation are equivalent, so we can say that the different representations correspond to isomorphisms of isomonodromy systems.

\begin{figure}[!ht]
\centering
\begin{tikzpicture}[scale=1.1]
\tikzstyle{vertex}=[circle,fill=cyan,minimum size=6pt,inner sep=0pt]
\tikzstyle{vertex_2}=[circle,fill=purple,minimum size=6pt,inner sep=0pt]
\tikzstyle{vertex_3}=[circle,fill=orange,minimum size=6pt,inner sep=0pt]
\tikzstyle{vertex_empty}=[circle,fill=white,minimum size=5pt,draw,inner sep=0pt]
\begin{scope}
\node[vertex] (A) at (-1,0){};
\node[vertex_2] (B) at (0,1){};
\node[vertex] (C) at (1,0){};
\node[vertex_2] (D) at (0,-1){};
\node[vertex_3] (E) at (0,0){};
\draw (A)--(B)--(C)--(D)--(A);
\draw (A)--(E)--(C);
\draw (B)--(E)--(D);
\end{scope}
\begin{scope}[xshift=3cm]
\node[vertex] (A) at (-1,0){};
\node[vertex_2] (B) at (0,1){};
\node[vertex] (C) at (1,0){};
\node[vertex_2] (D) at (0,-1){};
\node[vertex_3] (E) at (0,0){};

\node[vertex_empty] (B1) at (-0.5,0.5){};
\node[vertex_empty] (E1) at (-0.5,0){};
\node[vertex_empty] (D1) at (-0.5,-0.5){};

\node[vertex_empty] (B2) at (0.5,0.5){};
\node[vertex_empty] (E2) at (0.5,0){};
\node[vertex_empty] (D2) at (0.5,-0.5){};
\draw (A)--(B1);
\draw (A)--(E1);
\draw (A)--(D1);

\draw[dotted] (B1)--(B)--(B2);
\draw[dotted] (E1)--(E)--(E2);
\draw[dotted] (D1)--(D)--(D2);

\draw (C)--(B2);
\draw (C)--(E2);
\draw (C)--(D2);

\draw (B)--(E)--(D);
\end{scope}
\begin{scope}[xshift=6cm]
\node[vertex] (A) at (-1,0){};
\node[vertex_2] (B) at (0,1){};
\node[vertex] (C) at (1,0){};
\node[vertex_2] (D) at (0,-1){};
\node[vertex_3] (E) at (0,0){};

\node[vertex_empty] (A1) at (-0.5,0.5){};
\node[vertex_empty] (E1) at (0,0.5){};
\node[vertex_empty] (C1) at (0.5,0.5){};

\node[vertex_empty] (A2) at (-0.5,-0.5){};
\node[vertex_empty] (E2) at (0,-0.5){};
\node[vertex_empty] (C2) at (0.5,-0.5){};

\draw (B)--(A1);
\draw (B)--(E1);
\draw (B)--(C1);

\draw (D)--(A2);
\draw (D)--(E2);
\draw (D)--(C2);

\draw[dotted] (A1)--(A)--(A2);
\draw[dotted] (E1)--(E)--(E2);
\draw[dotted] (C1)--(C)--(C2);

\draw (A)--(E)--(C);

\end{scope}
\begin{scope}[xshift=9cm]
\node[vertex] (A) at (-1,0){};
\node[vertex_2] (B) at (0,1){};
\node[vertex] (C) at (1,0){};
\node[vertex_2] (D) at (0,-1){};
\node[vertex_3] (E) at (0,0){};

\node[vertex_empty] (A1) at (-0.5, 0){};
\node[vertex_empty] (B1) at (0, 0.5){};
\node[vertex_empty] (C1) at (0.5, 0){};
\node[vertex_empty] (D1) at (0, -0.5){};

\draw (A1)--(E)--(C1);
\draw (B1)--(E)--(D1);

\draw[dotted] (A)--(A1);
\draw[dotted] (B)--(B1);
\draw[dotted] (C)--(C1);
\draw[dotted] (D)--(D1);

\draw (A)--(B)--(C)--(D)--(A);

\end{scope}
\end{tikzpicture}
\caption{Example of different readings of a (core) diagram in the simply-laced case, cf.\ \cite[p.~11]{boalch2008irregular}. On the left, a complete 3-partite graph, with the different colours of the vertices corresponding to the~3 subsets of the partition. The generic reading corresponds to rank~5 connections with one pole of order 3. The three other figures on the right correspond to the nongeneric readings, each reading being obtained by singling out one of the subsets of the partition and interpreting as coming from regular singularities at finite distance. The nongeneric readings correspond to connections with ranks 3, 3, 4, and pole orders~${3+1+1}$, $3+1+1$, $3+1$ respectively from left to right.}
\label{fig:simply_laced_different_readings}
\end{figure}

\subsection[Basic operations on connections on P\^{}1 and weak representations]{Basic operations on connections on $\boldsymbol{\mathbb P^1}$ and weak representations}\label{subsec:weak_reps_intro}

From a more abstract point of view, Harnad duality is underlain by the Fourier--Laplace transform \cite{boalch2008irregular, boalch2012simply,harnad1994dual, kawakami2018degeneration}. It is an operation on (almost all) irreducible algebraic connections on Zariski open subsets of the affine line induced by the automorphism of the Weyl algebra of differential operators $A_1:=\mathbb C[z,\partial_z]$ defined by $z\mapsto -\partial_z$, $\partial_z\mapsto z$. Applied to such a connection, it gives another connection with in general different rank, number of singularities and pole orders.

Similarly, underlying the enlarged symmetry of the simply-laced case is the fact that the Fourier transform is actually part of a larger group of transformations acting on connections on the affine line: indeed, any matrix
\[A=\left(\begin{matrix}
a & b\\ c & d
\end{matrix}\right)
\]
in $\mathrm{SL}_2(\mathbb C)$ defines an automorphism of the Weyl algebra by $z\mapsto az+b\partial_z$, $\partial_z\mapsto cz+d\partial_z$ (this simply follows from the fact that this transformation preserves the defining commutation relation $[\partial_z,z]=1$ of $A_1$). This induces an action of the group $\mathrm{SL}_2(\mathbb C)$ on the set of isomorphism classes of irreducible algebraic connections on Zariski open subsets of the affine line $\mathbb C$ (excluding rank 1 connections with only a pole at infinity of order less than~2). The property of being of simply-laced type is preserved by this action, and the $k+1$ readings of the supernova quiver in the simply-laced case correspond to the different singularity patterns that appear in a given orbit.

More generally, there is a larger class of basic operations on irreducible irregular connections on~$\mathbb P^1$, consisting of
\begin{itemize}\itemsep=0pt
\item M\"obius transformations of $\mathbb P^1$,
\item Twists by rank one connections,
\item $\SL_2(\mathbb C)$ transformations.
\end{itemize}
We will often refer to the $\SL_2(\mathbb C)$ transformations as \emph{symplectic transformations} to avoid possible confusions with M\"obius transformations. These operations play an essential role in the extension due to Deligne--Arinkin \cite{arinkin2010rigid, deligne2006letter} of the Katz algorithm for rigid local systems~\cite{katz1996rigid} to the case of irregular connections, which provides a way to reduce any irreducible rigid connection on~$\mathbb P^1$ to the trivial rank one connection by repeated application of such operations.

M\"obius transformations and twists induce in a straightforward way isomorphisms of nonabelian Hodge spaces, in particular they do not change their Betti side description as wild character varieties. On the other hand, the Fourier transform is expected to induce nontrivial isomorphisms of wild character varieties.

It is known in the case of regular singularities and in the simply-laced case that the wild character varieties on both sides of the Fourier transform are isomorphic~\cite{boalch2016global}. In the case of regular singularities it has also been shown that the Fourier transform induces a hyperk\"ahler isometry~\cite{szabo2015plancherel} between the nonabelian Hodge moduli spaces.
Some work remains to be done to deal with the general case, due to the fact that determining the Stokes data of the Fourier transform of an irregular connection on $\mathbb P^1$ is quite involved \cite{dagnolo2020topological, doucot2025topological,hien2015local, hohl2022d_modules,malgrange1991equations, mochizuki2010note,mochizuki2018stokes, sabbah2016differential}. For this reason, since in the meantime we still want to view genus zero wild Riemann surfaces with boundary data related by basic operations as representations of the same nonabelian Hodge space, relying on the fact that basic operations on connections induce well-defined counterparts at the level of the corresponding wild Riemann surfaces with boundary data, we make the following definition.

\begin{Definition}
Let $\mathcal M=\mathcal M(\bm{\Sigma_0})$ be a genus zero nonabelian Hodge space, with $\bm{\Sigma_0}$ a wild Riemann surface with boundary data on $\mathbb P^1$. A \emph{weak representation} of $\mathcal M$ is a wild Riemann surface with boundary data $\bm{\Sigma}$ that can be obtained from $\bm{\Sigma_0}$ by successive application of basic operations or admissible deformations.
\end{Definition}

In the rest of the article, we will allow ourselves to simply use the term representation to refer to weak representations. Given a (weak) representation $\mathbf \Sigma$ of a nonabelian Hodge space $\mathcal M$, we say that the representations obtained from $\bm\Sigma$ by the action of a symplectic transformation defined by an element of $\SL_2(\mathbb C)$, are its \emph{nearby representations}. Finally, we will call \emph{basic representations} of a nonabelian Hodge space the nearby representations of a representation of minimal rank.

\subsection{Main results: from enriched trees to classes of nearby representations} The main purpose of this article is to develop a~combinatorial way to find many weak representations of any nonabelian Hodge space in genus zero. Namely, we generalise the picture of obtaining different representations of a~nonabelian Hodge space via different readings of a~diagram from the simply-laced case to the case of genus zero nonabelian Hodge spaces with arbitrary singularity data. More specifically, we define some combinatorial data generalising the supernova quivers of the simply-laced case, from which several (expected to be) isomorphic nonabelian Hodge spaces can be explicitly obtained, corresponding to the different types of singularities appearing in an arbitrary $\SL_2(\mathbb C)$ orbit, i.e., to different types of nearby representations of a~given wild nonabelian Hodge space.

We build on our previous work \cite{doucot2021diagrams}, in which a construction of a \emph{diagram} associated to an irregular connection on a Zariski open subset of $\mathbb P^1$ was formulated, for the general case of arbitrary irregular types. Furthermore it was shown that the dimension of the wild character variety is given by a formula involving the Cartan matrix of the diagram, in agreement with the philosophy of `global Lie theory' advocated in \cite{boalch2016global} suggesting to view nonabelian Hodge spaces as some global analogues of Lie groups. The main property of the diagram is that it is invariant under the symplectic $\SL_2(\mathbb C)$ action, so that connections with different types of singularities have the same diagram.

However, it turns out that, unlike for the simply-laced case, the diagram is not quite enough to reconstruct the corresponding nonabelian Hodge spaces, so some extra data will be needed. In particular, the core diagram is not a complete $k$-partite graph any more in the general case, and these extra data will include a canonical partition of the core diagram.

To obtain the appropriate extra data, as well as to give a precise meaning to the idea of `different types' of singularity data, as was already sketched in \cite[Section~5.3]{doucot2021diagrams}, we rely on our recent study of \textit{admissible deformations} of wild Riemann surfaces \cite{boalch2025twisted, doucot2023topology, doucot2022local, doucot2024moduli}. Being admissible deformations of each other defines an equivalence relation on the set of all possible wild Riemann surfaces with boundary data of such connections, such that up to isomorphism the wild character variety $\mathcal M_B(\Sigma, \mathbf a, \bm{\Theta},\bm{\mathcal C})$ is invariant under admissible deformations. This motivates the following definition.

\begin{Definition}
Let $\mathcal M$ be a nonabelian Hodge space. A \emph{class of representations} of $\mathcal M$ is an equivalence class of admissible deformations of weak representations of $\mathcal M$.
\end{Definition}

In \cite{boalch2025twisted}, we showed that a class of admissible deformations of wild Riemann surfaces is completely characterised by a pair $(g,\mathbf F)$, where $g$ is the genus of the underlying Riemann surface~$\Sigma$,~$\mathbf F$ is a \emph{fission forest}, i.e., a collection of \emph{fission trees} (see \cite[Definition~3.18]{boalch2025twisted}). We can thus write~${\mathcal M=\mathcal M(g, \mathbf F, \bm{\mathcal C})}$ up to isomorphism, where $\bm{\mathcal C}$ is the datum of a conjugacy class in~$\GL_n(\mathbb C)$ for each leaf of some fission tree of $\mathbf F$, with $n \geq 1$ an integer, the \emph{multiplicity} of the leaf. Here we only consider the case $g=0$, so we will drop the genus from the notation and write $\mathcal M=\mathcal M(\mathbf F, \bm{\mathcal C})$.

The data that will allow us to reconstruct the classes of nearby representations consist of an \textit{enriched tree}. To explain what it is, recall from \cite[Section~3.4]{boalch2025twisted} that a fission tree is a~tuple~${\mathbf T=(\mathcal T, \mathbb V, \mathbb A, \mathbb L,h,n)}$ where $\mathcal T$ is a metrised tree, $\mathbb V\subset \mathcal T$ is a subset whose points are the \textit{vertices} of the tree, $\mathbb A$ and $\mathbb L$ are two distinguished subsets of $\mathbb V$ such that $\mathbb L\subset \mathbb A$, $h\colon \mathcal T\to \mathbb R_{\geq 0}$ is a~function that we call the $height$, such that for any vertex $v\in \mathbb V$, $h(v)\in \mathbb Q$, and $h^{-1}(0)=:N$ is the set of leaves of the tree, and $n\colon N\to \mathbb Z_{\geq 1}$ consists of the datum of an integer multiplicity~${n_i\geq 1}$ for each leaf $i$ (see Section~\ref{subsec:admissible_defs_and_forests} for more details on fission trees).

\begin{Definition}
An enriched tree is a pair $\mathscr T=(\mathbf T, \bm{\mathcal C})$, where
\begin{itemize}\itemsep=0pt
\item $\mathbf T=(\mathcal T, \mathbb V, \mathbb A,\mathbb L, h, n)$ is a \textit{short fission tree}, i.e., a truncated fission tree in the sense of \cite[Section~3.5]{boalch2025twisted} such that all its vertices different from the root have height $\leq 2$.
\item If $N$ denotes the set of leaves of $\mathbf T$, $\bm{\mathcal C}$ is the datum of a conjugacy class $\mathcal C_i\in \GL_{n_i}(\mathbb C)$ for each $i\in N$, where $n_i\geq 1$ is the multiplicity of the leaf $i$.
\end{itemize}
\end{Definition}

We show that enriched trees provide an invariant of irregular connections on $\mathbb P^1$ under symplectic transformations.

\begin{Theorem}[Proposition~\ref{prop: definition_of_enriched_trees}, Definition~\ref{def:enriched_tree}, Corollary~\ref{cor:symplectic_invariance_enriched_tree}]
Any algebraic connection $(E,\nabla)$ on a Zariski open subset of $\mathbb P^1$ canonically determines an enriched tree $\mathscr T(E,\nabla)$, which depends only on its wild Riemann surface with boundary data. If $(E,\nabla)$ is irreducible, and is not a rank one connection with a pole of order less than~$2$ at infinity, then for $A\in \SL_2(\mathbb C)$, $A\cdot (E,\nabla)$ is well-defined and we have
$
\mathscr T(A\cdot (E,\nabla))=\mathscr T(E,\nabla)$.
\end{Theorem}

Our main result is that this invariant is powerful enough to allow us to reconstruct all classes of nearby representations.

\begin{Theorem}[Corollary~\ref{cor:we_get_classes_of_nearby_reps}]
\label{thm: intro_reconstruction_representations_from_tree}
Let $\mathscr T=(\mathbf T, \bm{\mathcal C})$ be an enriched tree, and $k$ be the number of principal subtrees of $\mathbf T$, i.e., the number of vertices of height $2$. Assuming that the nonabelian Hodge space $\mathcal M(\{\mathbf T\}, \bm{\mathcal C})$ is nonempty, $\mathscr T$ canonically determines in an explicit way $k$ classes of representations $(\mathbf F_i,\bm{\mathcal C_i})$, for $i\in \{1,\dots, k\}$, of $\mathcal M(\{\mathbf T\}, \bm{\mathcal C})$, such that, if $(\bm\Theta, \bm\Ccal)$ is any wild Riemann surface with boundary data on $\mathbb P^1$ with enriched tree $(\mathbf T, \bm{\mathcal C})$, then $(\{\mathbf T\}, \bm{\mathcal C})$ and $(\mathbf F_i,\bm{\mathcal C_i})$, for $i\in\{1,\dots, k\}$ are all the classes of its nearby representations of $\mathcal M(\{\mathbf T\}, \bm{\mathcal C})$.
\end{Theorem}

While we refer the reader to the body of the article for the detailed construction of the classes of representations (defined in Definition~\ref{def:nearby_representations}), let us now briefly summarise their properties and discuss the links with the diagrams of \cite{doucot2021diagrams}. Recall that if $(E,\nabla)$ is an algebraic connection on a Zariski open subset of $\mathbb P^1$, at each of its singularities it possesses a finite set of \emph{Stokes circles}, which basically correspond to the different exponential factors $e^q$ appearing in the asymptotics of the horizontal sections of the connection close to the singularity. In our language, they are encoded by the (modified) global irregular class $\bm{\breve\Theta}$ which is the collection of all Stokes circles of~$(E,\nabla)$, with some integer multiplicities.

\begin{Theorem}[Propositions~\ref{prop:rank_of_representations} and \ref{prop:number_of_singularities}] Keeping the setup of Theorem~{\rm\ref{thm: intro_reconstruction_representations_from_tree}}, for $i\in\{1,\dots,\allowbreak k\}$, let $\mathbf T_i$ be the $i$-th principal subtree of $\mathbf T$, and $N_i$ denote its set of leaves, so that we have a~partition $N=N_1 \sqcup \dots \sqcup N_k$. Furthermore, let $N_i^+\subset N_i$ be the subset of leaves in $N_i$ having an ancestor vertex in $\mathbf T_i$ of height $1<h<2$, and $N_i^-:=N_i\smallsetminus N_i^+$.

Let $(E,\nabla)$ be an irreducible irregular connection on $\mathbb P^1$ with enriched tree $(\mathbf T, \bm{\mathcal C})$, and $(E',\nabla')\allowbreak\in \SL_2(\mathbb C)\cdot (E,\nabla)$ an element in its orbit. We have
\begin{enumerate}\itemsep=0pt
\item[$(1)$] There is a bijection between $N$ and the set of Stokes circles of $(E',\nabla')$. Furthermore,
\begin{itemize}\itemsep=0pt
\item
If $(E',\nabla')$ is in the generic class of representations, all its Stokes circles are of \emph{generic form}, i.e., are at infinity and are of slope $\leq 2$.
\item Otherwise, if $(E',\nabla')$ is in the $i$-th nongeneric class of representations for $i\in\{1,\allowbreak\dots k\}$,
\begin{itemize}\itemsep=0pt
\item the Stokes circles corresponding to elements of $N_i^-$ via the bijection are at finite distance $\big($i.e., different from $\infty \in \mathbb P^1\big)$,
 \item the Stokes circles corresponding to elements of $N_i^+$ are at infinity and of slope~$>2$,
 \item the Stokes circles corresponding to $N_j$ for every $j\neq i$ are of generic form.
\end{itemize}
\end{itemize}
\item[$(2)$] If $(E',\nabla')$ is in the $i$-th nongeneric class of representations, the number of its singularities at finite distance is equal to the number of principal subtrees of $\mathbf T_i$ with root at height~$1$.
\item[$(3)$] The ranks of the classes of representations are as follows. For any leaf $I\in N_i$, let $r_I$ denote the lowest common multiple of the denominators $($when written in their reduced form$)$ of all the heights of its ancestors in $\mathbf T_i$, and let $n_I$ be the multiplicity of $I$. If $I\in N_i^+$, let $s_I$ be the integer such that $s_{I}/r_I$ is the greatest height $h$ of any ancestor of $I$ with $1<h<2$, and set $R_I:=s_I-r_I < r_I$.
\begin{itemize}\itemsep=0pt
\item The generic class of representations has rank
\[r_{\rm gen}=\sum_{i=1}^k \Bigg( \sum_{I\in N_i} n_I r_I\Bigg).\]
\item The $i$-th nongeneric class of representations has rank
\[r_i=\sum_{I\in N_i^+} n_I R_I+\sum_{\substack{l=1\\l\neq i}}^k \Bigg( \sum_{I\in N_{l}} n_I r_I\Bigg).\]
\end{itemize}
\end{enumerate}
\end{Theorem}

The link with the diagrams of \cite{doucot2021diagrams} is as follows. Any irregular connection on $\mathbb P^1$ canonically determines both a diagram $\Gamma$ and an enriched tree $\mathscr T$, both invariant under symplectic transformations. Our main result says that the datum of the enriched tree allows us to construct~${k+1}$ different \textit{readings} of the diagram, corresponding to the classes of representations in an $\SL_2(\mathbb C)$-orbit. The enriched tree is thus a finer invariant that contains in general strictly more information than the diagram, i.e., we can view it as an enriched diagram. Concretely, the enriched tree determines explicitly the diagram, as well as a partition of the set of its core vertices into $k$ subsets.

Notice that, from the perspective of our general framework, the simply-laced case of \cite{boalch2012simply} corresponds to the situation where all vertices of the generic fission tree have integer heights (belonging to the set $\{1,2\}$). In this case, the datum of the short fission tree is equivalent to the datum of the complete $k$-partite core diagram, which is consistent with the main result of~\cite{boalch2012simply} saying that everything can be reconstructed from the supernova quiver.

Let us also mention that the fission trees of \cite{boalch2025twisted} are closely related to combinatorial objects appearing in the study of the singularities of plane algebraic curves: they can be seen as a~variant of the Eggers--Wall tree \cite{wall2004singular}. In some sense, we can thus view the short fission tree as a global analogue of Eggers--Wall trees, for the case of several singularities located at different points on~$\mathbb P^1$.\looseness=-1

We apply our framework to the case of nonabelian Hodge spaces related to Painlev\'e equations, which have complex dimension 2 and constitute the simplest nontrivial examples of nonabelian Hodge spaces. Concretely, for all Painlev\'e moduli spaces we describe their classes of basic representations, which are classes of nearby representations of the standard rank 2 Lax representations known since the early 20th century \cite{fuchs1907lineare,garnier1919classe,jimbo1981monodromyI}. In this way we recover many known alternative Lax representations as basic representations different from the standard ones, and we can view the remaining basic representations as providing new Lax representations (in our abstract sense).\looseness=-1

While we refer the reader to the last section of the article for many more pictures of fission trees and basic representations for Painlev\'e moduli spaces, a concrete example is already summarized on Figure \ref{fig:PIII_tree_and_diagram}, which corresponds to the classes of basic representations of the Painlev\'e~III moduli space which are the classes of nearby representations of the standard rank 2 Painlev\'e~III Lax representation (see Section~\ref{subsec:painleve_III} for a more detailed discussion of this example). The figure illustrates how, starting from the short fission tree, the other classes of representations can be obtained, yielding different readings of the diagram.

\begin{figure}[!ht]
\centering
\begin{tikzpicture}
\tikzstyle{mandatory}=[circle,fill=black,minimum size=5pt,draw, inner sep=0pt]
\tikzstyle{authorised}=[circle,fill=white,minimum size=5pt,draw,inner sep=0pt]
\tikzstyle{empty}=[circle,fill=black,minimum size=0pt,inner sep=0pt]
\tikzstyle{root}=[fill=black,minimum size=5pt,draw,inner sep=0pt]
\tikzstyle{indeterminate}=[circle,densely dotted,fill=white,minimum size=5pt,draw, inner sep=0pt]
\tikzstyle{vertex}=[circle,fill=cyan,minimum size=6pt,inner sep=0pt]
\tikzstyle{vertex_2}=[circle,fill=purple,minimum size=6pt,inner sep=0pt]

\begin{scope}[xshift=-5.4cm,yshift=-4cm]

\draw (1,2.5) node {\small 1st nongeneric rep.};
 %\draw (-0.5,1) node {\scriptsize $1$};
 \draw (-0.5,0.5) node {\scriptsize $1/2$};
 \node[root] (A1) at (0,1){};
 \node[root] (B1) at (1,1){};
 \node[root] (C1) at (2,1){};
 \node[empty] (A0) at (0,0){};
 \node[empty] (B0) at (1,0){};
 \node[empty] (C0) at (2,0){};
 \node[mandatory] (C1/2) at (2,0.5){};
 \draw (A1)--(A0);
 \draw (B1)--(B0);
 \draw (C1)--(C1/2)--(C0);

\node[vertex] (A) at (0,-3){};
\node[vertex_2] (C) at (2,-3){};
\node[vertex] (B) at (1,-4.5) {};
\draw[double distance=2 pt] (A)--(C)--(B);
\draw[dashed] (C) to [out=10, in=-45] (2.7,-2.3) to[out=135,in=80] (C);

\draw[->,dashed] (0,-1)--(0,-2);
\draw[->,dashed] (1,-1)--(1,-3.5);
\draw[->,dashed] (2,-1)--(2,-2);

\draw (1,-5) node {\small rank 2};
\draw (1,-5.5) node {\small 3 singularities};
\draw (1,-6) node {\small Katz ranks 0, 0, 1/2};
\end{scope}

\begin{scope}
\draw (1,3.5) node {\small Generic rep.};
 \draw (-0.5,1) node {\scriptsize $1$};
 \draw (-0.5,2) node {\scriptsize $2$};
 \draw (-0.5,0.5) node {\scriptsize $1/2$};
 \node[root] (A3) at (1,3){};
 \node[authorised] (A2) at (0.5,2){};
 \node[authorised] (B2) at (2,2){};
 \node[authorised] (A1) at (0,1){};
 \node[authorised] (B1) at (1,1){};
 \node[authorised] (C1) at (2,1){};
 \node[mandatory] (C1/2) at (2,0.5){};
 \node[empty] (A0) at (0,0){};
 \node[empty] (B0) at (1,0){};
 \node[empty] (C0) at (2,0){};
 \draw (A3)--(A2);
 \draw (A3)--(B2);
 \draw (A2)--(A1)--(A0);
 \draw (A2)--(B1)--(B0);
 \draw (B2)--(C1)--(C1/2)--(C0);

\node[vertex] (A) at (0,-7){};
\node[vertex_2] (C) at (2,-7){};
\node[vertex] (B) at (1,-8.5) {};
\draw[double distance=2 pt] (A)--(C)--(B);
\draw[dashed] (C) to [out=10, in=-45] (2.7,-6.3) to[out=135,in=80] (C);

\draw[->,dashed] (0,-1)--(0,-6);
\draw[->,dashed] (1,-1)--(1,-7.5);
\draw[->,dashed] (2,-1)--(2,-6);

\draw (1,-9) node {\small rank 4};
\draw (1,-9.5) node {\small 1 singularity};
\draw (1,-10) node {\small Katz rank 2};

\draw[line width=1pt, double distance=3pt,
 arrows = {-Latex[length=0pt 3 0,open]}] (-1,0.5) -- (-3,-1);
 \draw[line width=1pt, double distance=3pt,
 arrows = {-Latex[length=0pt 3 0,open]}] (3,0.5) -- (5,-1);
\draw (-3.3,0.3) node {\scriptsize transform 1st subtree};
\draw (5.3,0.3) node {\scriptsize transform 2nd subtree};

\end{scope}

\begin{scope}[xshift=5.4cm,yshift=-4cm]
\draw (1,2.5) node {\small 2nd nongeneric rep.};
\draw (-0.5,1) node {\scriptsize $1$};
 \draw (-0.5,2) node {\scriptsize $2$};

 \node[root] (A2) at (0.5,2){};
 \node[root] (B2) at (2,2){};
 \node[authorised] (A1) at (0,1){};
 \node[authorised] (B1) at (1,1){};
 \node[authorised] (C1) at (2,1){};
 \node[empty] (A0) at (0,0){};
 \node[empty] (B0) at (1,0){};
 \node[empty] (C0) at (2,0){};

 \draw (A2)--(A1)--(A0);
 \draw (A2)--(B1)--(B0);
 \draw (B2)--(C1)--(C0);

\node[vertex] (A) at (0,-3){};
\node[vertex_2] (C) at (2,-3){};
\node[vertex] (B) at (1,-4.5) {};
\draw[double distance=2 pt] (A)--(C)--(B);
\draw[dashed] (C) to [out=10, in=-45] (2.7,-2.3) to[out=135,in=80] (C);

\draw[->,dashed] (0,-1)--(0,-2);
\draw[->,dashed] (1,-1)--(1,-3.5);
\draw[->,dashed] (2,-1)--(2,-2);

\draw (1,-5) node {\small rank 2};
\draw (1,-5.5) node {\small 2 singularities};
\draw (1,-6) node {\small Katz ranks 1, 1};
\end{scope}

\end{tikzpicture}\vspace{-1mm}

\caption{An example: classes of basic representations of the Painlev\'e~III moduli space, giving 3 readings of the same diagram. The rational numbers on the left of the fission forests are the heights of the vertices. The generic representation corresponds to a short fission tree with $k=2$ principal subtrees. Each leaf corresponds to a vertex in the diagram, the subset $N_1=N_1^-$ corresponds to blue vertices, and $N_2=N_2^-$ corresponds to the red vertex. There are thus 2 nongeneric classes of nearby representations, whose fission forests are drawn on the figure. For each $i\in \{1,2\}$ the forest is obtained from the short fission tree by changing the $i$-th principal subtree into its Fourier transform, and keeping the other principal subtree unchanged (up to some changes of truncations of the trees). Notice that the nongeneric representation on the right corresponds to the standard Painlev\'e~III Lax representation, with two second order poles (Poincar\'e--Katz rank 1), while the one on the left corresponds to the alternative Lax representation for Painlev\'e sometimes referred to as `degenerate Painlev\'e~V', with one twisted irregular singularity (of Poincar\'e--Katz rank $1/2$), and two regular singularities. See Section~\ref{subsec:painleve_III}, and in particular Remark~\ref{remark:fourier_and_other_operations}, for further discussion of this example, including explicit expressions of irregular classes realizing these representations.}
\label{fig:PIII_tree_and_diagram}\vspace{-2mm}
\end{figure}

To conclude this introduction, let us briefly sketch the main ideas leading to these results. An essential observation, already made in \cite{boalch2012simply} in the simply-laced case, and in \cite{doucot2021diagrams} in the general case, is that one can associate to any Stokes circle $I$ an element $\lambda(I)\in\mathbb P^1$ that we will call its \textit{Fourier sphere coefficient}, such that the $\SL_2(\mathbb C)$ action on the set of all possible Stokes circles (induced by the action on connections via the stationary phase formula) behaves in a very nice way with respect to those: it acts by homographies on the \textit{Fourier sphere} $\mathbf P^1$, the copy of~$\mathbb P^1$ where those coefficients live. Concretely, for
$A=\left(\begin{smallmatrix}
a & b\\
c & d
\end{smallmatrix}\right)\in \SL_2(\IC)$, letting $h_A$ denote the corresponding homography
$h_A \colon \ z\mapsto \frac{az+b}{cz+d}$,
we have that $A$ acts on the Fourier sphere coefficients by $h_A$.

Then, for any connection $(E,\nabla)$ one can consider the set $\{\lambda_1,\dots,\lambda_k\}$ of the distinct Fourier sphere coefficients of its Stokes circles, which defines the integer $k$.

Using the results of \cite{boalch2025twisted} on admissible deformations, it is then quite straightforward to see that up to admissible deformations, the irregular class of $A\cdot (E,\nabla)$ only depends on whether we have $h_A(\lambda_i)=\infty$ for some $i\in\{1,\dots, k\}$, and if so for which $i$, which leads to the nearby representations: the generic class of representations corresponds to the case where $h_A(\lambda_i)\neq \infty$ for all $i\in\{1,\dots k\}$ (hence the terminology ``generic''), and the $i$-th nongeneric class of representations to the case where ${h_A(\lambda_i)=\infty}$. Each distinct Fourier sphere coefficient $\lambda_i$ corresponds to one principal subtree of the fission tree of the generic class of representations.

\subsection{Outlook}
An interesting question for future work could be to investigate to which extent one can devise a general procedure to construct, for any such abstract Lax representation, a corresponding fully explicit Lax representation, that is suitable explicit parametrisations of connections with the singularity data determined by the abstract Lax representation, such that their isomonodromic deformations imply that some function $y$ appearing in the parametrisation satisfies the corresponding Painlev\'e equation. Indeed, the search for alternative Lax representations for Painlev\'e-type equations has been an active topic of research \cite{bobrova2023different,joshi2007linearization,joshi2009linearization,ormerod2017symmetric,suleimanov2008quantizations}.

Furthermore, a natural extension of this work is the question of actually finding explicitly the isomorphisms expected to exist between the nonabelian Hodge spaces corresponding to the different representations, and of understanding how the spaces of isomonodromic times for different representations are related, so as to obtain, as in the simply-laced case, a full correspondence between the isomonodromy systems.

\subsection{Contents}
The text is organised as follows. In Section~\ref{section:preliminaries}, we review the necessary facts about singularity data of irregular connections, admissible deformations, fission forests, and diagrams, and describe how these data are transformed under symplectic transformations, relying on the results of~\cite{doucot2021diagrams} and slightly extending them to deal with admissible deformations. In Section~\ref{section:basic_reps}, we define the generic form of a Stokes circle and a modified irregular class, Fourier sphere coefficients, and establish our main results, in particular we provide explicit formulas for the rank of all nearby representations. Finally, in Section~\ref{section:painleve_examples}, we apply this framework to describe in detail the basic representations of all Painlev\'e moduli spaces.

\section[Symplectic action on fission data of irregular connections on P\^{}1]{\!Symplectic action on fission data of irregular connections on~$\boldsymbol{\mathbb P^1}$}
\label{section:preliminaries}

In this section, we briefly review the necessary facts about formal data of connections with irregular singularities on the Riemann sphere, their transformation under Fourier transform and more general symplectic transformations, admissible deformations, and combinatorial objects (fission forests, diagrams) encoding in a useful way parts of the information contained in the formal data.

\subsection{Formal data of irregular connections}

Let us first recall the description of formal data of irregular connections on the Riemann sphere.

The well-known Turritin--Levelt theorem states that any connection on the formal punctured disc can be decomposed as a direct sum of elementary connections having only one exponential factor. This can be formulated in different ways; the one we will be using is the point of view \`a la Deligne--Malgrange (see, e.g., \cite{malgrange1982classification}), where one views exponential factors at a singularity as sections of a local system over the circle of directions around this singularity (see, e.g., \cite{boalch2021topology,doucot2021diagrams} for more details).

\subsubsection{The exponential local system}

Let $\Sigma=\mathbb P^1$ and $a\in \Sigma$. Let $\pi_a\colon \widehat{\Sigma}_a\to \Sigma$ be the real oriented blow-up at $a$ of $\Sigma$. The preimage~${\partial_a:=\pi^{-1}(a)}$ is a circle whose points correspond to the directions around $a$. An open subset of~$\partial_a$ corresponds to a germ of sector at $a$.

Let $z_a$ be a local complex coordinate on $\Sigma$ in a neighbourhood of $a$, vanishing at $a$. The exponential local system $\Ical_a$ is a local system of sets (i.e., a covering space) over~$\partial_a$ whose sections are germs of holomorphic functions on sectors of the form
$
\sum_{i} b_i z_a^{-k_i}$,
where $k_i\in \mathbb{Q}_{>0}$, and $b_i\in\IC$.
The connected component of such a local section is a finite order cover of the circle $\partial$. More precisely, let $r$ be the smallest integer such that the expression~${q=\sum_{i} b_i z_a^{-k_i}}$ is a polynomial in \smash{$z_a^{-1/r}$}. The corresponding holomorphic function is multivalued, and becomes single-valued when passing to a finite cover $t_a^r=z_a$. Therefore, the corresponding connected component, which we denote by $\cir{q}_a$, is a $r$-sheeted cover of $\partial_a$. As a~topological space, it is homeomorphic to a circle, and $\Ical_a$ is thus a disjoint union of (an infinite number of) such \textit{Stokes circles}.

The integer $r$ is the \textit{ramification order} of the Stokes circle \smash{$I=\cir{q}_a$}, which we denote by $r=:\ram(q)=\ram(I)$. If $r=1$, we say that the Stokes circle $\cir{q}_a$ is untwisted, or unramified. The degree $s$ of $q$ as a polynomial is the \textit{irregularity} $\irr(q)=\irr(I)$ of $q$ (or $I$), and the quotient~$s/r$ is the \textit{slope} of $q$ (or $I$), which we denote by $\slope(q)=\slope(I)$.

If $d\in\partial_a$, we denote by $(\Ical_a)_d$ the fibre of $\Ical_a$ over the direction $d$. Taking $d$ as a basepoint of $\partial_a$, the monodromy of $\Ical_a$ is the automorphism
$\rho_a\colon (\Ical_a)_d \to (\Ical_a)_d
$
of the fibre $(\Ical_a)_d$ obtained when going once around $\partial_a$ in the positive direction.

Finally, we define the global exponential local system as the disjoint union
\[
\mathcal I:=\bigsqcup_{a\in \Sigma} \mathcal I_a.
\]
It naturally comes with a projection $\pi\colon\mathcal I\to \partial$, where
\[
\partial:=\bigsqcup_{a\in \Sigma} \partial_a.
\]
If $I\subset \mathcal I_a$ is a Stokes circle at $a$ for some $a\in \mathbb P^1$, by a slight abuse of notation we will write~${\pi(I):=a}$.

\subsubsection{Local and global formal data}

\begin{Definition}
A (local) \emph{irregular class} at $a$ is a function $\Theta_a\colon \pi_0(\Ical_a)\to \mathbb N$ with compact support. The rank of an irregular class $\Theta_a$ is
\[
\rank(\Theta_a)=\sum_{I\in \pi_0(\Ical_a)} n_I \ram(I),
\]
where $n_I=\Theta_a(I)$ is the multiplicity of the circle $I$. If $n_I\neq 0$, we say that $I$ is an active circle of $\Theta_a$, or a Stokes circle of $\Theta_a$.

Furthermore, the \textit{Poincar\'e--Katz rank} of $\Theta_a$ is
\[
\operatorname{Katz}(\Theta_a):=\max\left\{ \slope(I) \mid n_I>0\right\}.
\]
\end{Definition}

The Turritin--Levelt theorem implies that any rank $r$ connection on the formal punctured disc~$\operatorname{Spec}(\mathbb C((z_a)))$ at $a$ uniquely determines a pair
$(\Theta_a, \mathcal C_a)$, where $\Theta_a$ is a rank $r$ irregular class, and~$\Ccal_a$ is the datum of a conjugacy class $\Ccal_{I}\in GL_{n_I}(\mathbb C)$ for any active circle~$I$ of $\Theta_a$ with multiplicity~$n_I$. The pair $(\Theta_a,\Ccal_a)$ constitutes the formal data of the connection and determines it completely up to isomorphism.

Passing to the global situation, to any algebraic connection $(E,\nabla)$ on a Zariski open subset~${\Sigma\smallsetminus \mathbf a}$ of $\Sigma$, by restricting the connection on a small disc around each singularity and considering its formalization there, one associates \emph{global formal data} $(\bm\Theta,\bm\Ccal):=(\Theta_a, \Ccal_a)_{a\in \mathbf a}$.

The Stokes circles of $\Theta_a$ at each singularity $a\in \Sigma$ correspond to the exponential factors of the connection in the Turritin--Levelt normal form at $a$, and the conjugacy classes $\Ccal_a$ encode its formal monodromy. In particular, the connection has a regular singularity at $a$ if and only if its only Stokes circle at $a$ is the tame circle $\cir{0}_a$.

\begin{Remark}
Notice that if we wanted to work with the full nonabelian Hodge package, we would need to consider weighted conjugacy classes (cf.\ \cite{boalch2011riemann}). Here for simplicity, since we will not be trying in this work to prove that the different classes of representations correspond as one expects to isomorphisms between the full nonabelian Hodge spaces, we will just work with standard (unweighted) conjugacy classes.
\end{Remark}

\subsubsection{Modified formal data} It turns out that the action of the Fourier transform on formal data of connections is more easily described in terms of \textit{modified formal data}, obtained as follows: if $i$ is a tame circle at finite distance (that is $i=\cir{0}_a$ with $a=\pi(i)\neq \infty$), then replace $n_i$ by the integer
$m_i := \text{rank}(A-1)$,
where $A\in \Ccal_i$.
Thus $A=1+vu$ for a surjective map $u\colon\IC^{n_i}\to \IC^{m_i}$ and an injective map~${v\colon\IC^{m_i}\to \IC^{n_i}}$.
Then replace $\Ccal_i$ by the conjugacy class $\breve\Ccal_i$ of $(1+uv)$ in $\GL_{m_i}(\IC)$;
this class is called the {\em child} of $\Ccal_i$ in \cite[Appendix]{boalch2016global}.
If $\pi(i)=\infty$ or if $i$ is not tame (\smash{$i\neq \cir{0}_a$} for some $a$),
then we do no modification: $m_i:=n_i$ and
$\breve\Ccal_i:=\Ccal_i$.
This defines modified formal data $\bigl(\bm{\breve\Theta}, \bm{\breve\Ccal}\bigr)$, where
\smash{$\bm{\breve\Theta}$} is the collection of all Stokes circles with multiplicities $m_i\ge 1$.

Notice that $\bm{\breve\Theta}$ satisfies the condition $\rank\bigl(\breve\Theta_a\bigr)\leq \rank\bigl(\breve\Theta_\infty\bigr)$ since for each $a\neq \infty$ passing to the modified formal data lowers the rank. If \smash{$\bm{\breve\Theta}$} is any global irregular class, we say it is \textit{compatible} if it satisfies this condition. If there exists an irreducible connection $(E,\nabla)$ on a~Zariski open subset of the affine line with modified formal data $\bigl(\bm{\breve\Theta}, \bm{\breve\Ccal}\bigr)$, we say that $\bigl(\bm{\breve\Theta},\bm{\breve\Ccal}\bigr)$ is \textit{effective}. In particular, \smash{$\bigl(\bm{\breve\Theta}, \bm{\breve\Ccal}\bigr)$} cannot be effective if \smash{$\bm{\breve\Theta}$} is not compatible.

Furthermore, if we know the modified formal data $\bigl(\bm{\breve\Theta}, \bm{\breve\Ccal}\bigr)$ of a connection $(E,\nabla)$, since we keep the non-modified formal data at infinity, we know the rank of $(E,\nabla)$: we have $\rank(\Theta_\infty)=\rank\bigl(\breve\Theta_\infty\bigr)$. In turn, this allows us to reconstruct from $\bigl(\bm{\breve\Theta}, \bm{\breve\Ccal}\bigr)$ the non-modified formal data $(\bm\Theta, \bm\Ccal)$, so that the datum of \smash{$\bigl(\bm{\breve\Theta},\bm{\breve\Ccal}\bigr)$} is thus equivalent to the one of $(\bm\Theta,\bm{\Ccal})$, cf.\ \cite[Section~5.1]{arinkin2010rigid}. We say that $(\bm\Theta,\bm{\Ccal})$ is effective if $\bigl(\bm{\breve\Theta},\bm{\breve\Ccal}\bigr)$ is effective.

\subsection{Admissible deformations, fission data, and forests}
\label{subsec:admissible_defs_and_forests}

We now recall a few definitions and results about admissible deformations of irregular connections, involving the notions of levels and fission trees, and about diagrams. We refer the reader to \cite{boalch2025twisted, doucot2021diagrams} for more details and for proofs of the statements.

\subsubsection{Local version} First, to define the notion of admissible deformation of an irregular class, we need to recall the notion of (full) irregular type.

\begin{Definition}
A \emph{full irregular type} at $a\in \Sigma$ of rank $n$ is a Galois-closed list $T=(q_1,\dots, q_n)$ of exponential factors in some direction $d\in \partial$ (i.e., $q_i$ is a section of $\mathcal I_a$ in a germ of sector around $d$).
\end{Definition}

Any full irregular type $T=(q_1,\dots, q_n)$ defines an associated irregular class $\Theta$ in the following way: its active circles are the connected components $\cir{q_i}$ in $\mathcal{I}$ of the exponential factors, and the multiplicity $n_i$ of $\cir{q_i}$ is the multiplicity of $q_i$ in the list $(q_1,\dots, q_n)$. Passing from a full irregular type to an irregular class amounts to forgetting the ordering of the exponential factors.

\begin{Definition}
Let $T=(q_1,\dots, q_n)$ and $T'=(q'_1,\dots, q'_n)$ be rank $n$ full irregular types at $a$ in some direction $d$. We say that $T$ and $T'$ are admissible deformations of each other, and denote this by $T\sim T'$, if we have
$
\slope(q_i-q_j)=\slope(q'_i-q'_j)$,
for all $i,j\in\{1,\dots, n\}$.

If $\Theta$ and $\Theta'$ are two rank $n$ irregular classes at $a\in \Sigma$, we say that they are admissible deformations of each other, and write $\Theta\sim \Theta'$, if there exist two full irregular types $T$ and $T'$ with respective associated irregular classes $\Theta$ and $\Theta'$ such that $T\sim T'$.
\end{Definition}

In particular, if $\Theta\sim\Theta'$, denoting by $N$ and $N'$ the sets of active circles of $\Theta$ and~$\Theta'$, respectively, there exists a bijection $\phi\colon N\to N'$ such that for each active circle $I\in N$, $\phi(I)$ is an admissible deformation of $I$, and $I$, $I'$ have the same multiplicity. We say that such a bijection is \textit{compatible}.

\begin{Definition}
Let $(\Theta, \mathcal C)$ and $(\Theta', \mathcal C')$ be two sets of local formal data. We say that they are admissible deformations of each other if $\Theta\sim\Theta'$ and there exists a compatible bijection~${\phi\colon N\to N'}$ between their respective sets of active circles $N$, $N'$ such that \smash{$\Ccal'_{\phi(I)}=\Ccal_I$} for all~${I\in N}$.
\end{Definition}

It turns out that there exist simple numerical criteria characterising whether two irregular classes are admissible deformations of each other.

\begin{Definition}

Let $\Theta$ be an irregular class at $a\in \Sigma$. The levels of $\Theta$ are the elements of the set
\[
\mathrm{Levels}(\Theta)=\{ \slope(q_\alpha-q_\beta)\mid \alpha,\beta\in \mathtt I_d\}\setminus \{0\}\subset \mathbb Q_{>0},
\]
where $\mathtt I\subset \mathcal I$ is the finite subcover underlying $\Theta$ and $\mathtt I_d$ is any fibre of $\mathtt I$.
\end{Definition}

It is clear that two irregular classes which are admissible deformations of each other have the same levels.

Let us first consider the case of an irregular class having only one Stokes circle. Let $I$ be a~Stokes circle at $a$. We can always write $I=\cir{q}_a$, with $q=\sum_{i=1}^p a_i z_a^{-k_i}$, for some integer $p$, with~${k_1>\dots >k_p}$ and $a_i\neq 0$ for every $i\in\{1,\dots, p\}$. The rational numbers $k_i$ are uniquely defined by $I$, and we call them the \textit{exponents} of $I$. We have the following.

\begin{Proposition}[{\cite[Proposition~3.1]{boalch2025twisted}}]~
\begin{itemize}\itemsep=0pt
\item[$(a)$] Two Stokes circles $I,J\subset \mathcal I_a$ are admissible deformations of each other if and only if $\Levels(I)=\Levels(J)\subset \mathbb Q$.

\item[$(b)$] A subset $(k_1>k_2>\cdots >k_m)\subset \mathbb Q_{>0}$ is the set of levels of some circle $I\subset \mathcal I_a$ if and only if
%\label{eq: condition on levels of a circle}
$k_1,k_2,\ldots,k_m$ have strictly increasing common denominators $>1$.
%\end{equation}
In other words if $d_i$ is the denominator of $k_i$ {\rm(}in lowest terms$)$ and
$
r_i$ is the lowest common multiple of~${d_1,d_2,\ldots,d_i}$
for each $i$ {\rm(}so that $r_i\mid r_{i+1})$, then $1<r_1<r_2<\cdots <r_m$.

\item[$(c)$] Let $I$ be a Stokes circle, and let $(k_1>\dots >k_p)$ denote its set of exponents. Then~$\Levels(I)$ is the largest possible subset of $\{k_1,\dots,k_p\}$ of the form $k_{i_1}>\dots >k_{i_m}$ for with $1=i_1<\dots <i_m\leq p$ such that $k_{i_1},\dots,k_{i_m}$ have strictly increasing common denominators~$>1$.
\end{itemize}
\end{Proposition}

Passing to the general case with irregular classes having several Stokes circles, to characterise admissible deformations, we need more data than the levels of each individual Stokes circle. To this end, we introduce the notion of common part and fission exponent of a pair of Stokes circles.

If \smash{$q=\sum_i b_i z_a^{-i/r}$} is an exponential factor at $a$, and $k\in \mathbb Q_{\geq 0}$, we write
\[
q_{\geq k}:=\sum_{i,i/r\geq k} b_i z^{-i/r}.
\]
In turn, if $I=\cir{q}$ is a Stokes circle, we set $I_{\geq k}:=\cir{q_{\geq k}}$ (this is well-defined since it does not depend on the choice of representative). We define in a similar way truncations $I_{>k}$, $I_{\leq k}$, $I_{<k}$.

\begin{Definition}
Let $I$, $\wh I$ be two Stokes circles at the same point with respective ramification orders $r$, $\wh r$. The \emph{common part} of $I$ and $\wh I$ is the circle $I_c=I_{\geq k}=\wh I_{\geq k}$, where
\[
k=\min\left( l\in \frac{1}{r\wh r}\mathbb Z_{\geq 0} \,\middle\vert\, I_{\geq l}=\wh I_{\geq l}\right).
\]
The \textit{fission exponent} of $I$ and $\widehat I$ is the rational number
\[
f_{I,\wh I}=\max\bigl(\slope(I_{<k}), \slope\bigl(\wh I_{<k}\bigr)\bigr)\in \frac{1}{r\wh r}\mathbb Q.
\]
\end{Definition}

If $I=\wh I$, then $I_c=I=\wh I$ and \smash{$f_{I,\wh I}=0$}, otherwise \smash{$f_{I,\wh I}>0$}. $I_c=0$ if and only if the leading terms of $I$ and $\wh I$ have different Galois orbits, in this case we say that $I$ and $\wh I$ have no common part.

\begin{Lemma}[{\cite[Lemma 3.11]{boalch2025twisted}}]
\label{lemma:slopes_differences}
Let $I$ and $\wh I$ be two Stokes circles with ramification orders~$r$,~$\wh{r}$. The set of slopes among the rational numbers $\slope\bigl({q_i}-{\wh{q}_j}\bigr)$, for $i\in\{0,\dots, r-1\}$, $j\in\big\{0,\dots, \allowbreak {\wh{r}-1}\big\}$ is of the form
\smash{$
\Levels(I_c) \sqcup \{f_{I,\wh I}\} \subset \mathbb Q_{\geq 0}$},
i.e., it consists of the levels of the common part of the two Stokes circles together with their fission exponent.
\end{Lemma}

This motivates the following definition.

\begin{Definition}
A (local) fission datum is a pair $\mathcal F=(\mathcal L,f)$
where $\mathcal L$ is a multiset $\mathcal L=L_1+\cdots + L_n$ of level data
and $f$, the fission exponents, consist of
the choice of a rational number
$f_{ij}=f_{ji}\in \mathbb Q_{\ge 0}$, for all
 $i,j\in \{1,\ldots,n\}$.
\end{Definition}

A (local) irregular class determines a fission datum as follows.
If $\Theta=\sum_1^n I_i$ is a rank $n$ irregular class at some point $a$ (where the Stokes circles $I_i$ are not necessarily distinct), then define
$L_i=\mathrm{Levels}(I_i)$ to be the level datum of $I_i$ for each $i$,
and define
$\mathcal L(\Theta) := \sum_1^n L_i$
to be the corresponding multiset of level data.
Then by taking $f_{ij}=f_{I_i,I_j}$ to be corresponding fission exponents, this determines the fission datum
$\mathcal F(\Theta)=(\mathcal L(\Theta),f)$
of the irregular class $\Theta$.
Note that the multiplicity of any
given Stokes circle $I_j$ in the class
$\Theta$ is determined from the fission data by the recipe
$\Theta(I_j) = \abs{\{i\in\{1,2,\ldots,n\} \mid f_{ij}=0\}}$.

The fission datum characterises classes of admissible deformations of (local) irregular classes.

\begin{Theorem}[{\cite[Theorem~3.17]{boalch2025twisted}}]
Two $($local$)$ irregular classes $\Theta$ and $\Theta'$ are admissible deformations of each other if and only if they have the same fission datum.
\end{Theorem}

Notice that not any fission datum $\mathcal F$ actually comes from an irregular class $\Theta$. To characterise such fission data, as well as provide a more graphical way of encoding them, the notion of (twisted) \textit{fission tree} was introduced in \cite{boalch2025twisted}. We will not be needing to repeat all the exact axioms for fission trees here, so we just briefly recall the main properties of the fission tree associated to an irregular class, referring the reader once more to loc.\ cit.\ for full details, and to Section~\ref{section:painleve_examples} for many concrete examples of fission trees.

{\samepage If $\Theta$ is an irregular class at infinity, its fission tree $\mathcal T(\Theta)$ is a tuple $(\mathcal T,\mathbb V,\mathbb A,\mathbb L, h,n)$, where
\begin{itemize}\itemsep=0pt
\item $\mathcal T$ is a metrized tree, with set of \textit{vertices} $\mathbb V$. An \textit{edge} of the tree is an element of $E:=\pi_0(\mathcal T \smallsetminus \mathbb V)$ of $\mathcal T$. Any vertex that is not a leaf nor the root is adjacent to $\geq 2$ edges, one of which is the parent edge, and the others are the descendant edges. The branch vertices~${\mathbb Y\subset \mathbb V}$ are those with $\geq 2$ descendants. The trunk of the
tree is the union of all edges and vertices above all the branch vertices.

If $i$ is a leaf of $\mathcal T$, we define the \emph{full branch} $\mathcal B_i$ of $i$ as the subset of $\mathcal T$ corresponding to the path from $i$ all the way to the end of the trunk.

\item $\mathbb L$ and $\mathbb A$ are two subsets of $\mathbb V$, whose elements are respectively called \textit{mandatory} and \textit{admissible} vertices, and $\mathbb L\subset \mathbb A$. The set $\mathbb L$ is finite and can be empty.

\item $h\colon\mathcal T\!\to\! \mathbb R_{\geq 0}$ is a length-preserving function, which we call the height, such that ${h^{-1}(0)\!=:\!\mathbb V_0}$ is the set of leaves of $\mathcal T$, $h$ maps isomorphically any edge to an interval, and for any leaf $i$, $h$ maps the full branch $\mathcal B_i$ isomorphically onto $\mathbb R_{\geq 0}$.
\item $n$ is a map $\mathbb V_0\to \mathbb Z_{\geq 1}$, i.e., consists of the datum of an integer multiplicity $n_i\geq 1$ for each leaf $i\in \mathbb V_0$.
\end{itemize}

}

Moreover, the fission tree $\mathcal T(\Theta)$ has the following properties:
\begin{itemize}\itemsep=0pt
\item The leaves are in one-to-one correspondence with the Stokes circles of $\Theta$: if $\Theta=\sum_{i\in N}n_i I_i$, we have $\mathbb V_0=N$, with the leaf $i$ corresponding to the Stokes circle $I_i$, and having multiplicity $n_i$.
\item If $i\in \mathbb V_0$ is a leaf and $I$ is the corresponding Stokes circle, the heights of the mandatory vertices of the full branch $\mathcal B_i$ are exactly the levels of $I$, that is
$
h(\mathbb L\cap \mathcal B_i)=\Levels(I)$.

\item If $i,j\in \mathbb V_0$ are two distinct leaves, $I$, $J$ are the corresponding Stokes circles, and $v_{ij}\in \mathbb V$ is closest common ancestor of $i$, $j$ then the height of all children of $v_{ij}$ is equal to the fission exponent $f_{I,J}$.

\end{itemize}

In brief, the leaves of the tree correspond to the Stokes circles of $\Theta$, and the tree contains the information about the level data as well as the fission exponents. In particular, the datum of the fission tree $\mathcal T(\Theta)$ is equivalent to the datum of $\mathcal F(\Theta)$.

Finally, if we include conjugacy classes of formal monodromies, a set of local formal data~$(\Theta,\mathcal C)$ determines a fission tree $\mathbf T:=\mathcal T(\Theta)$, and $\mathcal C$ amounts to the datum of a conjugacy class in~$\mathrm{GL}_n(\mathbb C)$ for each leaf of $\mathbf T$, where $n$ is its multiplicity. In turn, a class of admissible deformations of local formal data is characterised by such pairs $(\mathbf T, \mathcal C)$.

\subsubsection{Global version} Until now, in this paragraph we have dealt only with the local situation. The global case is described in a similar way, by collecting together the local data at each singularity.

For $\mathbf a\subset \Sigma$ a finite subset, and
$\mathbf{\Theta}=\{\Theta_a\mid a\in \mathbf a\}$ the datum of an irregular class for each marked point, we define the global fission datum of~${\bm\Theta}$ as the multiset of its local fission data~$\smash{\bm{\mathcal{F}}(\bm\Theta)=\sum_{a\in\mathbf a}\mathcal{F}(\Theta_a)}$. Similarly, we define the forest of $\bm\Theta$
as the multiset $\mathbf{F}(\bm\Theta)=\sum_{a\in \mathbf a} [\mathcal T(\Theta_a)]$ of isomorphism classes of
fission trees determined by all the $\Theta_a$,
as $a$ ranges over the marked points $\mathbf a \subset \Sigma$.
In general, it is a multiset rather than a set as some of the fission trees at distinct points may be isomorphic.
The notion of admissible deformations naturally generalises to the global case for global irregular classes $\bm{\Theta}$ (or modified ones), as well as for wild Riemann surfaces with boundary data $(\bm\Theta, \bm{\Ccal})$ and we have the following.

\begin{Corollary}[{\cite[Corollary~3.33]{boalch2025twisted}}]
 Two global irregular classes on $\Sigma$ are admissible deformations of each other if and only if they have the same global fission datum, or equivalently the same fission forest.
\end{Corollary}

As a consequence, a class of admissible deformations of wild Riemann surfaces with boundary data is characterised by a pair $(\mathbf F, \bm{\Ccal})$, where $\mathbf F$ is a fission forest, and $\bm{\Ccal}$ the datum of a conjugacy class in $\GL_n(\mathbb C)$ for each leaf of $\mathbf F$, where $n$ is the leaf multiplicity.

\subsection{Diagrams}

Another invariant that we can attach to any connection $(E,\nabla)$ on a Zariski open subset of $\mathbb P^1$ is a \textit{diagram} $\Gamma(E,\nabla)$ \cite{doucot2021diagrams}. As for the forest, the diagram only depends on the (modified) formal data $\bigl(\bm{\breve\Theta},\bm{\breve\Ccal}\bigr)$ of the connection.

The diagram \smash{$\Gamma(E,\nabla)=\Gamma\bigl(\bm{\breve\Theta},\bm{\breve\Ccal}\bigr)=(N,B)$} of the connection $(E,\nabla)$
has the following structure: it consists of a \textit{core diagram} $\Gamma_c(E,\nabla)=\Gamma_c\bigl(\bm{\breve\Theta}\bigr)$ with nodes $N$, to which is then glued, onto each core node, a \emph{leg} determined by a conjugacy class.

The set of core nodes $N$ is the set of Stokes circles of the modified irregular class $\bm{\breve{\Theta}}$. In terms of the non-modified irregular class $\bm\Theta$, $N$ is the disjoint union of three finite sets
\begin{align*}
N={}&\{\text{Stokes circles of $(E,\nabla)$ at $\infty$}\} \\
&\cup\{\text{wild Stokes circles of } (E,\nabla) \text{ not at $\infty$}\} \\
&\cup\{\text{tame Stokes circles of } (E,\nabla)\text{ not at $\infty$ with
nontrivial formal monodromy}\}.
\end{align*}

The matrix of edge/loop multiplicities
between the core nodes
is then defined as follows.
For any Stokes circle $I=\cir{q}$,
let $\alpha_I=\text{Irr}(I), \beta_I=\ram(I)\in \mathbb N$
denote its irregularity and
ramification numbers,
so that $I$ has slope $\alpha_I/\beta_I$.
We define an integer $B^\infty_{IJ}$ as follows:
\[%\label{eq: BY1}\label{eq: BY2}
B^\infty_{IJ}:=A_{IJ}-\beta_I\beta_J\quad \text{if}\ I\neq J,\qquad
B^\infty_{II}:=A_{II}-\beta_I^2+1\quad \text{otherwise,}
\]
where $A_{IJ}:=\irr(\Hom(I,J))$
is the irregularity of the irregular
class $\Hom(I,J)$.

\begin{Definition}\label{def:def_diagram_introduction}
Suppose $I,J\in N$ are Stokes circles, and write
$a_I=\pi(I),a_J=\pi(J)\in \mathbb P^1$ for the underlying
points of the Riemann sphere.
Define $B_{IJ}\in \mathbb Z$ as follows:
\begin{enumerate}\itemsep=0pt
\item[(1)] If $a_I=a_J=\infty$, then $B_{IJ}:=B^{\infty}_{IJ}$,
\item[(2)] If $a_I=a_J\neq \infty$, then
$B_{IJ}=B_{JI}:=B^{\infty}_{IJ}-\alpha_I\beta_J-\alpha_J\beta_I$,
\item[(3)] If $a_I\neq \infty$ and $a_J=\infty$,
then $B_{IJ}=B_{JI}:=\beta_J(\alpha_I+\beta_I)$,
\item[(4)] If $a_I\neq \infty, a_J\neq \infty$ and
$a_I\neq a_J$, then $B_{IJ}=B_{JI}:=0$.
\end{enumerate}
\end{Definition}

This completes the definition of the core diagram.
To obtain the full diagram, one glues to each core node $I\in N$ a \textit{leg} defined by the conjugacy class $\breve{\mathcal C}_i$.

\begin{Remark}
The core diagram $\Gamma_c\bigl(\bm{\breve\Theta}\bigr)$ is entirely determined by the fission data of $\bm{\breve\Theta}$ together with the slopes of all Stokes circles at finite distance, because these data determine the numbers~${A_{IJ}=\irr(\Hom(I,J))}$. In general, the diagram contains strictly less information, because~$A_{IJ}$ is a sum of different contributions, coming from the levels of the common parts of~$I$,~$J$ and the fission exponent~$f_{I,J}$. As we shall review below, the interest of the (core) diagram is that, unlike the fission data, it is invariant under symplectic transformations.
\end{Remark}

\subsection{Fourier transform and symplectic transformations}

We now describe how the $\SL_2(\mathbb C)$ transformations on irregular connections on $\mathbb P^1$ act on their formal data and fission data.

\subsubsection{Fourier transform}

The Fourier transform is the automorphism of the Weyl algebra $A_1=\IC[z]\cir{\partial_z}$ defined by $z\mapsto -\partial_z$, $\partial_z\mapsto z$. It induces a self-equivalence on the category of modules on the Weyl algebra: the Fourier--Laplace transform $F\cdot M$ of a $\IC[z]\langle \partial_z\rangle$-module $M$ is the $\IC[z]\langle \partial_z\rangle$-module obtained from $M$ by setting: $z\mapsto -\partial_z$, $\partial_z\mapsto z$. The Fourier transform induces a transformation on the set of isomorphism classes of irreducible connections \big(if we exclude the case of the trivial rank one connection on $\IP^1$\big), see \cite[Section~2.2]{arinkin2010rigid}.

The stationary phase formula \cite{fang2009calculation,sabbah2008explicit} states that the modified formal data of the Fourier transform $(E',\nabla')=F\cdot (E,\nabla)$ are determined by the modified formal data of $(E,\nabla)$.

\begin{Theorem}\label{th_stationary_phase} There exists a bijection, that we also denote by $F$, from the set of all modified formal data to itself, such that if $(E,\nabla)$ is a irreducible algebraic connection on a Zariski open subset of $\IP^1$ which is not a rank one connection with a single singularity at infinity of Poincar\'e--Katz rank $\leq 1$, $\bigl(\bm{\breve\Theta},\bm{\breve\Ccal}\bigr)$ its modified formal data, $(E',\nabla')$ is the Fourier transform of $(E,\nabla)$ and~\smash{$\bigl(\bm{\breve\Theta'},\bm{\breve \Ccal}'\bigr)$} its modified formal data, the following diagram commutes:
\[
\begin{tikzcd}
(E,\nabla) \arrow[r,"F",mapsto] \arrow[d,mapsto] & (E',\nabla') \arrow[d,mapsto]\\
\bigl(\bm{\breve\Theta},\bm{\breve\Ccal}\bigr) \arrow[r,"F",mapsto] & \bigl(\bm{\breve\Theta'},\bm{\breve \Ccal}'\bigr)
\end{tikzcd}
\]
\end{Theorem}

This map $F$ is the \emph{formal Fourier transform}. It acts independently on the Stokes circles: if we use the notation $\bigl(\bm{\breve\Theta},\bm{\breve\Ccal}\bigr)=\sum_{I\in N}\bigl(I,\breve{\mathcal C}_I\bigr)$ where $N$ is the set of active circles of $\bm{\breve\Theta}$, we have
\[
F\cdot\sum_{I\in N}\bigl(I,\breve{\mathcal C}_I\bigr)=\sum_{I\in N} F\cdot\bigl(I,\breve{\mathcal C}_I\bigr).
\]
 Furthermore, the irregular class $\bm{\breve{\Theta}}'$ only depends on $\bm{\breve{\Theta}}$: there is a self-bijection of the set of all Stokes circles $\pi_0(\mathcal I)$, that we will denote also by $F$, such that, if $\bm{\breve\Theta}=\sum n_i I_i$, then
\[
\bm{\breve\Theta}'=F\cdot \bm{\breve\Theta}=\sum_i n_i (F\cdot I_i).
\]

This bijection is actually given by a \textit{Legendre transform}. It can be computed explicitly, and takes slightly different forms depending on whether the Stokes circles are at finite distance or at infinity, and in the latter case on whether its slope is greater or less than~1, see \cite[Section~3.1]{doucot2021diagrams}. The next proposition describes the ramification orders, slopes and levels of the formal Fourier transform $F\cdot I$ of any Stokes circle $I$ as a function of those of~$I$, as well as the transformation of the conjugacy classes of formal monodromies.

\begin{Proposition}\label{prop:slopes_Fourier_transform}\samepage
Let $I$ be a Stokes circle, and $\breve{\mathcal C}_I\in \GL_n(\mathbb C)$ a conjugacy class. Let $\breve{\mathcal C}_{F\cdot I}\in \GL_n(\mathbb C)$ be the conjugacy class such that $F\cdot\bigl(I,\breve{\mathcal C}_I\bigr)=\bigl(F\cdot I, \breve{\mathcal C}_{F\cdot I}\bigr)$. We have the following:
\begin{enumerate}\itemsep=0pt
\item[$(1)$] If $I$ is of the form \smash{$\cir{\alpha z}_\infty$}, with $\alpha\in\IC$, then \smash{$F\cdot I=\cir{0}_{\alpha}$}, and \smash{$\breve{\mathcal C}_{F\cdot I}=\breve{\mathcal C}_I^{-1}$}.
\item[$(2)$] If $I$ is of slope $\leq 1$ at infinity, of the form \smash{$\cir{\alpha z +q}_\infty$}, with $\alpha\in\IC$, and $q\neq 0$ of slope~\mbox{$<1$}, then~${F\cdot I}$ is of the form \smash{$F\cdot I=\cir{\widetilde{q}}_{\alpha}$}, with $\irr(\tilde q)=s, \ram(\tilde q)=r-s$, $\Levels(\tilde q)=\frac{r}{r-s}\Levels(q)$. Also \smash{$\breve{\mathcal C}_{F\cdot I}=(-1)^s\breve{\mathcal C}_I^{-1}$}.
\item[$(3)$] If $I$ is of slope $> 1$ at infinity, with $\ram(I)=r$, $\irr(I)=s$, then $F\cdot I$ satisfies $\ram(F\cdot I)=s-r$, $\irr(F\cdot I)=s$, and $\Levels(F\cdot I)=\frac{r}{s-r}\Levels(I)$. Also $\breve{\mathcal C}_{F\cdot I}=(-1)^s\breve{\mathcal C}_I$.
\item[$(4)$] If $I=\cir{q}_a$ is an irregular circle at finite distance for $a\in\IC=\IP^1\smallsetminus\{\infty\}$, with ${q\neq 0}$, $\ram(q)=r$, and $\irr(q)=s$, then $F\cdot I$ is a circle of slope $\leq 1$, at infinity, of the form ${\cir{-az+\tilde{q}}_\infty}$, with $\ram(\tilde q)=r+s, \irr(\tilde q)=s$, and $\Levels(\tilde{q})=\frac{r}{r+s}\Levels(q)$. Also \smash{$\breve{\mathcal C}_{F\cdot I}=(-1)^s\breve{\mathcal C}_I^{-1}$}.
\item[$(5)$] If \smash{$I=\cir{0}_a$} is a tame circle at finite distance, then \smash{$F\cdot I=\cir{-az}_\infty$}, and \smash{$\breve{\mathcal C}_{F\cdot I}=\breve{\mathcal C}_I^{-1}$}.
\end{enumerate}
\end{Proposition}

\begin{proof}
All statements except those describing the sets of levels are already contained in the works \cite{malgrange1991equations,sabbah2008explicit} establishing the stationary phase formula. For the levels, some of the statements are proven in \cite{hiroe2017ramified} using results on singularities of curves. Let us discuss how to derive them all more directly from results in \cite{doucot2021diagrams}. In the case (3) where $I$ is at infinity and of slope $>1$, the statement for the levels can be deduced quite directly from \cite[Proposition~4.8 and the proof of Theorem~4.9]{doucot2021diagrams}, as follows. Proposition~4.8 says that all exponents of $F\cdot I$ belong to a certain set~${K\subset \mathbb Q_{>0}}$ determined explicitly from $I$.
Next, although it does not use the terminology of levels, the proof of Theorem~4.9 shows that, inside $K$, the subset $\Levels(F\cdot I)$ is exactly~${\frac{r}{s-r}\Levels(I)}$. This is because, in the explicit formula for the number of edge-loops~$B_{I,I}$ of a Stokes circle~$I$, given by Lemma~4.6 of loc.\ cit., the nonzero positive terms are in one-to-one correspondence with the levels of $I$ (by definition of the notion of levels). The proof of Theorem~4.9 establishes that there is a one-to-one correspondence between the nonzero positive terms in $B_{I,I}$ and $B_{F\cdot I, F\cdot I}$, i.e., there is a bijection between $\Levels(I)$, and $\Levels(F\cdot I)$, given by $\Levels(F\cdot I)=\frac{r}{s-r}\Levels(I)$.

The other cases are obtained in a completely similar way, by using the expressions for the Legendre transform for the different cases given, e.g., in \cite[Section~3.1]{doucot2021diagrams}. The point is that the difference between the cases concerns only the transformation of the leading terms of the Stokes circles, i.e., the change of slopes, without affecting what happens for the subleading terms.
\end{proof}

To fully determine the fission data of the Fourier transform, we must also describe the fission exponents between distinct Stokes circles of the Fourier transform which are at the same singularities. This can be done explicitly:

\begin{Proposition}\label{prop:fission_fourier_transform}
Let $I$, $J$ be two distinct Stokes circles. There are the following possibilities:
\begin{itemize}\itemsep=0pt
\item If $\pi(I)=\pi(J)=a\neq \infty$, then $\pi(F\cdot I)=\pi(F\cdot J)=\infty$ and
\begin{itemize}\itemsep=0pt
\item If $I$ and $J$ have no common part, then denoting $k=\slope(I)$, $l=\slope(k)$, we have $f_{I,J}=\max(\frac{k}{k+1},\frac{l}{l+1})$.
\item If $I$ and $J$ have a common part, denoting by $k$ its slope, we have $f_{F\cdot I, F\cdot J}=\frac{k}{k+1}f_{I,J}$.
\end{itemize}
\item If $\pi(I)\neq \infty, \pi(J)\neq \infty$ and $\pi(I)\neq \pi(J)$, then $\pi(F\cdot I)=\pi(F\cdot J)=\infty$ and $f_{F\cdot I, F\cdot J}=1$.
\item If $\pi(I)=\pi(J)=\infty $, $I, J$ are both of slope $\leq 1$ and are of the form $I=\cir{\alpha z + q}$, $J=\cir{\alpha z + q'}$ with $\alpha\in \IC$, and $q$, $q'$ of slope $<1$, then $\pi(F\cdot I)=\pi(F\cdot J)=\alpha$ and
\begin{itemize}\itemsep=0pt
\item If $q$ and $q'$ have no common part and respective slopes $k$, $l$, then $f_{F\cdot I, F\cdot J}=\max\bigl(\frac{k}{1-k},\allowbreak\frac{l}{1-l}\bigr)$.
\item If $q$ and $q'$ have a common part $q_c$, denoting by $k$ its slope, we have $f_{F\cdot I, F\cdot J}=\frac{k}{1-k}f_{q,q'}$.
\end{itemize}
\item If $\pi(I)=\pi(J)=\infty $, and $I$, $J$ are both of slope $>1$, then $\pi(F\cdot I)=\pi(F\cdot J)=\infty$ and
\begin{itemize}\itemsep=0pt
\item If $I$ and $J$ have no common part, then denoting $k=\slope(I)$, $l=\slope(J)$, we have $f_{I,J}=\max\bigl(\frac{k}{k-1},\frac{l}{l-1}\bigr)$.
\item If $I$ and $J$ have a common part, let $k$ denote its slope, then $f_{F\cdot I, F\cdot J}=\frac{k}{k-1}f_{I,J}$.
\end{itemize}
\item Otherwise, $\pi(I)\neq \pi(J)$.
\end{itemize}
\end{Proposition}

\begin{proof}
All statements where there is no common part involved follow directly from the previous proposition. For the statements in the cases where there is a common part, in the case where~$I$,~$J$ are at infinity and of slope $>1$ this is established in the proof of \cite[Theorem~5.9]{doucot2021diagrams} (the word `fission exponent' was not used there: in terms of the terminology used there it corresponds to the highest number among the exponents of the leading term of the \emph{different parts} of $F\cdot I$ and~$F\cdot J$). The other cases are similar, using again the expressions for the Legendre transform for the different cases given in \cite[Section~3.1]{doucot2021diagrams}, and noticing that this just changes the transformation of the slopes, without affecting the subleading terms (hence in particular the fission exponents when the common part is nonzero).
\end{proof}

\subsubsection[General symplectic SL\_2(C) transformations]{General symplectic $\boldsymbol{\SL_2(\IC)}$ transformations}

The Fourier--Laplace transform is part of a larger group of transformations acting on modules over the Weyl algebra. Indeed, to any matrix
$ A=\left(\begin{smallmatrix}
a & b\\
c & d
\end{smallmatrix}\right)
$
in $\SL_2(\IC)$, we can associate an automorphism of the Weyl algebra $A_1=\IC[z]\cir{\partial_z}$ given by $z\mapsto az+b\partial_z, \partial_z \mapsto cz+d\partial_z$. This induces an action of the group $\SL_2(\IC)$ of symplectic transformations on modules over the Weyl algebra (see \cite[p.~87]{malgrange1991equations}).

The group of symplectic transformations is generated by three types of elementary transformations:
\begin{itemize}\itemsep=0pt
\item the Fourier--Laplace transform $F$, corresponding to the matrix $\left(\begin{smallmatrix}
0 & 1\\
-1 & 0
\end{smallmatrix}\right)$;
\item twists at infinity $T_{\lambda}$, for $\lambda\in \IC$, corresponding to the matrix $\left(\begin{smallmatrix}
1 & \lambda\\
0 & 1
\end{smallmatrix}\right)$;
\item scalings $S_\lambda$, for $\lambda\in \IC^*$, corresponding to the matrix $\left(\begin{smallmatrix}
\lambda^{-1} & 0\\
0 & \lambda
\end{smallmatrix}\right)$.
\end{itemize}

The geometric interpretation of twists and scalings is the following. The twist $T_\lambda$ corresponds to taking the tensor product with the rank one module $(\IC[z],\partial_z+\lambda z)$, and the scaling $S_\lambda$ corresponds to perform the change of variable $z\mapsto z/\lambda$ on $\IP^1$.

As for the Fourier transform, any element of $\SL_2(\IC)$ induces a transformation on modules over the Weyl algebra, and on irreducible connections on Zariski open subsets of the affine line. Moreover, there exists again a corresponding formal transformation (that we will also denote by~$A$) on modified formal data with a similar structure. The action of twists and scalings on formal data and fission data is much simpler than for the Fourier transform.

\begin{Lemma}Let $I$ be a Stokes circle, and $\breve{\mathcal C}_I\in \GL_n(\mathbb C)$ a conjugacy class. Let $\breve{\mathcal C}_{T_\lambda\cdot I}\in \GL_n(\mathbb C)$ be the conjugacy class such that \smash{$T_\lambda\cdot\bigl(I,\breve{\mathcal C}_I\bigr)=(T_\lambda\cdot I, \breve{\mathcal C}_{T_\lambda\cdot I})$}. We have the following:
\begin{itemize}\itemsep=0pt
\item If $\pi(I)=\infty$, then $\pi(T_\lambda\cdot I)=\infty$, and if $I=\cir{q}_\infty$, we have $T_\lambda\cdot I=\cir{q-\frac{\lambda}{2}z^2}_\infty$. In particular, $T_\lambda\cdot I$ has the same ramification order and set of levels as $I$, $\slope(T_\lambda\cdot I)=\max(\slope(\cir{q}), \slope\bigl(\cir{q-\frac{\lambda}{2} z^2}\bigr)$, for $\lambda\in \IC$.
\item Otherwise, if $\pi(I)\neq\infty$, $T_\lambda\cdot I=I$.
\end{itemize}
In any case \smash{$\breve{\mathcal C}_{T_\lambda\cdot I}=\breve{\mathcal C}_I$}.
\end{Lemma}

In brief, all exponential factors at infinity are shifted by $-\frac{\lambda}{2} z^2$ and all other exponential factors are left unchanged. As a consequence, the twist $T_\lambda$ preserves the global fission data and classes of global admissible deformations. The twist may change the slopes of the exponential factors at infinity, but leaves the slopes of Stokes circles at finite distance unchanged.

\begin{Lemma}
Let $I$ be a Stokes circle such that $\pi(i)=a\in \IP^1$, and $\breve{\mathcal C}_I\in \GL_n(\mathbb C)$ a conjugacy class. Let $\breve{\mathcal C}_{S_\lambda\cdot I}\in \GL_n(\mathbb C)$ be the conjugacy class such that $S_\lambda\cdot\bigl(I,\breve{\mathcal C}_I\bigr)=\bigl(S_\lambda\cdot I, \breve{\mathcal C}_{T_\lambda\cdot I}\bigr)$. For $\lambda\in \mathbb C^*$, $\pi(S_\lambda\cdot I)= a/\lambda$, and writing $I=\cir{q\bigl(z_a^{-1}\bigr)}$ with \smash{$q\in \IC\big[z_a^{-1/r}\big]$}, using the local coordinate $z_a:=z-a$ if $a\neq \infty$, and $z_\infty:=z^{-1}$ if $a=\infty$, we have $S_\lambda \cdot I=\cir{q\bigl(\lambda z_{a/\lambda}^{-1}\bigr)}$. Furthermore~\smash{$\breve{\mathcal C}_{S_\lambda\cdot I}=\breve{\mathcal C}_I$}.
\end{Lemma}

As a consequence, for any $\lambda\in \mathbb C^*$, the scaling $S_\lambda$ preserves the global fission data and classes of global admissible deformations, as well as the slopes of all Stokes circles.

A crucial property of the diagram associated to a connection is that it is invariant under the~$\SL_2(\IC)$ action.

\begin{Theorem}[\cite{doucot2021diagrams}]
Let $(E,\nabla)$ be an irreducible connection on a Zariski open subset of~$\IP^1$ {\rm(}different from a rank one connection with only a pole at infinity of order less than $2)$, and~$A\in \SL_2(\mathbb C)$. Then
$
\Gamma(E,\nabla)=\Gamma(A\cdot (E,\nabla))$.
\end{Theorem}

\section{Classes of nearby representations}

\label{section:basic_reps}

In this section, we use the results of the previous paragraph to describe the formal data and classes of admissible deformations of the elements of the orbit of a connection $(E,\nabla)$ under the~$\SL_2(\mathbb C)$ action.

\subsection{Generic form}

\begin{Definition}
Let $I\in \pi_0(\Ical)$ be a Stokes circle. We say that $I$ is of \emph{generic form} if $\pi(I)=\infty$, and $\slope(q)\leq 2$.

If \smash{$I=\cir{q}_\infty$} is of generic form, the Fourier sphere coefficient $\lambda(I)$ of $I$ is the unique complex number $a\in \IC$ such that $\cir{\frac{a}{2}z^2+q}$ has slope $<2$. Explicitly, this means that $q$ is of the form
\[
q=-\frac{a}{2}z^2+q_{<2},
\]
where the exponential factor $q_{<2}$ has slope $<2$. We will denote $I_{<2}:=\cir{q_{<2}}$.

If $I$ is not of generic form, we set $\lambda(I)=\infty$.
\end{Definition}

For $A=\left(\begin{smallmatrix}
a & b\\
c & d
\end{smallmatrix}\right) \in \SL_2(\IC)$, let $h_A\colon \mathbb P^1\to \mathbb P^1$ denote the corresponding homography
\[
h_A \colon \ z\mapsto \frac{az+b}{cz+d}.
\]

We call the Fourier sphere, and denote by $\mathbf P^1$, the copy of $\mathbb P^1$ where the coefficients $\lambda(I)$ live (this terminology was introduced in \cite{boalch2012simply}).

The following result says that $\SL_2(\IC)$ acts on the Fourier sphere by homographies.

\begin{Lemma}
Let $I\in \pi_0(\Ical)$ be a Stokes circle, $a=\lambda(I)$ its Fourier sphere coefficient, and $A\in SL_2(\IC)$. Then
$
\lambda(A\cdot I)=h_A(\lambda(I))$.
\end{Lemma}

\begin{proof}In the case where $a\neq \infty$ and $h_A(a)\neq \infty$, this is the result of \cite[Proposition~4.13]{doucot2021diagrams}, and one checks that this is still true in the remaining cases.
\end{proof}

This immediately implies the following.

\begin{Corollary}
Let $I\in \pi_0(\Ical)$ be a Stokes circle. Then for $A$ in an open dense subset of $\SL_2(\mathbb C)$, $A\cdot I$ is of generic form.
\end{Corollary}

\begin{Definition}
Let $\bm{\breve\Theta}$ be a global modified irregular class. We say that it is of generic form if each of its Stokes circles is of generic form. Similarly if $(E,\nabla)$ is an irreducible connection on a~Zariski open subset of $\mathbb P^1$, we say that it is of generic form if its global modified irregular class is of generic form.
\end{Definition}

In particular, a connection of generic form only has a singularity at infinity, i.e., $\bm{\breve\Theta}=\Theta_\infty$. A global modified irregular class $\bm{\breve\Theta}$ is of generic form if and only if $\lambda\bigl(\bm{\breve\Theta}\bigr)\subset \IC=\mathbf{P}^1\smallsetminus\{\infty\}$, where~\smash{$\lambda\bigl(\bm{\breve\Theta}\bigr)$} denotes the set of Fourier sphere coefficients of all its Stokes circles. This immediately implies the following.

\begin{Lemma}\label{lemma3.5}
Let $\bm{\breve\Theta}$ be a global modified irregular class. Then for $A$ in an open dense subset of~$\SL_2(\IC)$, $A\cdot \bm{\breve\Theta}$ is of generic form.
\end{Lemma}

In particular, every $\SL_2(\IC)$ orbit contains an element of generic form. We can thus without loss of generality write any orbit as $\SL_2(\IC)\cdot (E,\nabla)$, with $(E,\nabla)$ of generic form.

As a preliminary step towards describing such orbits, let us first consider the case where $\Theta_\infty$ only has one Fourier sphere coefficient.

\begin{Proposition}\label{proposition:orbit_single_Fourier_coeff}
Let $(\Theta,\mathcal C)$ be formal data at infinity with $\Theta$ of generic form, and such that~${\lambda(\Theta)=\{a\}}$ for some $a\in \IC$. Let $A\in \SL_2(\IC)$. Let $\Theta^{<2}:=T_{-a}\cdot \Theta$ be the slope $<2$ irregular class obtained from $\Theta$ by subtracting $-\frac{a}{2}z^2$ to each of its Stokes circles. We have
\begin{itemize}\itemsep=0pt
\item If $h_A(a)=\infty$, then $A\cdot (\Theta,\mathcal C)$ is an admissible deformation of $F\cdot \bigl(\Theta^{<2},\mathcal C\bigr)$.
\item Otherwise, if $h_A(a)\neq\infty$, then $A\cdot (\Theta,\mathcal C)$ is an admissible deformation of $(\Theta,\mathcal C)$ , or equivalently of $\bigl(\Theta^{<2},\mathcal C\bigr)$.
\end{itemize}
\end{Proposition}

\begin{proof}This is basically a consequence of the fact that any element in $\SL_2(\IC)$ can be factorised as a product of elementary transformations, and among those only the Fourier transform may act nontrivially on global fission data. More precisely, we can use \cite[Lemma~4.15]{doucot2021diagrams} to write down explicit factorisations of~$A$. If $h_A(a)=\infty$, then we have $h_{AT_a F^{-1}}(\infty)=\infty$, so from this lemma, there exists $v\in \mathbb C^*, \rho\in \mathbb C$ such that $A T_a F^{-1}=S_v T_{\rho}$, so that $A=S_v T_\rho F T_{-a}$. This implies $A\cdot (\Theta,\mathcal C)=S_v T_\rho F T_{-a}\cdot (\Theta,\mathcal C)=S_v T_\rho F\cdot \bigl(\Theta^{<2},\mathcal C\bigr)$, and since twists and scalings do not change the equivalence class of admissible deformations $A\cdot (\Theta,\mathcal C)$ is an admissible deformation of $F\cdot \bigl(\Theta^{<2},\mathcal C\bigr)$.

Let us now consider the case $h_A(a)\neq\infty$. If $h_A(\infty)\neq \infty$ then from the lemma $A$ admits a~factorisation $A=S_{\nu}T_{\rho}$ with $v\in \mathbb C^*, \rho\in \mathbb C$, which implies that $A\cdot \Theta$ is an admissible deformation of~$\Theta$. Otherwise, from the lemma, $A$ admits a factorisation $A=T_{\mu}F S_{\nu} T_{\rho}$, with~${v\in \mathbb C^*}$, ${\mu,\rho\in \mathbb C}$. In turn $\Theta':=S_{\nu}T_{\rho}\cdot \Theta$ is an admissible deformation of $\Theta$, and $\lambda(\Theta')=\{a'\}$, with~${a'=h_{S_{\nu}T_{\rho}}(a)}$, and we have $a'\neq 0$ otherwise we would have $h_A(a)=h_{S_{\nu}F}(a')=\infty$.
Now, applying Propositions~\ref{prop:slopes_Fourier_transform} and~\ref{prop:fission_fourier_transform} to $\Theta'$, the fact that all exponential factors in $\Theta'$ are of slope~2 and have the same nonzero leading coefficient implies that this is also the case for $F\cdot \Theta'$, and that $F\cdot (\Theta',\mathcal C)$ is an admissible deformation of $(\Theta',\mathcal C)$. Since $T_{\mu}$ does not change the equivalence class of admissible deformations, we finally obtain that $A\cdot (\Theta,\mathcal C)$ is an admissible deformation of $(\Theta,\mathcal C)$, which concludes the proof.
\end{proof}

In particular, we can write down explicit formulas for the rank at infinity of the elements of the orbit of $\Theta$. Let $N$ be the set of Stokes circles of $\Theta$, and let us set
\begin{gather*}
N^+:=\left\{ \cir{q} \in N\mid \slope\left(q+\frac{a}{2}z^2\right)>1\right\}\subset N, \\ N^-:=\left\{ \cir{q} \in N\mid \slope\left(q+\frac{a}{2}z^2\right)\leq1\right\} \subset N,
\end{gather*}
so that $N=N^+\sqcup N^-$. Let us write the elements of $N^-$ as
\smash{$I^-_j=\cir{-\frac{a}{2}z^2+ b_{j}z+ q^-_j}
$}
for $j\in\{1,\dots, |N^-|\}$, where \smash{$\slope(q^-_j)<1$}, and set \smash{$s^-_j:=\irr(q^-_j)$}, \smash{$r^-_j:=\ram(q^-_j)$}. Similarly, we write
\smash{$I^+_j=\cir{-\frac{a}{2}z^2+q^+_j}$},
for $j\in\big\{1,\dots, \big|N^+\big|\big\}$, with $\slope\bigl(q^+_j\bigr)>1$, and set \smash{$s^+_j:=\irr\bigl(q^+_j\bigr)$}, \smash{$r^+_j:=\ram\bigl(q^+_j\bigr)$}.

\begin{Lemma}\label{lemma:rank_one_fourier_sphere_coeff}
Let $n^{\pm}_j$ denote the multiplicities of the circles in $N^{\pm}$, and set $R^+_j:=s^+_{j}-r^+_j$, for~\smash{$j\in\big\{1,\dots, \big|N^+\big|\big\}$}. Let $A\in \SL_2(\IC)$. We have
\begin{itemize}\itemsep=0pt
\item If $h_A(a)=\infty$, then $\rank_\infty(A\cdot \Theta)=
\sum_j n^+_j R^+_j$.
\item If $h_A(a)\neq\infty$, then $\rank_\infty(A\cdot \Theta)=
\sum_j n^+_j r^+_j + n^-_j r^-_j$.
\end{itemize}
\end{Lemma}

\begin{proof}
This follows directly from Lemma~\ref{lemma3.5} and Proposition~\ref{prop:slopes_Fourier_transform}, by observing that for the case $h_A(a)=\infty$, the Stokes circles in $N^-$ are sent by $A$ to finite distance (with the location of $A\cdot I_j^-$ determined by the coefficient $b_j$), while the circles in $N^+$ become circles of slope $>1$ at infinity, with \smash{$A\cdot I_j^+$} having ramification order \smash{$R_j^+$}.
\end{proof}

\subsection{The classes of nearby representations}

We now discuss general orbits. Let $\bigl(\bm{\breve\Theta},\bm{\breve\Ccal}\bigr)$ be effective global modified formal data, with $\bm{\breve\Theta}=\Theta_\infty$ of generic form. Let $a_1,\dots ,a_k\in \IC$ be the distinct Fourier sphere coefficients of its Stokes circles, i.e., we have $\lambda(\Theta_\infty)=\{a_1,\dots, a_k\}\subset \mathbf P^1$.

Let $N$ be the set of Stokes circles of $\bm{\breve\Theta}$. For $I\in N$, we denote by $n_I\in \mathbb N_{>0}$ the multiplicity of $I$ in $\bm{\breve\Theta}$. We partition $N$ as follows: for $i\in\{1,\dots, k\}$, we set
\begin{itemize}\itemsep=0pt
\item $N_i:=\{ I \in N \mid \lambda(I)=a_i\}$,
\item $N_i^+:=\big\{ \cir{q} \in N_i \mid \slope\bigl(q+\frac{a_i}{2}z^2\bigr) >1\big\}\subset N_i$,
\item $N_i^-:=\big\{ \cir{q} \in N_i \mid \slope\bigl(q+\frac{a_i}{2}z^2\bigr) \leq 1\big\}\subset N_i$,
\end{itemize}
so that
$
N=N_1\sqcup \dots \sqcup N_k$, $ N_i=N_i^{+}\sqcup N_i^{-}$.
We also set $N^+:=N_1^+\sqcup \dots \sqcup N_k^+$, and $N^-:=N_1^-\sqcup \cdots \allowbreak \sqcup N_k^-$. For $i\in\{1,\dots, k\}$, let $\Theta_i:=\sum_{I\in N_i} n_I I$ be the irregular class obtained from~$\Theta_\infty$ obtained by discarding all circles with Fourier sphere coefficient different from~$a_i$, $\Ccal_i$ the corresponding collection of conjugacy classes, and let $\Theta^{<2}_i$ be the irregular class at infinity obtained from $\Theta_i$ by subtracting $-\frac{a_i}{2}z^2$ to each Stokes circle (in particular all circles of $\Theta_i$ have slope $<2$).

\begin{Theorem}\label{theorem:def_of_fundamental_representations}
Let $A\in \SL_2(\IC)$, then
\begin{itemize}\itemsep=0pt
\item If $h_A(a_i)=\infty$ for some $i\in \{1,\dots, k\}$, then $A\cdot \bigl(\bm{\breve\Theta},\bm{\breve \Ccal}\bigr)$ is an admissible deformation of
\[
F\cdot \bigl(\Theta_i^{<2},\Ccal_i\bigr) + \sum_{l\neq i} (\Theta_l,\Ccal_l).
\]

\item Otherwise, if $h_A(a_i)\neq \infty$ for each $i\in\{1,\dots, k\}$, then $A\cdot \bigl(\bm{\breve\Theta},\bm{\breve \Ccal}\bigr)$ is an admissible deformation of $\bigl(\bm{\breve\Theta},\bm{\breve \Ccal}\bigr)$.
\end{itemize}
\end{Theorem}

\begin{proof}
This follows from applying Proposition~\ref{proposition:orbit_single_Fourier_coeff} to every $(\Theta_i,\mathcal C_i)$.
\end{proof}

We notice that the corresponding classes of admissible deformations are in fact completely determined by the (short) fission tree of $\bm{\breve\Theta}$.

\begin{Proposition}\label{prop:nearby_reps_determined_from_tree}
Keeping previous notations, let $\mathbf T$ be the fission tree of $\bm{\breve \Theta}$. Then for $i\in\{1,\dots, k\}$, the class of admissible deformations $(\mathbf F_i, \bm{\mathcal C}_i)$ of the nonmodified formal data associated to the formal data
\[
F\cdot \bigl(\Theta_i^{<2},\mathcal C_i\bigr) + \sum_{l\neq i} (\Theta_l,\mathcal C_l)
\] is fully determined by the pair $\bigl(\mathbf T, \bm{\breve\mathcal C}\bigr)$, in an explicit way.
\end{Proposition}

\begin{proof}
This follows from Propositions~\ref{prop:slopes_Fourier_transform} and~\ref{prop:fission_fourier_transform}. Indeed, applying them to determine the invariants of $F\cdot \bigl(\Theta_i^{<2},\mathcal C_i\bigr)$ in terms of those of $\bigl(\Theta_i^{<2},\mathcal C_i\bigr)$, we obtain that
\begin{itemize}\itemsep=0pt
\item The fission datum of the global modified formal data $F\cdot \bigl(\Theta_i^{<2},\mathcal C_i\bigr) + \sum_{l\neq i} (\Theta_l,\mathcal C_l)$ is fully determined by $\mathbf T$, in an explicit way.
\item The slope of any Stokes circle $I_a$ at a singularity $a\neq \infty$ at finite distance of $F\cdot \bigl(\Theta_i^{<2}, \mathcal C_i\bigr)$ is determined by $\mathbf T$ in an explicit way. In turn, this determines the fission exponent~\smash{$f_{I_a, \langle0\rangle_a}$} between $I_a$ and the tame circle at $a$, and thus determines the fission datum of the nonmodified formal data associated to \smash{$F\cdot \bigl(\Theta_i^{<2},\mathcal C_i\bigr) + \sum_{l\neq i} (\Theta_l,\mathcal C_l)$}.\hfill $\qed$
\end{itemize}\renewcommand{\qed}{}
\end{proof}

This makes it possible to take this as a definition of the classes of nearby representations of the enriched tree $\bigl(\mathbf T,\bm{\breve\mathcal C}\bigr)$.

\begin{Definition}\label{def:nearby_representations}
Let $\bigl(\mathbf T, \bm{\breve\mathcal C}\bigr)$ be an enriched tree, i.e., $\mathbf T$ is a short fission tree, $N$ its set of leaves and \smash{$\bm{\breve\mathcal C}$} the datum of a conjugacy class $\mathcal C_i\in \GL_{n_i}(\mathbb C)$ for each leaf $i\in N$, where $n_i$ is the multiplicity of $i$. Assume that there exists an irregular class $\Theta_\infty$ of generic form with fission tree $\mathbf T$ such that \smash{$\bigl(\Theta_\infty,\bm{\breve\Ccal}\bigr)$} are effective modified global formal data. For $i\in\{1,\dots, k\}$, keeping previous notations, we define the $i$-th nongeneric class of nearby representations of \smash{$\bigl(\mathbf T, \bm{\breve\mathcal C}\bigr)$} as the class $(\mathbf F_i, \bm{\mathcal C}_i)$ of admissible deformations of the nonmodified formal data associated to the formal data
\[
F\cdot \bigl(\Theta_i^{<2},\mathcal C_i\bigr) + \sum_{l\neq i} (\Theta_l,\mathcal C_l).
\]
(This is well-defined since by the previous proposition this does not depend on the choice of~$\Theta_\infty$.)
\end{Definition}

Theorem~\ref{theorem:def_of_fundamental_representations} can then be reformulated in the following way, which justifies our terminology.

\begin{Corollary}\label{cor:we_get_classes_of_nearby_reps}
Let $(\mathbf \Theta, \bm{\Ccal})$ be an effective wild Riemann surface with boundary data on $\mathbb P^1$ of generic form, with class of admissible deformations $(\{\mathbf T\}, \bm{\mathcal C})$, and let $(\mathbf F_i, \bm{\Ccal_i})$, with $i\in\{1,\dots, k\}$, be the classes of nearby representations of $(\mathbf T, \bm{\Ccal})$ {\rm(}in the sense of the previous definition{\rm)}. Then the classes of representations of the nonabelian Hodge space $\mathcal M(\mathbf \Theta, \bm{\Ccal})$ which are nearby $(\mathbf \Theta, \bm{\Ccal})$ {\rm(}in the sense of Section~{\rm\ref{subsec:weak_reps_intro})} are exactly $(\{\mathbf T\}, \bm{\Ccal})$ and $(\mathbf F_i, \bm{\mathcal C}_i)$ for $i\in\{1,\dots, k\}$.
\end{Corollary}

By abuse of language, we will sometimes allow ourselves to say that $\{\mathbf T\}$ and $\mathbf F_i$, for $i\in\{1,\dots, k\}$ are the nearby representations of $\mathbf T$.

Finally, all this allows us to have a well-defined notion of enriched tree associated to any irregular connection on $\mathbb P^1$.

\begin{Proposition}\label{prop: definition_of_enriched_trees}
Let $(E,\nabla)$ be an irreducible algebraic connection on a Zariski open subset of $\mathbb P^1$, such that $(E,\nabla)$ is not a rank one connection with a unique singularity at infinity with Poincar\'e--Katz rank $\leq 1$. Let $\bigl(\bm{\breve{\Theta}}, \bm{\breve{\mathcal C}}\bigr)$ be its modified global formal data, and $A\in\SL_2(\mathbb C)$ such that \smash{$A\cdot \bigl(\bm{\breve\Theta}, \bm{\breve\mathcal C}\bigr)$} is of generic form. Let $\mathbf T$ be the $($short$)$ fission tree of \smash{$A\cdot \bm{\breve\Theta}$}, and $\bm{\mathcal C}$ the collection of conjugacy classes of \smash{$A\cdot \bigl(\bm{\breve\Theta}, \bm{\breve\mathcal C}\bigr)$}. Then the pair $(\mathbf T, \bm{\mathcal C})$ does not depend on the choice of $A$.
\end{Proposition}

\begin{proof}
This is an immediate consequence of Theorem~\ref{theorem:def_of_fundamental_representations}.
\end{proof}

\begin{Definition}\label{def:enriched_tree}
Keeping the setup of the proposition, we say that $(\mathbf T, \bm{\mathcal C})$ is the \emph{enriched tree} of $(E,\nabla)$ and write $\mathscr T(E,\nabla):=(\mathbf T, \bm{\mathcal C})$.
\end{Definition}

Proposition~\ref{prop: definition_of_enriched_trees} implies that the enriched tree is invariant under symplectic transformations.

\begin{Corollary}\label{cor:symplectic_invariance_enriched_tree}
If $(E,\nabla)$ is irreducible and is not a rank one connection with a unique singularity at infinity of Poincar\'e--Katz rank $\leq 1$, so that $A\cdot (E,\nabla)$ is well-defined for any $A\in \SL_2(\mathbb C)$, we have
$
\mathscr T(A\cdot (E,\nabla))=\mathscr T(E,\nabla)$.
\end{Corollary}

\subsection{Properties of classes of nearby representations}

Proposition~\ref{prop:nearby_reps_determined_from_tree} gives an explicit method to determine the nearby classes of representations from an enriched tree, using Propositions~\ref{prop:slopes_Fourier_transform}
and~\ref{prop:fission_fourier_transform}. We will compute many examples in Section~\ref{section:painleve_examples}. Before this, in this paragraph, let us state several general properties satisfied by the classes of nearby representations.

First, we can get explicit formulas for the rank of the nearby representations. For $i\in\{1,\dots, k\}$, let us denote by \smash{$\cir{q^+_{i,j}}$} the elements of \smash{$N_i^+$}, and write
\smash{$\cir{q^+_{i,j}}=\cir{-\frac{a_i}{2} z^2+ q'^+_{ij}}$}, with~\smash{$q'^+_{ij}$} having ramification \smash{$r^+_{ij}$} and slope \smash{$s^+_{ij}/r^+_{ij}>1$}. Similarly, let us denote by \smash{$\cir{q^-_{i,j}}$} the elements of~$N_i^-$, and write
\smash{$\cir{q^-_{i,j}}=\cir{-\frac{a_i}{2} z^2+b_{ij}z+ q'^-_{ij}}$}, with $q'^-_{ij}$ having ramification $r^-_{ij}$ and slope $s^-_{ij}/r^-_{ij}<1$. Let \smash{$n_{ij}^\pm$} be the multiplicities of the active circles.

\begin{Proposition}\label{prop:rank_of_representations}
 Let us set $R^+_{ij}=s^+_{ij}-r^+_{ij}$.
\begin{enumerate}\itemsep=0pt
\item[$(1)$] The generic class of representations has rank
\[
r_{\rm gen}=\sum_i \Bigg(\sum_j n^+_{ij}r^+_{ij} + \sum_j n^-_{ij}r^-_{ij}\Bigg).
\]
\item[$(2)$] The $i$-th nongeneric class of representations has rank
\[
r_{i}=\sum_{j} n^+_{ij} R^+_{ij}+ \sum_{l\neq i} \Bigg(\sum_j n^+_{l j}r^+_{l j} + \sum_j n^-_{l j}r^-_{l j}\Bigg).
\]
\end{enumerate}
\end{Proposition}

\begin{proof}
This follows directly from Lemma \ref{lemma:rank_one_fourier_sphere_coeff}.
\end{proof}

We can also express the number of singularities at finite distance of the classes of nearby representations.

\begin{Proposition}
\label{prop:number_of_singularities}
Let $(E,\nabla)$ of generic form, and $(E',\nabla')\in \SL_2(\mathbb C)\cdot (E,\nabla)$. If $(E',\nabla')$ is in the $i$-th nongeneric class of representations, its number of singularities at finite distance is equal to the number of vertices of height 1 in the subtree of $\mathbf T$ with set of leaves $N_i^-$.
\end{Proposition}

\begin{proof}
By definition of fission trees, keeping previous notations, this number of vertices is equal to the number of distinct coefficients $b_{ij}$ among the Stokes circles of $N_i^-$.
The result is then a~direct consequence of item (2) of Proposition~\ref{prop:slopes_Fourier_transform}: the position of the singularity above which~$F\cdot \cir{q^-_{i,j}}$ lies is determined by the coefficient~$b_{ij}$.
\end{proof}

In particular, if $N_i^-=\varnothing$, then the $i$-th nongeneric class of representations only has a singularity at infinity. In general, the class of nearby representations all have different ranks and number of singularities.

Finally, we can also read from the short fission tree when two nongeneric representations have the same fission forest.

\begin{Proposition} Keeping previous notations, for $i,j\in\{1,\dots k\}$, with $i\neq j$, the fission forests $\mathbf F_i$ and $\mathbf F_j$ are isomorphic if and only if the $i$-th and $j$-th principal subtrees of the short fission tree $\mathbf T$ are isomorphic.
\end{Proposition}

\begin{proof}
Since isomorphism classes of fission trees characterise admissible deformations, it is equivalent to show that the $i$-th and the $j$-th nongeneric forests are isomorphic if and only if~$\Theta_i^{<2}$ and $\Theta_j^{<2}$ are admissible deformations of each other. Let us first assume that $\Theta_i^{<2}$ and~\smash{$\Theta_j^{<2}$} are admissible deformations of each other. We have to be careful that, in general, if two irregular class $\Theta$ and $\Theta'$ at infinity are admissible deformations of each other, this does not imply that their Fourier transforms are also admissible deformations of each other. Indeed, if $\cir{q}$ is a~Stokes circle for $\Theta$, and $\cir{q'}$ is a Stokes circle for $\Theta'$ corresponding to $\cir{q}$ under some admissible deformation, their slopes may differ, and in turn $F\cdot \cir{q}$ and $F\cdot \cir{q'}$ may not be admissible deformations of each other. However here, because all Stokes circles of $\Theta_i^{<2}$ and $\Theta_j^{<2}$ are of slope $<2$ this situation cannot happen, and we have that \smash{$F\cdot \Theta_i^{<2}$} and \smash{$F\cdot \Theta_j^{<2}$} are admissible deformations of each other. In turn, the irregular classes $\Theta_i^{<2} + \sum_{l\neq i} \Theta_l$ and $\Theta_j^{<2} + \sum_{l\neq j} \Theta_l$ are admissible deformations of each other, so the $i$-th and $j$-th representations coincide. Conversely, if the $i$-th and $j$-th classes of representations coincide, then in particular $\sum_{l\neq i} \Theta_l$ and $\sum_{l\neq j} \Theta_l$ are admissible deformations of each other, which implies that $\Theta_i$ and $\Theta_j$ are admissible deformations of each other, and the conclusion follows.
\end{proof}

\begin{Remark}
\label{remark:fourier_and_other_operations}
The classes of nearby representations arise from the action of $\SL_2(\mathbb C)$ on connections on the affine line. However, as already mentioned in the introduction, there is a larger class of natural operations on connections: indeed one can consider M\"obius transformations, which have the effect of changing the point at infinity in the identification $\mathbb P^1=\mathbb C\cup\{\infty\}$, as well as the operation consisting in tensoring a connection with an arbitrary rank one connection on a~Zariski open subset of $\mathbb P^1$. This larger class of operations on connections on~$\mathbb P^1$ plays an essential role in the extension due to Deligne--Arinkin \cite{arinkin2010rigid,deligne2006letter} of the Katz algorithm for rigid local systems to the case of irregular connections, which provides a way to reduce any irreducible rigid connection on~$\mathbb P^1$ to the trivial rank one connection by repeated application of one of these operations. Since the symplectic $\SL_2(\mathbb C)$ transformations do not commute with these other operations in general, the latter do not preserve the diagram or the classes of nearby representations.
\end{Remark}

\section{Basic representations for Painlev\'e moduli spaces}
\label{section:painleve_examples}

In this section, as a set of examples illustrating our general discussion, we describe the basic representations of all Painlev\'e moduli spaces obtained from the standard rank 2 Lax representation of each Painlev\'e equation. Recall that, although they were not first discovered in this way historically, all Painlev\'e equations can be obtained from isomonodromic deformations. Each Painlev\'e equation arises from isomonodromic deformations of connections with a specific type of singularity data, i.e., different Painlev\'e equations correspond to different nonabelian Hodge spaces. In the standard way to derive the Painlev\'e equations from isomonodromy \cite{jimbo1981monodromyI}, all Lax representations are of rank 2, and in all cases except for Painlev\'e~VI, these Lax representations have at least one irregular singularity.

For each case, starting from this standard Lax representation (in its abstract form, i.e., consisting of the datum of a wild Riemann surface), we will draw the corresponding short fission tree of minimal rank and describe explicitly the classes of nearby representations. While the cases of Painlev\'e~VI,~V,~IV all fall in the simply-laced case and were discussed by Boalch in~\cite{boalch2012simply}, it is still worth discussing them briefly from the point of view of short fission trees. On the other hand, describing all basic representations for Painlev\'e~III,~II,~I as well as for the moduli spaces corresponding to the so-called degenerate versions of Painlev\'e~III really requires using our general framework. Doing this, we recover among the classes of nearby representations many known alternative Lax representations for various Painlev\'e equations, e.g., from \cite{flaschka1980monodromy, joshi2007linearization, joshi2009linearization} and the other representations can be viewed as new abstract Lax representations.

\begin{Remark}
Notice that there is the following subtlety: since taking the Fourier transform does not always commute with taking admissible deformations, different connections having the same fission forest in general do not have the same classes of nearby representations, nor the same diagram. In the cases where this happens, among all possible wild Riemann surfaces with boundary data having the fission forest of the standard Painlev\'e Lax representation, we choose one such that its generic representation has minimal rank. Making this choice is necessary to allow us to recover several known alternative Lax representations for Painlev\'e equations and to recover the diagrams for Painlev\'e equations obtained in \cite{boalch2020diagrams}, all equal or related to affine Dynkin diagrams. Other choices would lead to different diagrams (typically with more vertices), and for each choice all nearby representations could also be explicitly determined, yielding even more abstract Lax representations, typically of higher rank.
\end{Remark}

To draw a fission tree $(\mathcal T,\mathbb V,\mathbb A,\mathbb L, h,n)$, we follow the conventions of \cite{boalch2025twisted}. The heights of vertices of the trees are indicated to the left of the figures. The different types of (nonempty) vertices are indicated by different colours: the \emph{mandatory} vertices, i.e., the elements of $\mathbb L$, are drawn in black. The \emph{inconsequential} vertices, i.e., the elements of $\mathbb A\smallsetminus \mathbb L$, are drawn in white. Concretely, a black vertex means that in an irregular class with fission tree $\mathcal T$ the corresponding coefficient is not allowed to vanish (hence the name mandatory), while a white vertex means that the coefficient is allowed to vanish. Unless otherwise specified, the multiplicity of the leaves in a fission tree is equal to 1.

In \cite{boalch2025twisted}, fission trees were first defined with an infinite trunk, then suitable truncations of the tree were introduced (see Section~3.5 there) so as to obtain rooted trees with a finite number of vertices: if $\Theta$ is an irregular class, the root of the truncated tree is defined as the (unique) point of the trunk with height $\lfloor \operatorname{Katz}(\Theta)+1\rfloor$.

Here, we will also draw truncated fission trees, and our convention for the height of the root is consistent with loc.\ cit., as follows. If $\mathcal F$ is the fission datum of some irregular class, we define its \textit{minimal Poincar\'e--Katz rank} as
\[
\operatorname{Katz}_{\min}(\mathcal F):=\min \{\operatorname{Katz}(\Theta) \mid \Theta \text{ irregular class with }\mathcal F(\Theta)=\mathcal F \}.
\]
In turn, the height of the root of the fission tree $\mathbf T$ corresponding to $\mathcal F$ is taken to be
\[
\eta:=\lfloor \operatorname{Katz}_{\min}(\mathcal F)+1\rfloor.
\]

To draw the diagrams, we follow the conventions of \cite{doucot2021diagrams}: unless otherwise specified the multiplicity of each vertex is 1, negative edges or loops are drawn in dashed lines, and for multiple edges/loops the multiplicity is indicated at their middle. Furthermore, here we use colours to draw the vertices, each colour corresponding to a subset in the partition defined by the short fission tree.

\subsection{Summary of basic Painlev\'e representations}

The short fission trees as well as the rank, number of singularities, and minimal Poincar\'e--Katz ranks of each class of basic representations for each Painlev\'e moduli space are listed below. In the next section, we will explain in more detail how to obtain these classes of representations starting from the standard Lax representations.

\subsubsection*{Painlev\'e I}

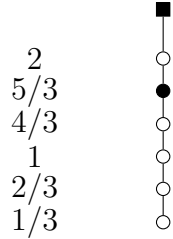
\begin{center}
\begin{tikzpicture}[scale=1.3]
\tikzstyle{mandatory}=[circle,fill=black,minimum size=5pt,draw, inner sep=0pt]
\tikzstyle{authorised}=[circle,fill=white,minimum size=5pt,draw,inner sep=0pt]
\tikzstyle{empty}=[circle,fill=black,minimum size=0pt,inner sep=0pt]
\tikzstyle{root}=[fill=black,minimum size=5pt,draw,inner sep=0pt]
\tikzstyle{indeterminate}=[circle,densely dotted,fill=white,minimum size=5pt,draw, inner sep=0pt]
 \node[root] (R) at (0,2){};
 \node[mandatory] (A5) at (0,1.67){};
 \node[authorised] (A4) at (0,1.33){};
 \node[authorised] (A3) at (0,1){};
 \node[authorised] (A2) at (0,0.67){};
 \node[authorised] (A1) at (0,0.33){};
 \node[empty] (A0) at (0,0){};
 \draw (R)--(A5)--(A4)--(A3)--(A2)--(A1)--(A0);
 %\draw (-1.3,2) node {$2$};
 \draw (-1.3,1.67) node {$5/3$};
 \draw (-1.3,1.33) node {$4/3$};
 \draw (-1.3,1) node {$1$};
 \draw (-1.3,0.67) node {$2/3$};
 \draw (-1.3,0.33) node {$1/3$};
\end{tikzpicture}
\end{center}

\begin{center}
\begin{tabular}{|c|c|c|c|}
\hline
 & Rank & Number of singularities & min. Katz ranks\\
 \hline
 Generic rep. & 3 & 1 & 5/3 \\
 \hline
 Nongeneric rep. & 2 & 1 & 5/2\\
 \hline
\end{tabular}
\end{center}

\subsubsection*{Painlev\'e~II}

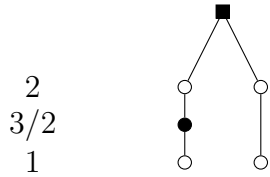
\begin{center}
\begin{tikzpicture}
\tikzstyle{mandatory}=[circle,fill=black,minimum size=5pt,draw, inner sep=0pt]
\tikzstyle{authorised}=[circle,fill=white,minimum size=5pt,draw,inner sep=0pt]
\tikzstyle{empty}=[circle,fill=black,minimum size=0pt,inner sep=0pt]
\tikzstyle{root}=[fill=black,minimum size=5pt,draw,inner sep=0pt]
\tikzstyle{indeterminate}=[circle,densely dotted,fill=white,minimum size=5pt,draw, inner sep=0pt]
 \draw (-2,1) node {$1$};
 \draw (-2,2) node {$2$};
 \draw (-2,1.5) node {$3/2$};
 \draw (-2,0.5) node {$1/2$};
 \node[root] (A3) at (0.5,3){};
 \node[authorised] (A2) at (0,2){};
 \node[authorised] (B2) at (1,2){};
 \node[authorised] (A1) at (0,1){};
 \node[authorised] (B1) at (1,1){};
 \node[authorised] (A1/2) at (0,0.5){};
 \node[mandatory] (A3/2) at (0,1.5){};
 \node[empty] (A0) at (0,0){};
 \node[empty] (B0) at (1,0){};
 \draw (A3)--(A2);
 \draw (A3)--(B2);
 \draw (A2)--(A3/2)--(A1)--(A1/2)--(A0);
 \draw (B2)--(B1)--(B0);
\end{tikzpicture}
\end{center}

\begin{center}
\begin{tabular}{|c|c|c|c|}
\hline
 & Rank & Number of singularities & min. Katz ranks\\
 \hline
 Generic rep. & 3 & 1 & 2\\
 \hline
 1st nongeneric rep. & 2 & 1 & 3 \\
 \hline
 2nd nongeneric rep. & 2 & 2 & 3/2, 0 \\
 \hline
\end{tabular}
\end{center}

\subsubsection*{Painlev\'e~III}

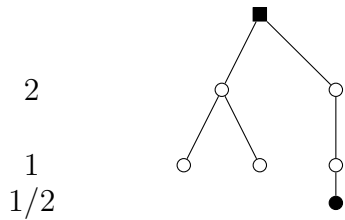
\begin{center}
\begin{tikzpicture}
\tikzstyle{mandatory}=[circle,fill=black,minimum size=5pt,draw, inner sep=0pt]
\tikzstyle{authorised}=[circle,fill=white,minimum size=5pt,draw,inner sep=0pt]
\tikzstyle{empty}=[circle,fill=black,minimum size=0pt,inner sep=0pt]
\tikzstyle{root}=[fill=black,minimum size=5pt,draw,inner sep=0pt]
\tikzstyle{indeterminate}=[circle,densely dotted,fill=white,minimum size=5pt,draw, inner sep=0pt]
 \draw (-2,1) node {$1$};
 \draw (-2,2) node {$2$};
 \draw (-2,0.5) node {$1/2$};
 \node[root] (A3) at (1,3){};
 \node[authorised] (A2) at (0.5,2){};
 \node[authorised] (B2) at (2,2){};
 \node[authorised] (A1) at (0,1){};
 \node[authorised] (B1) at (1,1){};
 \node[authorised] (C1) at (2,1){};
 \node[mandatory] (C1/2) at (2,0.5){};
 \node[empty] (A0) at (0,0){};
 \node[empty] (B0) at (1,0){};
 \node[empty] (C0) at (2,0){};
 \draw (A3)--(A2);
 \draw (A3)--(B2);
 \draw (A2)--(A1)--(A0);
 \draw (A2)--(B1)--(B0);
 \draw (B2)--(C1)--(C1/2)--(C0);
\end{tikzpicture}
\end{center}

\begin{center}
\begin{tabular}{|c|c|c|c|}
\hline
 & Rank & Number of singularities & min. Katz ranks\\
 \hline
 Generic rep. & 4 & 1 & 2\\
 \hline
 1st nongeneric rep. & 2 & 3 & 1/2, 0, 0\\
 \hline
 2nd nongeneric rep. & 2 & 2 & 1, 1\\
 \hline
\end{tabular}
\end{center}

\subsubsection*{Degenerate Painlev\'e~III}

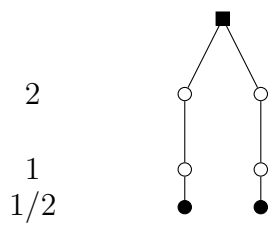
\begin{center}
\begin{tikzpicture}
\tikzstyle{mandatory}=[circle,fill=black,minimum size=5pt,draw, inner sep=0pt]
\tikzstyle{authorised}=[circle,fill=white,minimum size=5pt,draw,inner sep=0pt]
\tikzstyle{empty}=[circle,fill=black,minimum size=0pt,inner sep=0pt]
\tikzstyle{root}=[fill=black,minimum size=5pt,draw,inner sep=0pt]
\tikzstyle{indeterminate}=[circle,densely dotted,fill=white,minimum size=5pt,draw, inner sep=0pt]
 \draw (-2,1) node {$1$};
 \draw (-2,2) node {$2$};
 \draw (-2,0.5) node {$1/2$};
 \node[root] (A3) at (0.5,3){};
 \node[authorised] (A2) at (0,2){};
 \node[authorised] (B2) at (1,2){};
 \node[authorised] (A1) at (0,1){};
 \node[authorised] (B1) at (1,1){};
 \node[mandatory] (B1/2) at (1,0.5){};
 \node[mandatory] (A1/2) at (0,0.5){};
 \node[empty] (A0) at (0,0){};
 \node[empty] (B0) at (1,0){};
 \draw (A3)--(A2);
 \draw (A3)--(B2);
 \draw (A2)--(A1)--(A1/2)--(A0);
 \draw (B2)--(B1)--(B1/2)--(B0);
\end{tikzpicture}
\end{center}

\begin{center}
\begin{tabular}{|c|c|c|c|}
\hline
 & Rank & Number of singularities & min. Katz ranks\\
 \hline
 Generic rep. & 4 & 1 & 2\\
 \hline
 Nongeneric reps. & 2 & 2 & 1, 1/2\\
 \hline
\end{tabular}
\end{center}

\subsubsection*{Doubly degenerate Painlev\'e~III}

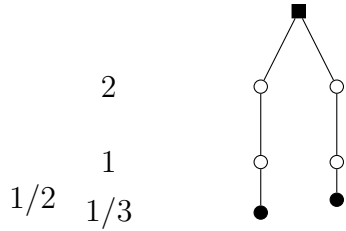
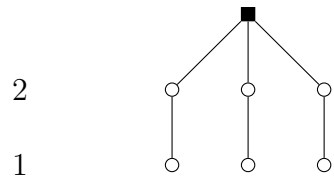
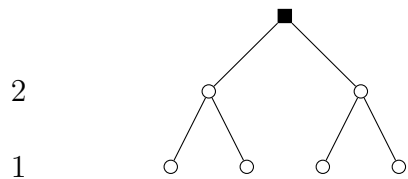
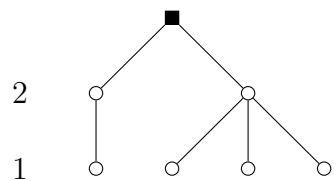
\begin{center}
\begin{tikzpicture}
\tikzstyle{mandatory}=[circle,fill=black,minimum size=5pt,draw, inner sep=0pt]
\tikzstyle{authorised}=[circle,fill=white,minimum size=5pt,draw,inner sep=0pt]
\tikzstyle{empty}=[circle,fill=black,minimum size=0pt,inner sep=0pt]
\tikzstyle{root}=[fill=black,minimum size=5pt,draw,inner sep=0pt]
\tikzstyle{indeterminate}=[circle,densely dotted,fill=white,minimum size=5pt,draw, inner sep=0pt]
 \draw (-2,1) node {$1$};
 \draw (-2,2) node {$2$};
 \draw (-3,0.5) node {$1/2$};
 \draw (-2,0.33) node {$1/3$};
 \node[root] (A3) at (0.5,3){};
 \node[authorised] (A2) at (0,2){};
 \node[authorised] (B2) at (1,2){};
 \node[authorised] (A1) at (0,1){};
 \node[authorised] (B1) at (1,1){};
 \node[mandatory] (A1/3) at (0,0.33){};
 \node[mandatory] (B1/2) at (1,0.5){};
 \node[empty] (A0) at (0,0){};
 \node[empty] (B0) at (1,0){};
 \draw (A3)--(A2);
 \draw (A3)--(B2);
 \draw (A2)--(A1)--(A1/3)--(A0);
 \draw (B2)--(B1)--(B1/2)--(B0);
\end{tikzpicture}
\end{center}

\begin{center}
\begin{tabular}{|c|c|c|c|}
\hline
 & Rank & Number of singularities & min. Katz ranks\\
 \hline
 Generic rep. & 5 & 1 & 2\\
 \hline
 1st nongeneric rep. & 2 & 2 & 1/2, 1/2\\
 \hline
 2nd nongeneric rep. & 3 & 2 & 1, 1/3\\
 \hline

\end{tabular}
\end{center}

\subsubsection*{Painlev\'e~IV}

\begin{center}
\begin{tikzpicture}
\tikzstyle{mandatory}=[circle,fill=black,minimum size=5pt,draw, inner sep=0pt]
\tikzstyle{authorised}=[circle,fill=white,minimum size=5pt,draw,inner sep=0pt]
\tikzstyle{empty}=[circle,fill=black,minimum size=0pt,inner sep=0pt]
\tikzstyle{root}=[fill=black,minimum size=5pt,draw,inner sep=0pt]
\tikzstyle{indeterminate}=[circle,densely dotted,fill=white,minimum size=5pt,draw, inner sep=0pt]
 \draw (-2,2) node {$2$};
 \draw (-2,1) node {$1$};
 \node[root] (A3) at (1,3){};
 \node[authorised] (A2) at (0,2){};
 \node[authorised] (B2) at (1,2){};
 \node[authorised] (C2) at (2,2){};
 \node[authorised] (A1) at (0,1){};
 \node[authorised] (B1) at (1,1){};
 \node[authorised] (C1) at (2,1){};
 \node[empty] (A0) at (0,0){};
 \node[empty] (B0) at (1,0){};
 \node[empty] (C0) at (2,0){};
 \draw (A3)--(A2);
 \draw (A3)--(B2);
 \draw (A3)--(C2);
 \draw (A2)--(A1)--(A0);
 \draw (B2)--(B1)--(B0);
 \draw (C2)--(C1)--(C0);

\end{tikzpicture}
\end{center}

\begin{center}
\begin{tabular}{|c|c|c|c|}
\hline
 & Rank & Number of singularities & min. Katz ranks\\
 \hline
 Generic rep. & 3 & 1 & 2\\
 \hline
 Nongeneric reps. & 2 & 2 & 2, 0\\
 \hline
\end{tabular}
\end{center}

\subsubsection*{Painlev\'e~V}

\begin{center}
\begin{tikzpicture}
\tikzstyle{mandatory}=[circle,fill=black,minimum size=5pt,draw, inner sep=0pt]
\tikzstyle{authorised}=[circle,fill=white,minimum size=5pt,draw,inner sep=0pt]
\tikzstyle{empty}=[circle,fill=black,minimum size=0pt,inner sep=0pt]
\tikzstyle{root}=[fill=black,minimum size=5pt,draw,inner sep=0pt]
\tikzstyle{indeterminate}=[circle,densely dotted,fill=white,minimum size=5pt,draw, inner sep=0pt]
 \draw (-2,1) node {$1$};
 \draw (-2,2) node {$2$};
 \node[root] (A3) at (1.5,3){};
 \node[authorised] (A2) at (0.5,2){};
 \node[authorised] (B2) at (2.5,2){};
 \node[authorised] (A1) at (0,1){};
 \node[authorised] (B1) at (1,1){};
 \node[authorised] (C1) at (2,1){};
 \node[authorised] (D1) at (3,1){};
 \node[empty] (A0) at (0,0){};
 \node[empty] (B0) at (1,0){};
 \node[empty] (C0) at (2,0){};
 \node[empty] (D0) at (3,0){};
 \draw (A2)--(A1)--(A0);
 \draw (A2)--(B1)--(B0);
 \draw (B2)--(C1)--(C0);
 \draw (B2)--(D1)--(D0);
 \draw (A3)--(A2);
 \draw (A3)--(B2);
\end{tikzpicture}
\end{center}

\begin{center}
\begin{tabular}{|c|c|c|c|}
\hline
 & Rank & Number of singularities & min. Katz ranks\\
 \hline
 Generic rep. & 4 & 1 & 2\\
 \hline
 Nongeneric reps. & 2 & 3 & 1, 0, 0\\
 \hline
\end{tabular}
\end{center}

\subsubsection*{Painlev\'e~VI}

\begin{center}
\begin{tikzpicture}
\tikzstyle{mandatory}=[circle,fill=black,minimum size=5pt,draw, inner sep=0pt]
\tikzstyle{authorised}=[circle,fill=white,minimum size=5pt,draw,inner sep=0pt]
\tikzstyle{empty}=[circle,fill=black,minimum size=0pt,inner sep=0pt]
\tikzstyle{root}=[fill=black,minimum size=5pt,draw,inner sep=0pt]
\tikzstyle{indeterminate}=[circle,densely dotted,fill=white,minimum size=5pt,draw, inner sep=0pt]
 \draw (-2,1) node {$1$};
 \draw (-2,2) node {$2$};
 \node[root] (A3) at (0,3){};
 \node[authorised] (A2) at (1,2){};
 \node[authorised] (B2) at (-1,2){};
 \node[authorised] (A1) at (0,1){};
 \node[authorised] (B1) at (1,1){};
 \node[authorised] (C1) at (2,1){};
 \node[authorised] (D1) at (-1,1){};
 \node[empty] (A0) at (0,0){};
 \node[empty] (B0) at (1,0){};
 \node[empty] (C0) at (2,0){};
 \node[empty] (D0) at (-1,0){};
 \draw (A2)--(A1)--(A0);
 \draw (A2)--(B1)--(B0);
 \draw (A2)--(C1)--(C0);
 \draw (B2)--(D1)--(D0);
 \draw (A3)--(A2);
 \draw (A3)--(B2);
 \draw (-1,-0.5) node {{\scriptsize $2$}};
\end{tikzpicture}
\end{center}

\begin{center}
\begin{tabular}{|c|c|c|c|}
\hline
 & Rank & Number of singularities & min. Katz ranks\\
 \hline
 Generic rep. & 5 & 1 & 2\\
 \hline
 1st nongeneric rep. & 3 & 2 & 1, 0\\
 \hline
 2nd nongeneric rep. & 2 & 4 & 0, 0, 0, 0\\
 \hline
\end{tabular}
\end{center}

\subsection{Painlev\'e I}

\subsubsection{Standard representation} The standard rank~2 Painlev\'e I Lax representation has a single irregular singularity given by a pole of order~4, but such that the leading term is nilpotent, which leads to a Turritin--Levelt normal form having a single twisted Stokes circle of slope~$5/2$.

{\samepage This corresponds to the following fission forest $\mathbf F^{\rm I}_{\rm std}$ (having a single tree):
\begin{center}
\begin{tikzpicture}
\tikzstyle{mandatory}=[circle,fill=black,minimum size=5pt,draw, inner sep=0pt]
\tikzstyle{authorised}=[circle,fill=white,minimum size=5pt,draw,inner sep=0pt]
\tikzstyle{empty}=[circle,fill=black,minimum size=0pt,inner sep=0pt]
\tikzstyle{root}=[fill=black,minimum size=5pt,draw,inner sep=0pt]
\tikzstyle{indeterminate}=[circle,densely dotted,fill=white,minimum size=5pt,draw, inner sep=0pt]
 \node[root] (R) at (0,3){};
 \node[mandatory] (A5) at (0,2.5){};
 \node[authorised] (A4) at (0,2){};
 \node[authorised] (A3) at (0,1.5){};
 \node[authorised] (A2) at (0,1){};
 \node[authorised] (A1) at (0,0.5){};
 \node[empty] (A0) at (0,0){};
 \draw (R)--(A5)--(A4)--(A3)--(A2)--(A1)--(A0);
 \draw (-1.3,2.5) node {$5/2$};
 \draw (-1.3,2) node {$2$};
 \draw (-1.3,1.5) node {$3/2$};
 \draw (-1.3,1) node {$1$};
 \draw (-1.3,0.5) node {$1/2$};

\end{tikzpicture}
\end{center}

}

The global irregular classes with fission forest $\mathbf F^{\rm I}_{\rm std}$ and minimal Poincar\'e--Katz rank (equal to $5/2$) are exactly those of the form
\[
\bm{{\Theta}}^{\rm I}_{\rm std}=\bigl\langle \alpha_{5/2} z_a^{-5/2}+\alpha_{2} z_a^{-2}+\alpha_{3/2} z_a^{-3/2}+\alpha_{1} z_a^{-1}+\alpha_{1/2} z_a^{-1/2}\bigr\rangle_a,
\]
with $a\in \mathbb P^1$, $\alpha_{5/2}, \alpha_2, \alpha_{3/2}, \alpha_1, \alpha_{1/2}\in \mathbb C$, and $\alpha_{5/2}\neq 0$.

Given an irregular connection on $\mathbb P^1$ with such an irregular class, the rank of its generic class of nearby representations is minimal when $a=\infty$ (to which we can always reduce up to applying a M\"obius transformation).

Let us thus consider a (modified) irregular class of the form
\[
\bm{\breve{\Theta}}^{\rm I}_{\rm std}=\bigl\langle\alpha_{5/2} z^{5/2}+\alpha_{2} z^{2}+\alpha_{3/2} z^{3/2}+\alpha_{1} z+\alpha_{1/2} z^{1/2}\bigr\rangle_\infty,
\]
with $\alpha_{5/2}\neq 0$, and determine the classes of its nearby representations.

\subsubsection{Generic representation}

The short fission tree $\mathbf T^{\rm I}$ of $\bm{\breve{\Theta}}^{\rm I}_{\rm std}$ is the following:
$$
\begin{tikzpicture}[scale=1.3]
\tikzstyle{mandatory}=[circle,fill=black,minimum size=5pt,draw, inner sep=0pt]
\tikzstyle{authorised}=[circle,fill=white,minimum size=5pt,draw,inner sep=0pt]
\tikzstyle{empty}=[circle,fill=black,minimum size=0pt,inner sep=0pt]
\tikzstyle{root}=[fill=black,minimum size=5pt,draw,inner sep=0pt]
\tikzstyle{indeterminate}=[circle,densely dotted,fill=white,minimum size=5pt,draw, inner sep=0pt]
 \node[root] (R) at (0,2){};
 \node[mandatory] (A5) at (0,1.67){};
 \node[authorised] (A4) at (0,1.33){};
 \node[authorised] (A3) at (0,1){};
 \node[authorised] (A2) at (0,0.67){};
 \node[authorised] (A1) at (0,0.33){};
 \node[empty] (A0) at (0,0){};
 \draw (R)--(A5)--(A4)--(A3)--(A2)--(A1)--(A0);
 %\draw (-1.3,2) node {$2$};
 \draw (-1.3,1.67) node {$5/3$};
 \draw (-1.3,1.33) node {$4/3$};
 \draw (-1.3,1) node {$1$};
 \draw (-1.3,0.67) node {$2/3$};
 \draw (-1.3,0.33) node {$1/3$};
\end{tikzpicture}
$$

An irregular class of generic form has fission tree $\mathbf T^I$ if and only if it is of the form
\[
\bm{\breve{\Theta}}^{\rm I}_{\rm gen}=\cir{\lambda' z^2+\alpha'_{5/3}+\alpha'_{4/3} z^{4/3}+ \alpha'_{1} z+ \alpha'_{2/3} z^{1/3}+ \alpha'_{1/3} z^{1/3}}_\infty,
\]
with $\alpha'_{5/3}\neq 0$.

\begin{Remark}
Notice that the coefficient $\lambda'$ of exponent 2 may vanish, so the minimal Poincar\'e--Katz rank for such an irregular class is $5/3$, hence following our conventions for truncations of fission trees, we draw the root of the tree at height $\big\lfloor \frac{5}{3}+1\big\rfloor=2$.
\end{Remark}

The generic class of representations thus has rank 3, a single Fourier sphere coefficient, so $N=N_1=N_1^+$, and the standard representation corresponds to the sole nongeneric class of representations.

Notice that the generic class of representations corresponds to the (abstract version of the more explicit) alternative Lax representation for Painlev\'e I found in \cite{joshi2009linearization}.

\subsubsection{Diagram} The corresponding diagram $\Gamma^I$ has one vertex with one loop \cite[Section~7]{doucot2021diagrams}, as drawn below

\begin{center}
\begin{tikzpicture}[scale=0.6]
\tikzstyle{vertex}=[circle,fill=cyan,minimum size=6pt,inner sep=0pt]
\tikzstyle{vertex_2}=[circle,fill=purple,minimum size=6pt,inner sep=0pt]
\node[vertex] (A) at (0,0){} ;
\draw (A) to[out=45, in=0] (0,1.3) to[out=180,in=135] (A);
\end{tikzpicture}
\end{center}
The two classes of representations correspond to different readings of the diagram.

\subsection{Painlev\'e~II}

\subsubsection{Standard representation} The standard rank 2 Painlev\'e~II Lax representation has one irregular untwisted singularity corresponding to a pole of order 4.

This corresponds to the following fission forest $\mathbf F^{\rm II}_{\rm std}$:
$$
\begin{tikzpicture}
\tikzstyle{mandatory}=[circle,fill=black,minimum size=5pt,draw, inner sep=0pt]
\tikzstyle{authorised}=[circle,fill=white,minimum size=5pt,draw,inner sep=0pt]
\tikzstyle{empty}=[circle,fill=black,minimum size=0pt,inner sep=0pt]
\tikzstyle{root}=[fill=black,minimum size=5pt,draw,inner sep=0pt]
\tikzstyle{indeterminate}=[circle,densely dotted,fill=white,minimum size=5pt,draw, inner sep=0pt]
 \draw (-2,1) node {$1$};
 \draw (-2,2) node {$2$};
 \draw (-2,3) node {$3$};
 \node[root] (A4) at (0.5,4){};
 \node[authorised] (A3) at (0,3){};
 \node[authorised] (B3) at (1,3){};
 \node[authorised] (A2) at (0,2){};
 \node[authorised] (B2) at (1,2){};
 \node[authorised] (A1) at (0,1){};
 \node[authorised] (B1) at (1,1){};
 \node[empty] (A0) at (0,0){};
 \node[empty] (B0) at (1,0){};
 \draw (A4)--(A3);
 \draw (A4)--(B3);
 \draw (A3)--(A2)--(A1)--(A0);
 \draw (B3)--(B2)--(B1)--(B0);
\end{tikzpicture}
$$

The irregular classes with fission forest $\mathbf F^{\rm II}_{\rm std}$ and minimal Poincar\'e--Katz rank (equal to 3) are all irregular classes of the form
\[
\bm{\Theta}_{\rm std}^{\rm II}=\cir{\alpha_3 z_a^{-3}+\alpha_2 z_a^{-2}+ \alpha_1 z_a^{-1}}_a+\cir{\beta_3 z_a^{-3}+\beta_2 z_a^{-2}+\beta_1 z_a^{-1}}_a,
\]
with $a\in \mathbb P^1$ and $\alpha_3\neq \beta_3$.

Given a connection on $\mathbb P^1$ with such an irregular class, the rank of its generic class of nearby representations is minimal when $a=\infty$ and (say) $\beta_3=0$ (to which we can always reduce up to applying a M\"obius transformation and tensoring with a rank one connection).

Let us thus consider a (modified) irregular class of the form
\[
\bm{\breve\Theta}^{\rm II}_{\rm std}=\bm{\Theta}^{\rm II}_{\rm std}=\cir{\alpha_3 z^{3}+\alpha_2 z^{2}+ \alpha_1 z}_\infty+\cir{\beta_2 z^{2}+\beta_1 z}_\infty,
\]
with $\alpha_3\neq 0$ and determine the classes of its nearby representations.

\subsubsection{Generic representation}

The short fission tree $\mathbf T^{\rm II}$ of $\bm{\breve\Theta}^{\rm II}_{\rm std}$ is
$$
\begin{tikzpicture}
\tikzstyle{mandatory}=[circle,fill=black,minimum size=5pt,draw, inner sep=0pt]
\tikzstyle{authorised}=[circle,fill=white,minimum size=5pt,draw,inner sep=0pt]
\tikzstyle{empty}=[circle,fill=black,minimum size=0pt,inner sep=0pt]
\tikzstyle{root}=[fill=black,minimum size=5pt,draw,inner sep=0pt]
\tikzstyle{indeterminate}=[circle,densely dotted,fill=white,minimum size=5pt,draw, inner sep=0pt]
 \draw (-2,1) node {$1$};
 \draw (-2,2) node {$2$};
 \draw (-2,1.5) node {$3/2$};
 \draw (-2,0.5) node {$1/2$};
 \node[root] (A3) at (0.5,3){};
 \node[authorised] (A2) at (0,2){};
 \node[authorised] (B2) at (1,2){};
 \node[authorised] (A1) at (0,1){};
 \node[authorised] (B1) at (1,1){};
 \node[authorised] (A1/2) at (0,0.5){};
 \node[mandatory] (A3/2) at (0,1.5){};
 \node[empty] (A0) at (0,0){};
 \node[empty] (B0) at (1,0){};
 \draw (A3)--(A2);
 \draw (A3)--(B2);
 \draw (A2)--(A3/2)--(A1)--(A1/2)--(A0);
 \draw (B2)--(B1)--(B0);

\end{tikzpicture}
$$

An irregular class of generic form has fission tree $\mathbf T^{\rm II}$ if and only if it is of the form
\[
\bm{\breve\Theta}^{\rm II}_{\rm gen}=\cir{\lambda' z^2+\alpha'_{3/2}z^{3/2}+\alpha'_1 z+\alpha'_{1/2}z^{1/2}}_\infty + \cir{\mu' z^2+ \beta'_1 z}_\infty,
\]
with $\lambda'\neq \mu'$ and $\alpha'_{3/2}\neq 0$.

The generic class of representations has thus rank 3, $k=2$ Fourier sphere coefficients, and the corresponding partition of the set of Stokes circles of $\bm{\breve\Theta}^{\rm II}_{\rm gen}$ is given by
\begin{gather*}
N_1=N_1^+=\big\{\cir{\lambda' z^2+\alpha'_{3/2}z^{3/2}+\alpha'_1 z+\alpha'_{1/2}z^{1/2}}_\infty\big\},
\\
N_2=N_2^-= \big\{\cir{\mu' z^2+ \beta'_1 z}_\infty\big\}.
\end{gather*}

\subsubsection{Nongeneric representations} The two principal subtrees of $\mathbf T^{\rm II}$ are not isomorphic, so the two nongeneric classes of representations have distinct forests. The standard representation corresponds to the nongeneric class of representations for $i=1$.

The second nongeneric class of representations has the following forest $\mathbf F^{\rm II}_2$:
$$
\begin{tikzpicture}
\tikzstyle{mandatory}=[circle,fill=black,minimum size=5pt,draw, inner sep=0pt]
\tikzstyle{authorised}=[circle,fill=white,minimum size=5pt,draw,inner sep=0pt]
\tikzstyle{empty}=[circle,fill=black,minimum size=0pt,inner sep=0pt]
\tikzstyle{root}=[fill=black,minimum size=5pt,draw,inner sep=0pt]
\tikzstyle{indeterminate}=[circle,densely dotted,fill=white,minimum size=5pt,draw, inner sep=0pt]
\node[root] (A2)at (0,2){};
\node[mandatory] (A3/2) at (0,1.5){};
\node[authorised] (A1) at(0,1){};
\node[authorised] (A1/2) at (0,0.5){};
\node[empty] (A0) at (0,0){};
\node[root] (B1) at (2,1){};
\node[empty] (B0) at (2,0){};
\draw (A2)--(A3/2)--(A1)--(A1/2)--(A0);
\draw (B1)--(B0);
\draw (-2,1) node {$1$};
\draw (-2,1.5) node {$3/2$};
\draw (-2,0.5) node {$1/2$};
\draw (2,-0.5) node {{\scriptsize 2}};
\end{tikzpicture}
$$
It is of rank 2, with one irregular singularity (with Poincar\'e--Katz rank $3/2$) and one regular singularity.

Explicitly, elements of the orbit $\SL_2(\mathbb C)\cdot \bm{\breve\Theta}^{\rm II}_{\rm gen}$ belonging to this class of representations are of the form
\[
\bm{\breve\Theta}^{\rm II}_{2}=\cir{\lambda'' z^2+\alpha''_{3/2} z^{3/2}+\alpha''_1 z + \alpha''_{1/2} z^{1/2}}_\infty + \cir{0}_{a''},
\]
with $\alpha''\neq 0$, and $a''\in \mathbb P^1\smallsetminus\{\infty\}$, hence the corresponding non-modified irregular class is
\[
\bm{\Theta}^{\rm II}_{2}=\cir{\lambda'' z^2+\alpha''_{3/2} z^{3/2}+\alpha''_1 z + \alpha''_{1/2} z^{1/2}}_\infty + 2\cir{0}_{a''}.
\]

We notice that it corresponds precisely the other known Lax representation for Painlev\'e~II due to Flaschka--Newell \cite{flaschka1980monodromy}, also sometimes referred to as degenerate Painlev\'e~IV.

\subsubsection{Diagram}

The corresponding diagram $\Gamma^{\rm II}$ is the affine $A_1$ Dynkin diagram \cite[Section~7]{doucot2021diagrams}
$$
\begin{tikzpicture}[scale=0.8]
\tikzstyle{vertex}=[circle,fill=cyan,minimum size=6pt,inner sep=0pt]
\tikzstyle{vertex_2}=[circle,fill=purple,minimum size=6pt,inner sep=0pt]
\node[vertex] (A) at (0,0){} ;
\node[vertex_2] (B) at (3,0){} ;
\draw[double distance=2 pt] (A)--(B);
\end{tikzpicture}
$$
The classes of representations correspond to different readings of the diagram, as indicated on the figure below
$$
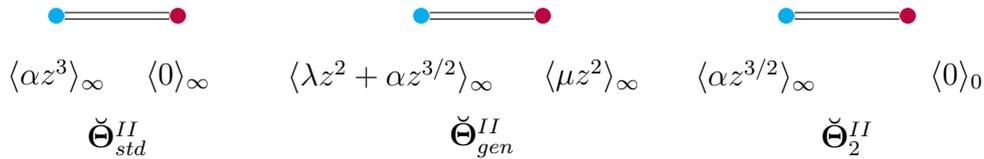
\begin{tikzpicture}[scale=0.758]
\tikzstyle{vertex}=[circle,fill=cyan,minimum size=6pt,inner sep=0pt]
\tikzstyle{vertex_2}=[circle,fill=purple,minimum size=6pt,inner sep=0pt]
\begin{scope}
\node[vertex] (A) at (0,0){} ;
\node[vertex_2] (B) at (3,0){} ;
\draw[double distance=2 pt] (A)--(B);
\draw (0,-1) node {{\scriptsize $\cir{\alpha_3 z^{3}+\cdots}_\infty$}};
\draw (3,-1) node {{\scriptsize $\cir{\beta_2 z^2+\cdots}_\infty$}};
\draw (1.5,-2) node {$\bm{\breve\Theta}^{\rm II}_{\rm std}$};
\end{scope}\hspace{-4mm}
\begin{scope}[xshift=8cm]
\node[vertex] (A) at (0,0){} ;
\node[vertex_2] (B) at (3,0){} ;
\draw[double distance=2 pt] (A)--(B);
\draw (-0.5,-1) node {{\scriptsize $\cir{\lambda' z^2+\alpha'_{3/2} z^{3/2}+\cdots}_\infty$}};
\draw (3.2,-1) node {{\scriptsize $\cir{\mu' z^2+\beta'_1 z}_\infty$}};
\draw (1.5,-2) node {$\bm{\breve\Theta}^{\rm II}_{\rm gen}$};
\end{scope}\hspace{-4mm}
\begin{scope}[xshift=16cm]
\node[vertex] (A) at (0,0){} ;
\node[vertex_2] (B) at (3,0){} ;
\draw[double distance=2 pt] (A)--(B);
\draw (-0.5,-1) node {{\scriptsize $\cir{\lambda'' z^2+\alpha''_{3/2} z^{3/2}+\cdots}_\infty$}};
\draw (3,-1) node {{\scriptsize $\cir{0}_{a''}$}};
\draw (1.5,-2) node {$\bm{\breve\Theta}^{\rm II}_{2}$};
\end{scope}
\end{tikzpicture}
$$

\subsection{Painlev\'e~III}
\label{subsec:painleve_III}

\subsubsection{Standard representation}

The standard rank 2 Painlev\'e~III Lax representation has two irregular singularities, both corresponding to a second order pole.

The corresponding fission forest $\mathbf F_{\rm std}^{\rm III(2)}$ is the following:\footnote{We will denote Painlev\'e~III, degenerate Painlev\'e~III and doubly degenerate Painlev\'e~III by $\rm III(2)$, $\rm III(1)$ and $\rm III(0)$, respectively, the numbers 2, 1, 0 referring to the number of parameters of each equation.}
$$
\begin{tikzpicture}
\tikzstyle{mandatory}=[circle,fill=black,minimum size=5pt,draw, inner sep=0pt]
\tikzstyle{authorised}=[circle,fill=white,minimum size=5pt,draw,inner sep=0pt]
\tikzstyle{empty}=[circle,fill=black,minimum size=0pt,inner sep=0pt]
\tikzstyle{root}=[fill=black,minimum size=5pt,draw,inner sep=0pt]
\tikzstyle{indeterminate}=[circle,densely dotted,fill=white,minimum size=5pt,draw, inner sep=0pt]
 \draw (-2,1) node {$1$};
 \node[root] (A2) at (0.5,2){};
 \node[root] (B2) at (2.5,2){};
 \node[authorised] (A1) at (0,1){};
 \node[authorised] (B1) at (1,1){};
 \node[authorised] (C1) at (2,1){};
 \node[authorised] (D1) at (3,1){};
 \node[empty] (A0) at (0,0){};
 \node[empty] (B0) at (1,0){};
 \node[empty] (C0) at (2,0){};
 \node[empty] (D0) at (3,0){};
 \draw (A2)--(A1)--(A0);
 \draw (A2)--(B1)--(B0);
 \draw (B2)--(C1)--(C0);
 \draw (B2)--(D1)--(D0);
\end{tikzpicture}
$$

The irregular classes with this fission forest and minimal Poincar\'e--Katz ranks (both equal to 1) are the irregular classes of the form
\[
\bm{\Theta}^{\rm III(2)}_{\rm std}=\cir{\alpha z_b^{-1}}_b+ \cir{\beta z_b^{-1}}_b+ \cir{\gamma z^{-1}_a}_a +\cir{\delta z_a^{-1}}_a,
\]
with $\alpha,\beta, \gamma,\delta\in \mathbb C$, $a, b\in \mathbb P^1$ distinct, $\alpha\neq \beta$, and $\gamma\neq \delta$.

If $(E,\nabla)$ is a connection on $\mathbb P^1$ of this type, the rank of its generic class of nearby representations is minimal when one of the singularities, say $b$, is at infinity, and, at the other singularity~$a$ at finite distance, and one of the Stokes circles at finite distance, say $\cir{\delta z_a^{-1}}_a$, is the tame circle, with trivial formal monodromy (we can always reduce to this case by applying a M\"obius transformation to send $b$ to infinity, and a twist by a rank one connection to set $\delta$ to zero and make the formal monodromy of $\cir{\delta z_a^{-1}}_a$ trivial).

This corresponds to having a modified irregular class of the form
\[
\bm{\breve\Theta}^{\rm III(2)}_{\rm std}=\cir{\alpha z}_\infty+ \cir{\beta z}_\infty+\cir{\gamma z_a^{-1}}_a.
\]
with $\alpha\neq \beta$, and $\gamma\neq 0$. Let us determine the classes of nearby representations of \smash{$\bm{\breve\Theta}^{\rm III(2)}_{\rm std}$}.

\subsubsection{Generic representation}

The short fission tree $\mathbf T^{\rm III(2)}$ of $\bm{\breve\Theta}^{\rm III(2)}_{\rm std}$ is the following:
$$
\begin{tikzpicture}
\tikzstyle{mandatory}=[circle,fill=black,minimum size=5pt,draw, inner sep=0pt]
\tikzstyle{authorised}=[circle,fill=white,minimum size=5pt,draw,inner sep=0pt]
\tikzstyle{empty}=[circle,fill=black,minimum size=0pt,inner sep=0pt]
\tikzstyle{root}=[fill=black,minimum size=5pt,draw,inner sep=0pt]
\tikzstyle{indeterminate}=[circle,densely dotted,fill=white,minimum size=5pt,draw, inner sep=0pt]
 \draw (-2,1) node {$1$};
 \draw (-2,2) node {$2$};
 \draw (-2,0.5) node {$1/2$};
 \node[root] (A3) at (1,3){};
 \node[authorised] (A2) at (0.5,2){};
 \node[authorised] (B2) at (2,2){};
 \node[authorised] (A1) at (0,1){};
 \node[authorised] (B1) at (1,1){};
 \node[authorised] (C1) at (2,1){};
 \node[mandatory] (C1/2) at (2,0.5){};
 \node[empty] (A0) at (0,0){};
 \node[empty] (B0) at (1,0){};
 \node[empty] (C0) at (2,0){};
 \draw (A3)--(A2);
 \draw (A3)--(B2);
 \draw (A2)--(A1)--(A0);
 \draw (A2)--(B1)--(B0);
 \draw (B2)--(C1)--(C1/2)--(C0);
\end{tikzpicture}
$$

An irregular class of generic form has fission tree $\mathbf T^{\rm III(2)}$ if and only if it is of the form
\[\bm{\breve{\Theta}}^{\rm III(2)}_{\rm gen}=\cir{\lambda' z^2+\alpha' z}_\infty+\cir{\lambda' z^2+\beta' z}_{\infty}+\cir{\mu' z^2 +\gamma'_1 z+ \gamma'_{1/2} z^{1/2}}_\infty,\]
with $\lambda'\neq \mu'$, $\alpha'\neq \beta'$, and $\gamma'_{1/2}\neq 0$.

The generic class of representations has thus rank 4, $k=2$ Fourier sphere coefficients, and the corresponding partition of the set of Stokes circles of \smash{$\bm{\breve{\Theta}}^{\rm III(2)}_{\rm gen}$} is given by
\begin{gather*}
N_1=N_1^-=\big\{ \cir{\lambda' z^2+\alpha' z}_\infty,\cir{\lambda' z^2+\beta' z}_{\infty}\big\},
\\
N_2=N_2^-=\big\{\cir{\mu' z^2 + \gamma'_1 z+\gamma'_{1/2} z^{1/2}}_\infty\big\}.
\end{gather*}

\subsubsection{Nongeneric representations}

Since the two principal subtrees are not isomorphic, the forests of the two corresponding nongeneric classes of representations do not coincide. The standard Lax representation is the nongeneric class of representations for $i=2$.

The other nongeneric class of representations, corresponding to $i=1$, has the following forest \smash{$\mathbf F^{\rm III(2)}_1$}, with one irregular singularity (with Poincar\'e--Katz rank 1/2), and two regular singularities:
$$
\begin{tikzpicture}
\tikzstyle{mandatory}=[circle,fill=black,minimum size=5pt,draw, inner sep=0pt]
\tikzstyle{authorised}=[circle,fill=white,minimum size=5pt,draw,inner sep=0pt]
\tikzstyle{empty}=[circle,fill=black,minimum size=0pt,inner sep=0pt]
\tikzstyle{root}=[fill=black,minimum size=5pt,draw,inner sep=0pt]
\tikzstyle{indeterminate}=[circle,densely dotted,fill=white,minimum size=5pt,draw, inner sep=0pt]
\tikzstyle{vertex}=[circle,fill=cyan,minimum size=6pt,inner sep=0pt]
\tikzstyle{vertex_2}=[circle,fill=purple,minimum size=6pt,inner sep=0pt]
 %\draw (-0.5,1) node {\scriptsize $1$};
 \draw (-0.5,0.5) node {\scriptsize $1/2$};
 \node[root] (A1) at (0,1){};
 \node[root] (B1) at (1,1){};
 \node[root] (C1) at (2,1){};
 \node[empty] (A0) at (0,0){};
 \node[empty] (B0) at (1,0){};
 \node[empty] (C0) at (2,0){};
 \node[mandatory] (C1/2) at (2,0.5){};
 \draw (A1)--(A0);
 \draw (B1)--(B0);
 \draw (C1)--(C1/2)--(C0);
\end{tikzpicture}
$$

Explicitly, elements of the orbit \smash{$\SL_2(\mathbb C)\cdot \bm{\breve\Theta}^{\rm III(2)}_{\rm gen}$} belonging to this class of representations are of the form
\[
\bm{\breve{\Theta}}^{\rm III(2)}_1=\cir{0}_{a''}+ \cir{0}_{b''}+\cir{\alpha'' z^{1/2}}_\infty
\]
with $a'', b''\in \mathbb P^1\smallsetminus \{\infty\}$, $a''\neq b''$, $\alpha''\neq 0$, hence the corresponding non-modified irregular class~is
\[
\bm{\Theta}^{\rm III(2)}_1=2\cir{0}_{a''}+ 2\cir{0}_{b''}+\cir{\alpha'' z^{1/2}}_\infty .
\]

Notice that this class of representations corresponds to the alternative Lax representation for Painlev\'e~III sometimes referred to as degenerate Painlev\'e~V, cf.\ \cite{joshi2007linearization}.

\begin{Remark}\label{rem:about_fig_intro}
We obtain in this way Figure \ref{fig:PIII_tree_and_diagram} of the introduction: the fission forest on the left of the figure is \smash{$\mathbf F_1^{\rm III(2)}$}, and the fission forest on the right of the figure is the forest associated to the \textit{modified} irregular class \smash{$\bm{\breve\Theta}^{\rm III(2)}_{\rm std}$}, i.e., compared to \smash{$\mathbf F_{\rm std}^{\rm III(2)}$} one full branch has been removed.
\end{Remark}

\subsubsection{Diagram}

The corresponding diagram $\Gamma^{\rm III(2)}$ is the following \cite[Section~7]{doucot2021diagrams}:
$$
\begin{tikzpicture}[scale=0.8]
\tikzstyle{vertex}=[circle,fill=cyan,minimum size=6pt,inner sep=0pt]
\tikzstyle{vertex_2}=[circle,fill=purple,minimum size=6pt,inner sep=0pt]
\node[vertex] (A) at (0,0){};
\node[vertex_2] (B) at (2,0){};
\node[vertex] (C) at (4,0){};
\draw[double distance=2 pt] (A)--(B)--(C);
\draw[dashed] (B) to[out=45, in=0] (2,1.3) to[out=180,in=135] (B);
\end{tikzpicture}
$$
The classes of representations correspond to different readings, as represented on the figure below
$$
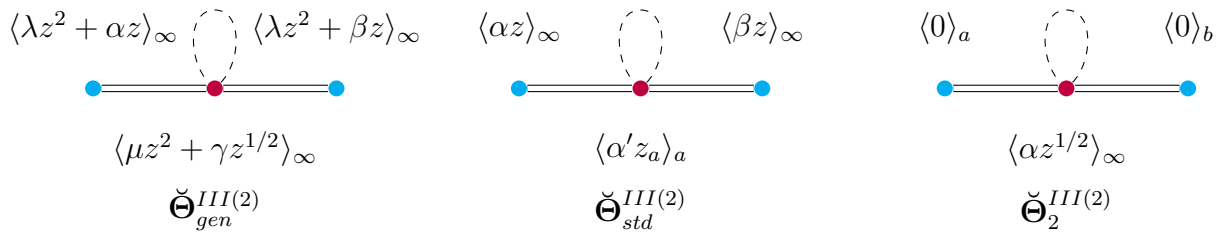
\begin{tikzpicture}[scale=0.78]
\tikzstyle{vertex}=[circle,fill=cyan,minimum size=6pt,inner sep=0pt]
\tikzstyle{vertex_2}=[circle,fill=purple,minimum size=6pt,inner sep=0pt]
\begin{scope}
\node[vertex] (A) at (0,0){};
\node[vertex_2] (B) at (2,0){};
\node[vertex] (C) at (4,0){};
\draw[double distance=2 pt] (A)--(B)--(C);
\draw[dashed] (B) to[out=45, in=0] (2,1.3) to[out=180,in=135] (B);
\draw (0,1) node {{\scriptsize $\cir{\lambda' z^2+\alpha' z}_\infty$}};
\draw (2,-1) node {{\scriptsize $\cir{\mu' z^2+\gamma'_1 z+\gamma'_{1/2} z^{1/2}}_\infty$}};
\draw (4,1) node {{\scriptsize $\cir{\lambda' z^2+\beta' z}_\infty$}};
\draw (2,-2) node {$\bm{\breve{\Theta}}^{\rm III(2)}_{\rm gen}$};
\end{scope}
\begin{scope}[xshift=7cm]
\tikzstyle{vertex_nofill}=[draw,circle,minimum size=8pt,inner sep=0pt]
\node[vertex] (A) at (0,0){};
\node[vertex_2] (B) at (2,0){};
\node[vertex] (C) at (4,0){};
\draw[double distance=2 pt] (A)--(B)--(C);
\draw[dashed] (B) to[out=45, in=0] (2,1.3) to[out=180,in=135] (B);
\draw (0,1) node {$\cir{\alpha z}_\infty$};
\draw (2,-1) node {$\cir{\gamma z_a^{-1}}_a$};
\draw (4,1) node {$\cir{\beta z}_\infty$};
\draw (2,-2) node {$\bm{\breve{\Theta}}^{\rm III(2)}_{\rm std}$};
\end{scope}
\begin{scope}[xshift=14cm]
\tikzstyle{vertex_nofill}=[draw,circle,minimum size=8pt,inner sep=0pt]
\node[vertex] (A) at (0,0){};
\node[vertex_2] (B) at (2,0){};
\node[vertex] (C) at (4,0){};
\draw[double distance=2 pt] (A)--(B)--(C);
\draw[dashed] (B) to[out=45, in=0] (2,1.3) to[out=180,in=135] (B);
\draw (0,1) node {$\cir{0}_{a''}$};
\draw (2,-1) node {$\cir{\alpha'' z^{1/2}}_\infty$};
\draw (4,1) node {$\cir{0}_{b''}$};
\draw (2,-2) node {$\bm{\breve{\Theta}}^{\rm III(2)}_{1}$};
\end{scope}
\end{tikzpicture}
$$

\subsection{Degenerate Painlev\'e~III}

\subsubsection{Standard representation}

The standard Lax representation for degenerate Painlev\'e~III, as for nondegenerate Painlev\'e~III, has two irregular singularities, both corresponding to poles of order 2, but one of them with a~nilpotent leading term giving rise in the Turritin--Levelt normal form to a~single Stokes circle of slope $1/2$ \cite{ohyama2006studiesV}.

This corresponds to the following fission forest \smash{$\mathbf F^{\rm III(1)}_{\rm std}$}:
$$
\begin{tikzpicture}
\tikzstyle{mandatory}=[circle,fill=black,minimum size=5pt,draw, inner sep=0pt]
\tikzstyle{authorised}=[circle,fill=white,minimum size=5pt,draw,inner sep=0pt]
\tikzstyle{empty}=[circle,fill=black,minimum size=0pt,inner sep=0pt]
\tikzstyle{root}=[fill=black,minimum size=5pt,draw,inner sep=0pt]
\tikzstyle{indeterminate}=[circle,densely dotted,fill=white,minimum size=5pt,draw, inner sep=0pt]
 \draw (-2,1) node {$1$};
 \draw (-2,0.5) node {$1/2$};
 \node[root] (A2) at (0.5,2){};
 \node[authorised] (A1) at (0,1){};
 \node[authorised] (B1) at (1,1){};
 \node[root] (C1) at (2,1){};
 \node[mandatory] (C1/2) at (2,0.5){};
 \node[empty] (A0) at (0,0){};
 \node[empty] (B0) at (1,0){};
 \node[empty] (C0) at (2,0){};
 \draw (A2)--(A1)--(A0);
 \draw (A2)--(B1)--(B0);
 \draw (C1)--(C1/2)--(C0);
\end{tikzpicture}
$$

The irregular classes with fission forest $\mathbf F^{\rm III(1)}_{\rm std}$ and minimal Poincar\'e--Katz ranks (equal to 1 and 1/2) are those of the form
\[
\bm{\Theta}^{\rm III(1)}_{\rm std}=\cir{\alpha z_a^{-1}}_a +\cir{\beta z_a^{-1}}_a+\cir{\gamma z_b^{-1/2}}_b
\]
with $a,b\in \mathbb P^1$, $a\neq b$, $\gamma\neq 0$, and $\alpha\neq \beta$.

Given a connection of this type, the rank of its generic class of nearby representations is minimal if $b=\infty$ and one of the Stokes circles at $a$, say $\cir{\beta z_a^{-1}}_a$ is the tame circle, with trivial formal monodromy (we can always reduce to this case by applying a M\"obius transformation to have $b=\infty$, and a twist by a rank one connection to set $\beta=0$ and make the formal monodromy of $\cir{\beta z_a^{-1}}_a$ trivial).
This corresponds to having a modified irregular class of the form
\[
\bm{\breve\Theta}^{\rm III(1)}_{\rm std}=\cir{\alpha z_a^{-1}}_a+\cir{\gamma z^{1/2}}_\infty
\]
with $\alpha\neq 0$, $\gamma\neq 0$, $a\neq\infty$. Let us thus determine the classes of nearby representations of~\smash{$\bm{\breve\Theta}^{\rm III(1)}_{\rm std}$}.

\subsubsection{Generic representation}

The short fission tree $\mathbf T^{\rm III(1)}$ of \smash{$\bm{\breve\Theta}^{\rm III(1)}_{\rm std}$} is the following:
$$
\begin{tikzpicture}
\tikzstyle{mandatory}=[circle,fill=black,minimum size=5pt,draw, inner sep=0pt]
\tikzstyle{authorised}=[circle,fill=white,minimum size=5pt,draw,inner sep=0pt]
\tikzstyle{empty}=[circle,fill=black,minimum size=0pt,inner sep=0pt]
\tikzstyle{root}=[fill=black,minimum size=5pt,draw,inner sep=0pt]
\tikzstyle{indeterminate}=[circle,densely dotted,fill=white,minimum size=5pt,draw, inner sep=0pt]
 \draw (-2,1) node {$1$};
 \draw (-2,2) node {$2$};
 \draw (-2,0.5) node {$1/2$};
 \node[root] (A3) at (0.5,3){};
 \node[authorised] (A2) at (0,2){};
 \node[authorised] (B2) at (1,2){};
 \node[authorised] (A1) at (0,1){};
 \node[authorised] (B1) at (1,1){};
 \node[mandatory] (B1/2) at (1,0.5){};
 \node[mandatory] (A1/2) at (0,0.5){};
 \node[empty] (A0) at (0,0){};
 \node[empty] (B0) at (1,0){};
 \draw (A3)--(A2);
 \draw (A3)--(B2);
 \draw (A2)--(A1)--(A1/2)--(A0);
 \draw (B2)--(B1)--(B1/2)--(B0);
\end{tikzpicture}
$$
An irregular class of generic form has fission tree $\mathbf T^{\rm III(1)}$ if and only if it is of the form
\[
\bm{\breve\Theta}^{\rm III(1)}_{\rm gen}=\cir{\lambda' z^2+ \alpha'_1 z+\alpha'_{1/2} z^{1/2}}_\infty+\cir{\mu' z^2+\beta'_1 z+\beta'_{1/2} z^{1/2}}_\infty
\]
with $\lambda'\neq \mu'$, $\alpha'_{1/2}\neq 0$ and $\beta'_{1/2}\neq 0$.

The generic class of representations has thus rank 4, $k=2$ Fourier sphere coefficients, and the corresponding partition of the set of Stokes circles of \smash{$\bm{\breve\Theta}^{\rm III(1)}_{\rm gen}$} is given by
\begin{gather*}
 N_1=N_1^-=\big\{\cir{\lambda' z^2+ \alpha'_1 z+\alpha'_{1/2} z^{1/2}}_\infty\big\},
\\
N_2=N_2^-=\big\{\cir{\mu' z^2+\beta'_1 z+\beta'_{1/2} z^{1/2}}_\infty\big\}.
\end{gather*}

Furthermore, the two principal subtrees are isomorphic, so the two nongeneric classes of representations have the same forest, and there are only two distinct forests, the generic one and the nongeneric one, which coincides with the standard Lax representation.

\subsubsection{Diagram}

The corresponding diagram $\Gamma^{\rm III(1)}$ is the following \cite[Section~7]{doucot2021diagrams}:
$$
\begin{tikzpicture}[scale=0.8]
\tikzstyle{vertex}=[circle,fill=cyan,minimum size=6pt,inner sep=0pt]
\tikzstyle{vertex_2}=[circle,fill=purple,minimum size=6pt,inner sep=0pt]
\node[vertex] (A) at (-2,0){} ;
\node[vertex_2] (C) at (2,0){};
\draw (C)-- node[midway,above]{$4$} (A);
\draw[dashed] (A) to[out=-135, in=-90] (-3.4,0) to[out=90,in=135] (A);
\draw[dashed] (C) to[out=-45, in=-90] (3.4,0) to[out=90,in=45] (C);
\draw (4,0) node {$-1$};
\draw (-4,0) node {$-1$};
\end{tikzpicture}
$$
The classes of representations correspond to different readings of the diagram, as indicated on the figure below
$$
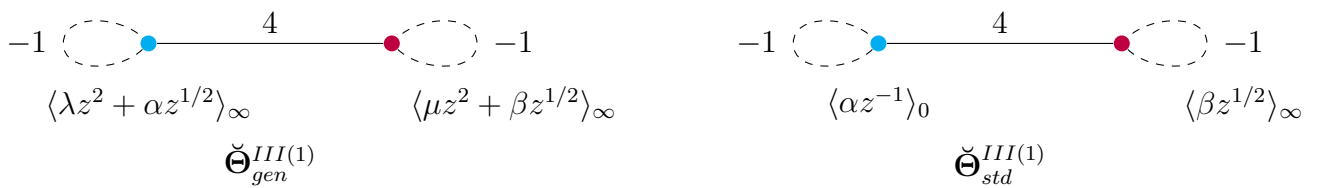
\begin{tikzpicture}[scale=0.712]
\tikzstyle{vertex}=[circle,fill=cyan,minimum size=6pt,inner sep=0pt]
\tikzstyle{vertex_2}=[circle,fill=purple,minimum size=6pt,inner sep=0pt]
\begin{scope}
\node[vertex] (A) at (-2,0){} ;
\node[vertex_2] (C) at (2,0){};
\draw (C)-- node[midway,above]{$4$} (A);
\draw[dashed] (A) to[out=-135, in=-90] (-3.4,0) to[out=90,in=135] (A);
\draw[dashed] (C) to[out=-45, in=-90] (3.4,0) to[out=90,in=45] (C);
\draw (4,0) node {$-1$};
\draw (-4,0) node {$-1$};
\draw (-2,-1) node {$\cir{\lambda' z^2+ \alpha'_1 z+\alpha'_{1/2} z^{1/2}}_\infty$};
\draw (4,-1) node {$\cir{\mu' z^2+\beta'_1 z+\beta'_{1/2} z^{1/2}}_\infty$};
\draw (0,-2) node {$\bm{\breve\Theta}^{\rm III(1)}_{\rm gen}$};
\end{scope}\hspace{-2mm}
\begin{scope}[xshift=12cm]
\node[vertex] (A) at (-2,0){} ;
\node[vertex_2] (C) at (2,0){};
\draw (C)-- node[midway,above]{$4$} (A);
\draw[dashed] (A) to[out=-135, in=-90] (-3.4,0) to[out=90,in=135] (A);
\draw[dashed] (C) to[out=-45, in=-90] (3.4,0) to[out=90,in=45] (C);
\draw (4,0) node {$-1$};
\draw (-4,0) node {$-1$};
\draw (-2,-1) node {$\cir{\alpha z_a^{-1}}_a$};
\draw (4,-1) node {$\cir{\gamma z^{1/2}}_\infty$};
\draw (0,-2) node {$\bm{\breve\Theta}^{\rm III(1)}_{\rm std}$};;
\end{scope}
\end{tikzpicture}
$$

\subsection{Doubly degenerate Painlev\'e~III}

\subsubsection{Standard representation}

The standard Lax representation for the doubly degenerate Painlev\'e~III equation has two irregular singularities, both corresponding to a second order pole, with a nilpotent leading term at each pole so that in the Turritin--Levelt decomposition, there is a single twisted Stokes circle of slope 1/2 \cite{ohyama2006studiesV}.

This corresponds to the following fission forest \smash{$\mathbf F^{\rm III(0)}_{\rm std}$}:
$$
\begin{tikzpicture}
\tikzstyle{mandatory}=[circle,fill=black,minimum size=5pt,draw, inner sep=0pt]
\tikzstyle{authorised}=[circle,fill=white,minimum size=5pt,draw,inner sep=0pt]
\tikzstyle{empty}=[circle,fill=black,minimum size=0pt,inner sep=0pt]
\tikzstyle{root}=[fill=black,minimum size=5pt,draw,inner sep=0pt]
\tikzstyle{indeterminate}=[circle,densely dotted,fill=white,minimum size=5pt,draw, inner sep=0pt]
 \draw (-2,0.5) node {$1/2$};
 \node[root] (A2) at (0,1){};
 \node[root] (B2) at (2,1){};
 \node[mandatory] (A1) at (0,0.5){};
 \node[mandatory] (B1) at (2,0.5){};
 \node[empty] (A0) at (0,0){};
 \node[empty] (B0) at (2,0){};
 \draw (A2)--(A1)--(A0);
 \draw (B2)--(B1)--(B0);
\end{tikzpicture}
$$

The irregular classes with fission forest \smash{$\mathbf F^{\rm III(0)}_{\rm std}$} and minimal Poincar\'e--Katz ranks (both equal to 1/2) are those of the form
\[
\bm{\Theta}^{\rm III(0)}_{\rm std}=\cir{\alpha z_a^{-1/2}}_a +\cir{\beta z_b^{-1/2}}_b,
\]
with $a, b\in \mathbb P^1$, $a\neq b$, $\alpha\neq 0$, $\beta\neq 0$.

Given a connection of this type, the rank of its generic class of nearby representations is minimal if one of the singularities, say $\beta$, is at infinity (to which we can reduce by applying a~M\"obius transformation).

Let us thus consider a (modified) irregular class of the form
\[
\bm{\breve\Theta}^{\rm III(0)}_{\rm std}=\cir{\alpha z_a^{-1/2}}_a +\cir{\beta z^{1/2}}_\infty,
\]
and determine its classes of nearby representations.

\subsubsection{Generic representation}

The short fission tree $\mathbf T^{\rm III(0)}$ of \smash{$\bm{\breve\Theta}^{\rm III(0)}_{\rm std}$} is the following:
$$
\begin{tikzpicture}
\tikzstyle{mandatory}=[circle,fill=black,minimum size=5pt,draw, inner sep=0pt]
\tikzstyle{authorised}=[circle,fill=white,minimum size=5pt,draw,inner sep=0pt]
\tikzstyle{empty}=[circle,fill=black,minimum size=0pt,inner sep=0pt]
\tikzstyle{root}=[fill=black,minimum size=5pt,draw,inner sep=0pt]
\tikzstyle{indeterminate}=[circle,densely dotted,fill=white,minimum size=5pt,draw, inner sep=0pt]
 \draw (-2,1) node {$1$};
 \draw (-2,2) node {$2$};
 \draw (-3,0.5) node {$1/2$};
 \draw (-2,0.33) node {$1/3$};
 \node[root] (A3) at (0.5,3){};
 \node[authorised] (A2) at (0,2){};
 \node[authorised] (B2) at (1,2){};
 \node[authorised] (A1) at (0,1){};
 \node[authorised] (B1) at (1,1){};
 \node[mandatory] (A1/3) at (0,0.33){};
 \node[mandatory] (B1/2) at (1,0.5){};
 \node[empty] (A0) at (0,0){};
 \node[empty] (B0) at (1,0){};
 \draw (A3)--(A2);
 \draw (A3)--(B2);
 \draw (A2)--(A1)--(A1/3)--(A0);
 \draw (B2)--(B1)--(B1/2)--(B0);

\end{tikzpicture}
$$
An irregular class of generic form has fission tree $\mathbf T^{\rm III(0)}$ if and only if it is of the form
\[
\bm{\breve\Theta}^{\rm III(0)}_{\rm gen}=\cir{\lambda' z^2+\alpha'_1 z+\alpha'_{1/3} z^{1/3}}_\infty +\cir{\mu' z^2+\beta'_1 z+\beta'_{1/2} z^{1/2}}_\infty,
\]
with $\lambda'\neq \mu'$, $\alpha'_{1/3}\neq 0$, $\beta'_{1/2}\neq 0$.

The generic representation has thus rank 5, $k=2$ Fourier sphere coefficients, and the corresponding partition of the set of Stokes circles of \smash{$\bm{\breve\Theta}^{\rm III(0)}_{\rm gen}$} is given by
\begin{gather*}
N_1=N_1^{-}=\big\{\cir{\lambda' z^2+\alpha'_1 z+\alpha'_{1/3} z^{1/3}}_\infty\big\},
\qquad
N_2=N_2^-=\big\{\cir{\mu' z^2+\beta'_1 z+\beta'_{1/2} z^{1/2}}_\infty\big\}.
\end{gather*}

\subsubsection{Nongeneric representations}

 Since the principal subtrees are not isomorphic, the two nongeneric forests are distinct. The standard representation corresponds to the nongeneric class of representations with $i=1$.

The other nongeneric representation, for $i=2$, has the following fission forest \smash{$\mathbf F^{\rm III(0)}_2$}:
$$
\begin{tikzpicture}
\tikzstyle{mandatory}=[circle,fill=black,minimum size=5pt,draw, inner sep=0pt]
\tikzstyle{authorised}=[circle,fill=white,minimum size=5pt,draw,inner sep=0pt]
\tikzstyle{empty}=[circle,fill=black,minimum size=0pt,inner sep=0pt]
\tikzstyle{root}=[fill=black,minimum size=5pt,draw,inner sep=0pt]
\tikzstyle{indeterminate}=[circle,densely dotted,fill=white,minimum size=5pt,draw, inner sep=0pt]
 \draw (-3,0.33) node {$1/3$};
 \draw (-3,1) node {$1$};
 \node[root] (R) at (0.5, 2){};
 \node[authorised] (A1) at (0,1){};
 \node[authorised] (B1) at (1,1){};
 \node[root] (C1) at (-2, 1){};
 \node[mandatory] (C1/3) at (-2, 0.33){};
 \node[empty] (A0) at (0,0){};
 \node[empty] (B0) at (1,0){};
 \node[empty] (C0) at (-2,0){};
 \draw (R)--(A1)--(A0);
 \draw (R)--(B1)--(B0);
 \draw (C1)--(C0);
 \draw (1,-0.5) node {{\scriptsize 2}};
\end{tikzpicture}
$$

Explicitly, elements of the orbit \smash{$\SL_2(\mathbb C)\cdot \bm{\breve\Theta}^{\rm III(0)}_{\rm gen}$} belonging to this class of representations are of the form
\[
\bm{\breve\Theta}^{\rm III(0)}_{2}=\cir{\beta'' z_{b''}^{-1}}_{b''}+\cir{\alpha'' z^{1/3}}_\infty
\]
with $b''\in \mathbb P^1\smallsetminus\{\infty\}$, $\alpha''\neq 0$ and $\beta''\neq 0$, hence the corresponding non-modified irregular class is
\[
\bm{\Theta}^{\rm III(0)}_{2}=\cir{\beta'' z_{b''}^{-1}}_{b''}+2\cir{0}_{b''}+\cir{\alpha'' z^{1/3}}_\infty,
\]
and its fission forest is indeed \smash{$\mathbf F^{\rm III(0)}_2$}.

The first nongeneric class of representations has thus rank~3, and minimal Poincar\'e--Katz ranks $1$, $1/3$.

\subsubsection{Diagram}

The corresponding diagram \smash{$\Gamma^{\rm III(0)}$} is the following \cite[Section~7]{doucot2021diagrams}:
$$
\begin{tikzpicture}[scale=0.8]
\tikzstyle{vertex}=[circle,fill=cyan,minimum size=6pt,inner sep=0pt]
\tikzstyle{vertex_2}=[circle,fill=purple,minimum size=6pt,inner sep=0pt]
\node[vertex] (A) at (-2,0){} ;
\node[vertex_2] (C) at (2,0){};
\draw (C)-- node[midway,above]{$6$} (A);
\draw[dashed] (A) to[out=-135, in=-90] (-3.4,0) to[out=90,in=135] (A);
\draw[dashed] (C) to[out=-45, in=-90] (3.4,0) to[out=90,in=45] (C);
\draw (4,0) node {$-1$};
\draw (-4,0) node {$-3$};
\end{tikzpicture}
$$
The three classes of representations correspond to different readings of the diagram, as indicated on the figure below
$$
\begin{tikzpicture}[scale=0.74]
\tikzstyle{vertex}=[circle,fill=cyan,minimum size=6pt,inner sep=0pt]
\tikzstyle{vertex_2}=[circle,fill=purple,minimum size=6pt,inner sep=0pt]
\begin{scope}
\node[vertex] (A) at (-2,0){} ;
\node[vertex_2] (C) at (2,0){};
\draw (C)-- node[midway,above]{$6$} (A);
\draw[dashed] (A) to[out=-135, in=-90] (-3.4,0) to[out=90,in=135] (A);
\draw[dashed] (C) to[out=-45, in=-90] (3.4,0) to[out=90,in=45] (C);
\draw (4,0) node {$-1$};
\draw (-4,0) node {$-3$};
\draw (-3,-1) node {{\scriptsize $\cir{\lambda' z^2+\alpha'_1 z+\alpha'_{1/3} z^{1/3}}_\infty$}};
\draw (3,-1) node {{\scriptsize $\cir{\mu' z^2+\beta'_1 z+\beta'_{1/2} z^{1/2}}_\infty$}};
\draw (0,-2) node {$\bm{\breve\Theta}^{\rm III(0)}_{\rm gen}$};
\end{scope}\hspace{-8mm}
\begin{scope}[xshift=11.5cm]
\node[vertex] (A) at (-2,0){} ;
\node[vertex_2] (C) at (2,0){};
\draw (C)-- node[midway,above]{$6$} (A);
\draw[dashed] (A) to[out=-135, in=-90] (-3.4,0) to[out=90,in=135] (A);
\draw[dashed] (C) to[out=-45, in=-90] (3.4,0) to[out=90,in=45] (C);
\draw (4,0) node {$-1$};
\draw (-4,0) node {$-3$};
\draw (-2,-1) node {$\cir{\alpha'' z^{1/3}}_\infty$};
\draw (2,-1) node {$\cir{\beta'' z^{-1}}_{b''}$};
\draw (0,-2) node {$\bm{\breve\Theta}^{\rm III(0)}_{2}$};
\end{scope}
\begin{scope}[xshift=7cm, yshift=-4cm]
\node[vertex] (A) at (-2,0){} ;
\node[vertex_2] (C) at (2,0){};
\draw (C)-- node[midway,above]{$6$} (A);
\draw[dashed] (A) to[out=-135, in=-90] (-3.4,0) to[out=90,in=135] (A);
\draw[dashed] (C) to[out=-45, in=-90] (3.4,0) to[out=90,in=45] (C);
\draw (4,0) node {$-1$};
\draw (-4,0) node {$-3$};
\draw (-2,-1) node {$\cir{\alpha z_a^{-1/2}}_a$};
\draw (2,-1) node {$\cir{\beta z^{1/2}}_\infty$};
\draw (0,-2) node {$\bm{\breve\Theta}^{\rm III(0)}_{\rm std}$};
\end{scope}
\end{tikzpicture}
$$

\subsection{Painlev\'e~IV}

\subsubsection{Standard representation}

The standard rank 2 Painlev\'e~IV Lax representation has an untwisted irregular singularity corresponding to a pole of order 3, and a regular singularity.

This corresponds to the following fission forest $\mathbf F^{\rm IV}_{\rm std}$:
$$
\begin{tikzpicture}
\tikzstyle{mandatory}=[circle,fill=black,minimum size=5pt,draw, inner sep=0pt]
\tikzstyle{authorised}=[circle,fill=white,minimum size=5pt,draw,inner sep=0pt]
\tikzstyle{empty}=[circle,fill=black,minimum size=0pt,inner sep=0pt]
\tikzstyle{root}=[fill=black,minimum size=5pt,draw,inner sep=0pt]
\tikzstyle{indeterminate}=[circle,densely dotted,fill=white,minimum size=5pt,draw, inner sep=0pt]
 \draw (-2,1) node {$1$};
 \draw (-2,2) node {$2$};
 \node[root] (A1) at (0,1){};
 \node[root] (B3) at (2.5,3){};
 \node[authorised] (B2) at (2,2){};
 \node[authorised] (C2) at (3,2){};
 \node[authorised] (B1) at (2,1){};
 \node[authorised] (C1) at (3,1){};
 \node[empty] (A0) at (0,0){};
 \node[empty] (B0) at (2,0){};
 \node[empty] (C0) at (3,0){};
 \draw (A1)--(A0);
 \draw (B3)--(B2)--(B1)--(B0);
 \draw (B3)--(C2)--(C1)--(C0);
 \draw (0,-0.5) node {{\scriptsize 2}};
\end{tikzpicture}
$$
The irregular classes with fission forest $\mathbf F^{\rm IV}_{\rm std}$ and minimal Poincar\'e--Katz ranks are those of the form
\[
\bm{{\Theta}}^{\rm IV}_{\rm std}=2\cir{0}_a+\cir{\alpha_2 z_b^{-2}+\alpha_1 z_b^{-1}}_b + \cir{\beta_2 z_b^{-2}+\beta_1 z_b^{-1}}_b
\]
with $a,b\in \mathbb P^1$, $a\neq b$, $\alpha_2, \alpha_1, \beta_2, \beta_1\in\mathbb C$, and $\alpha_2\neq \beta_2$.

Given a connection of this type, its generic class of nearby representations is of minimal rank if $b=\infty$, and its formal monodromy for $\cir{0}_b$ admits 1 as an eigenvalue (we can always reduce to this situation by applying a M\"obius transformation and a twist). This corresponds to having a modified irregular class of the form
\[
\bm{\breve{\Theta}}^{\rm IV}_{\rm std}=\cir{0}_a+\cir{\alpha_2 z^2+\alpha_1 z}_\infty + \cir{\beta_2 z^2+ \beta_1 z }_\infty
\]
with $a\in \mathbb P^1\smallsetminus \{\infty\}$, and $\alpha_2 \neq \beta_2$. Let us thus determine the classes of nearby representations of~$\bm{\breve{\Theta}}^{\rm IV}_{\rm std}$.

\subsubsection{Generic representation}

The short fission tree $\mathbf T^{\rm IV}$ of $\bm{\breve{\Theta}}^{\rm IV}_{\rm std}$ is the following:
$$
\begin{tikzpicture}
\tikzstyle{mandatory}=[circle,fill=black,minimum size=5pt,draw, inner sep=0pt]
\tikzstyle{authorised}=[circle,fill=white,minimum size=5pt,draw,inner sep=0pt]
\tikzstyle{empty}=[circle,fill=black,minimum size=0pt,inner sep=0pt]
\tikzstyle{root}=[fill=black,minimum size=5pt,draw,inner sep=0pt]
\tikzstyle{indeterminate}=[circle,densely dotted,fill=white,minimum size=5pt,draw, inner sep=0pt]
 \draw (-2,2) node {$2$};
 \draw (-2,1) node {$1$};
 \node[root] (A3) at (1,3){};
 \node[authorised] (A2) at (0,2){};
 \node[authorised] (B2) at (1,2){};
 \node[authorised] (C2) at (2,2){};
 \node[authorised] (A1) at (0,1){};
 \node[authorised] (B1) at (1,1){};
 \node[authorised] (C1) at (2,1){};
 \node[empty] (A0) at (0,0){};
 \node[empty] (B0) at (1,0){};
 \node[empty] (C0) at (2,0){};
 \draw (A3)--(A2);
 \draw (A3)--(B2);
 \draw (A3)--(C2);
 \draw (A2)--(A1)--(A0);
 \draw (B2)--(B1)--(B0);
 \draw (C2)--(C1)--(C0);
\end{tikzpicture}
$$
An irregular class of generic form has fission tree $\mathbf T^{\rm IV}$ if and only if it is of the form
\[
\bm{\breve{\Theta}}^{\rm IV}_{\rm gen}=\cir{\lambda' z^2+\alpha' z}_\infty + \cir{\mu' z^2 +\beta' z}_\infty +\cir{\nu' z^2+\gamma' z}_\infty,
\]
with $\lambda'$, $\mu'$, $\nu'$ pairwise distinct.

Hence the generic representation has rank $3$, $k=3$ Fourier sphere coefficients, and the corresponding partition of the set of Stokes circles of $\bm{\breve{\Theta}}^{\rm IV}_{\rm gen}$ is given by
\begin{gather*}
N_1=N_1^-=\cir{\lambda' z^2+\alpha' z}_\infty,
\qquad\!\!\!
N_2=N_2^-=\cir{\mu' z^2 +\beta' z}_\infty,\!\!\!
\qquad
N_3=N_3^-=\cir{\nu' z^2+\gamma' z}_\infty.
\end{gather*}

Notice that the generic class of representations corresponds to the rank~3 Lax representation for Painlev\'e~IV found in \cite{joshi2007linearization}.

We observe that the three principal subtrees of the generic fission forest are isomorphic. This implies that the forests of the three nongeneric representations coincide so there are only two distinct forests, the generic one and a nongeneric one, corresponding to the standard representation.

\subsubsection{Diagram}

The corresponding diagram $\Gamma^{\rm IV}$ is the affine $A_2$ Dynkin diagram \cite{boalch2018wild, doucot2021diagrams}
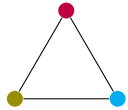
$$
\begin{tikzpicture}[scale=0.6]
\tikzstyle{vertex}=[circle,fill=cyan,minimum size=6pt,inner sep=0pt]
\tikzstyle{vertex_2}=[circle,fill=purple,minimum size=6pt,inner sep=0pt]
\tikzstyle{vertex_3}=[circle,fill=olive,minimum size=6pt,inner sep=0pt]
\node[vertex] (A) at (-30:1.3){} ;
\node[vertex_2] (B) at (90:1.3){} ;
\node[vertex_3] (C) at (210:1.3){};
\draw (A)--(B);
\draw (A)--(C);
\draw (B)--(C);
\end{tikzpicture}
$$
The two classes of representations correspond to different readings of the diagram, as indicated on the figure below
$$
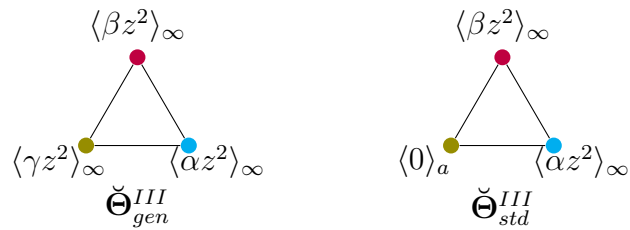
\begin{tikzpicture}[scale=0.6]
\tikzstyle{vertex}=[circle,fill=cyan,minimum size=6pt,inner sep=0pt]
\tikzstyle{vertex_2}=[circle,fill=purple,minimum size=6pt,inner sep=0pt]
\tikzstyle{vertex_3}=[circle,fill=olive,minimum size=6pt,inner sep=0pt]
\begin{scope}
\node[vertex] (A) at (-30:1.3){} ;
\node[vertex_2] (B) at (90:1.3){} ;
\node[vertex_3] (C) at (210:1.3){};
\draw (A)--(B);
\draw (A)--(C);
\draw (B)--(C);
\draw (-30:2.5) node {$\cir{\lambda' z^2+\alpha' z}_\infty$};
\draw (90:2.1) node {$\cir{\mu' z^2+\beta' z}_\infty$};
\draw (210:2.5) node {$\cir{\nu' z^2+\gamma' z}_\infty$};
\draw (0,-2.3) node {$\bm{\breve{\Theta}}^{\rm IV}_{\rm gen}$};
\end{scope}
\begin{scope}[xshift=8cm]
\node[vertex] (A) at (-30:1.3){} ;
\node[vertex_2] (B) at (90:1.3){} ;
\node[vertex_3] (C) at (210:1.3){};
\draw (A)--(B);
\draw (A)--(C);
\draw (B)--(C);
\draw (-30:2.5) node {$\cir{\alpha_2 z^2+\alpha_1 z}_\infty$};
\draw (90:2.1) node {$\cir{\beta_2 z^2+\beta_1 z}_\infty$};
\draw (210:2.5) node {$\cir{0}_a$};
\draw (0,-2.3) node {$\bm{\breve{\Theta}}^{\rm IV}_{\rm std}$};
\end{scope}
\end{tikzpicture}
$$

\subsection{Painlev\'e~V}

\subsubsection{Standard representation}

The standard Painlev\'e~V Lax representation has rank 2, one irregular singularity corresponding to a pole of order 2, and two regular singularities.

The corresponding fission forest $\mathbf F^{\rm V}_{\rm std}$ is the following:
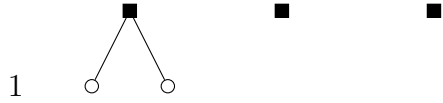
$$
\begin{tikzpicture}
\tikzstyle{mandatory}=[circle,fill=black,minimum size=5pt,draw, inner sep=0pt]
\tikzstyle{authorised}=[circle,fill=white,minimum size=5pt,draw,inner sep=0pt]
\tikzstyle{empty}=[circle,fill=black,minimum size=0pt,inner sep=0pt]
\tikzstyle{root}=[fill=black,minimum size=5pt,draw,inner sep=0pt]
\tikzstyle{indeterminate}=[circle,densely dotted,fill=white,minimum size=5pt,draw, inner sep=0pt]
 \draw (-1,1) node {$1$};
 \node[root] (R1) at (0.5,2){};
 \node[authorised] (A1) at (0,1){};
 \node[authorised] (B1) at (1,1){};
 \node[empty] (A0) at (0,0){};
 \node[empty] (B0) at (1,0){};
 \draw (R1)--(A1)--(A0);
 \draw (R1)--(B1)--(B0);
 \node[root] (R2) at (2.5,1){};
 \node[empty] (C0) at (2.5,0){};
 \node[root] (R3) at (4.5,1){};
 \node[empty] (D0) at (4.5,0){};
 \draw (R2)--(C0);
 \draw (R3)--(D0);
 \draw (2.5,-0.5) node {{\scriptsize 2}};
 \draw (4.5,-0.5) node {{\scriptsize 2}};
\end{tikzpicture}
$$
The irregular classes with fission forest $\mathbf F^{\rm V}_{\rm std}$ and minimal Poincar\'e--Katz rank are exactly those of the form
\[
\bm{{\Theta}}^{\rm V}_{\rm std}=2\cir{0}_a + 2\cir{0}_b +\cir{\alpha z_c^{-1}}_c + \cir{\beta z_c^{-1}}_c
\]
with $a, b, c\in \mathbb P^1$ pairwise distinct, $\alpha, \beta\in \mathbb C$, $\alpha\neq \beta$.

If $(E,\nabla)$ is a connection of this type, its generic class of nearby representations is of minimal rank if $c=\infty$, and if the conjugacy class in $\GL_2(\mathbb C)$ of its formal monodromy for the tame circles~$\cir{0}_a$ and~$\cir{0}_b$ has $1$ as an eigenvalue (we can always reduce to this case by applying a~M\"obius transformation and a twist). This corresponds to having a modified irregular class of the form
\[
\bm{\breve\Theta}^{\rm V}_{\rm std}=\cir{0}_a + \cir{0}_b+\cir{\alpha z}_\infty + \cir{\beta z}_\infty,
\]
with $\alpha\neq \beta$. Let us thus describe the classes of nearby representations of $\bm{\breve{\Theta}}^{\rm V}_{\rm std}$.

\subsubsection{Generic representation}

The short fission tree $\mathbf T^{\rm V}$ of $\bm{\breve\Theta}^{\rm V}_{\rm std}$ is

\begin{center}
\begin{tikzpicture}
\tikzstyle{mandatory}=[circle,fill=black,minimum size=5pt,draw, inner sep=0pt]
\tikzstyle{authorised}=[circle,fill=white,minimum size=5pt,draw,inner sep=0pt]
\tikzstyle{empty}=[circle,fill=black,minimum size=0pt,inner sep=0pt]
\tikzstyle{root}=[fill=black,minimum size=5pt,draw,inner sep=0pt]
\tikzstyle{indeterminate}=[circle,densely dotted,fill=white,minimum size=5pt,draw, inner sep=0pt]
 \draw (-2,1) node {$1$};
 \draw (-2,2) node {$2$};
 \node[root] (A3) at (1.5,3){};
 \node[authorised] (A2) at (0.5,2){};
 \node[authorised] (B2) at (2.5,2){};
 \node[authorised] (A1) at (0,1){};
 \node[authorised] (B1) at (1,1){};
 \node[authorised] (C1) at (2,1){};
 \node[authorised] (D1) at (3,1){};
 \node[empty] (A0) at (0,0){};
 \node[empty] (B0) at (1,0){};
 \node[empty] (C0) at (2,0){};
 \node[empty] (D0) at (3,0){};
 \draw (A2)--(A1)--(A0);
 \draw (A2)--(B1)--(B0);
 \draw (B2)--(C1)--(C0);
 \draw (B2)--(D1)--(D0);
 \draw (A3)--(A2);
 \draw (A3)--(B2);
\end{tikzpicture}
\end{center}

An irregular class of generic form has fission tree $\mathbf T^{\rm V}$ if and only if it is of the form
\[
\bm{\breve{\Theta}}^{\rm V}_{\rm gen}=\cir{\lambda' z^2 +\alpha' z}_\infty + \cir{\lambda' z^2+\beta' z}_\infty +\cir{\mu' z^2 +\gamma' z}_\infty + \cir{\mu' z^2+\delta' z}_\infty
\]
with $\lambda'\neq \mu'$, $\alpha' \neq \beta'$, $\gamma'\neq \delta'$.

The generic class of representations has thus rank 4, $k=2$ Fourier sphere coefficients, and the corresponding partition of the set of Stokes circles of $\bm{\breve{\Theta}}^{\rm V}_{\rm gen}$ is given by
\begin{gather*}
N_1=N_1^-=\big\{\cir{\lambda' z^2 +\alpha' z}_\infty, \cir{\lambda' z^2+\beta' z}_\infty\big\},
\\
N_2=N_2^-=\big\{\cir{\mu' z^2 +\gamma' z}_\infty,\cir{\mu' z^2+\delta' z}_\infty\big\}.
\end{gather*}
Furthermore, the two principal subtrees corresponding to $N_1$ and $N_2$ are isomorphic, so the two nongeneric classes of representations have in fact the same forest. In turn, there are only two distinct classes of basic representations, the generic one and the nongeneric one, corresponding to the standard Lax representation.

\subsubsection{Diagram}

The corresponding diagram $\Gamma^V$ is the affine $A_3$ Dynkin diagram \cite{boalch2018wild, doucot2021diagrams}

\begin{center}
\begin{tikzpicture}[scale=0.6]
\tikzstyle{vertex}=[circle,fill=cyan,minimum size=6pt,inner sep=0pt]
\tikzstyle{vertex_2}=[circle,fill=purple,minimum size=6pt,inner sep=0pt]
\node[vertex_2] (A) at (1,1){} ;
\node[vertex] (B) at (-1,1){} ;
\node[vertex_2] (C) at (-1,-1){};
\node[vertex] (D) at (1,-1){};
\draw (A)--(B)--(C)--(D)--(A);
\end{tikzpicture}
\end{center}
The two classes of representations correspond to different readings of the diagram, as indicated on the figure below
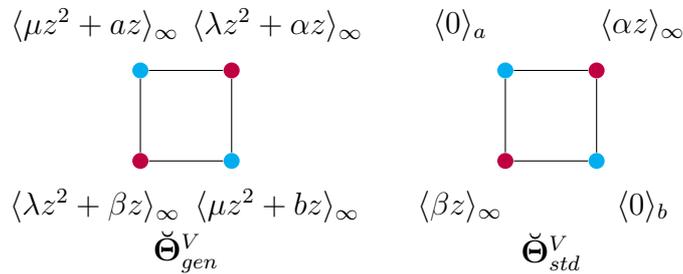
\begin{center}
\begin{tikzpicture}[scale=0.6]

\tikzstyle{vertex}=[circle,fill=cyan,minimum size=6pt,inner sep=0pt]
\tikzstyle{vertex_2}=[circle,fill=purple,minimum size=6pt,inner sep=0pt]
\begin{scope}
\node[vertex_2] (A) at (1,1){} ;
\node[vertex] (B) at (-1,1){} ;
\node[vertex_2] (C) at (-1,-1){};
\node[vertex] (D) at (1,-1){};
\draw (A)--(B)--(C)--(D)--(A);
\draw (2,2) node {$\cir{\lambda' z^2 +\alpha' z}_\infty$};
\draw (-2,-2) node {$\cir{\lambda' z^2 +\beta' z}_\infty$};
\draw (-2,2) node {$\cir{\mu' z^2 +\gamma' z}_\infty$};
\draw (2,-2) node {$\cir{\mu' z^2 +\delta' z}_\infty$};
\draw (0,-3) node {$\bm{\breve\Theta}^{\rm V}_{\rm gen}$};
\end{scope}
\begin{scope}[xshift=8cm]
\node[vertex_2] (A) at (1,1){} ;
\node[vertex] (B) at (-1,1){} ;
\node[vertex_2] (C) at (-1,-1){};
\node[vertex] (D) at (1,-1){};
\draw (A)--(B)--(C)--(D)--(A);
\draw (2,2) node {$\cir{\alpha z}_\infty$};
\draw (-2,-2) node {$\cir{\beta z}_\infty$};
\draw (-2,2) node {$\cir{0}_a$};
\draw (2,-2) node {$\cir{0}_b$};
\draw (0,-3) node {$\bm{\breve\Theta}^{\rm V}_{\rm std}$};
\end{scope}
\end{tikzpicture}
\end{center}

\subsection{Painlev\'e~VI}

\subsubsection{Standard representation}

The standard rank 2 Lax representation for Painlev\'e~VI has 4 regular singularities.

The corresponding fission forest $\mathbf F^{\rm VI}_{\rm std}$ is the following (all fission trees are trivial):
\begin{center}
\begin{tikzpicture}
\tikzstyle{mandatory}=[circle,fill=black,minimum size=5pt,draw, inner sep=0pt]
\tikzstyle{authorised}=[circle,fill=white,minimum size=5pt,draw,inner sep=0pt]
\tikzstyle{empty}=[circle,fill=black,minimum size=0pt,inner sep=0pt]
\tikzstyle{root}=[fill=black,minimum size=5pt,draw,inner sep=0pt]
\tikzstyle{indeterminate}=[circle,densely dotted,fill=white,minimum size=5pt,draw, inner sep=0pt]
 \node[root] (R1) at (0,1){};
 \node[root] (R2) at (2,1){};
 \node[root] (R3) at (4,1){};
 \node[root] (R4) at (6,1){};
 \node[empty] (A0) at (0,0){};
 \node[empty] (B0) at (2,0){};
 \node[empty] (C0) at (4,0){};
 \node[empty] (D0) at (6,0){};
 \draw (R1)--(A0);
 \draw (R2)--(B0);
 \draw (R3)--(C0);
 \draw (R4)--(D0);
 \draw (0,-0.5) node {{\scriptsize 2}};
 \draw (2,-0.5) node {{\scriptsize 2}};
 \draw (4,-0.5) node {{\scriptsize 2}};
 \draw (6,-0.5) node {{\scriptsize 2}};
 %\draw (0,-1) node {$\cir{0}_\infty$};
 %\draw (2,-1) node {$\cir{0}_a$};
 %\draw (4,-1) node {$\cir{0}_b$};
 %\draw (6,-1) node {$\cir{0}_c$};
\end{tikzpicture}
\end{center}
The irregular classes with fission forest $\mathbf F^{\rm VI}_{\rm std}$ are those of the form
\[
\bm{\Theta}^{\rm VI}_{\rm std}=2\cir{0}_a + 2\cir{0}_b+ 2\cir{0}_c + 2\cir{0}_d.
\]
with $a,b,c,d\in \mathbb P^1$ pairwise distinct.

Given a connection of this type, its generic class of nearby representations has minimal rank if one of the singularities, say $d$, is at infinity, and if for the other singularities $a$, $b$, $c$, the formal monodromy of the connection has 1 as an eigenvalue (once again we can always reduce to this situation by applying a M\"obius transformation and a twist). This corresponds to a modified irregular class of the form
\[
\bm{\breve{\Theta}}^{\rm VI}_{\rm std}=2\cir{0}_\infty + \cir{0}_a+ \cir{0}_b + \cir{0}_c,
\]
with $a, b,c\in \mathbb P^1\smallsetminus\{\infty\}$ pairwise distinct. Let us thus determine the classes of nearby representations of $\bm{\breve{\Theta}}^{\rm VI}_{\rm std}$.

\subsubsection{Generic representation}

The short fission tree $\mathbf T^{\rm VI}$ of \smash{$\bm{\breve{\Theta}}^{\rm VI}_{\rm std}$} is
\begin{center}
\begin{tikzpicture}
\tikzstyle{mandatory}=[circle,fill=black,minimum size=5pt,draw, inner sep=0pt]
\tikzstyle{authorised}=[circle,fill=white,minimum size=5pt,draw,inner sep=0pt]
\tikzstyle{empty}=[circle,fill=black,minimum size=0pt,inner sep=0pt]
\tikzstyle{root}=[fill=black,minimum size=5pt,draw,inner sep=0pt]
\tikzstyle{indeterminate}=[circle,densely dotted,fill=white,minimum size=5pt,draw, inner sep=0pt]
 \draw (-2,1) node {$1$};
 \draw (-2,2) node {$2$};
 \node[root] (A3) at (0,3){};
 \node[authorised] (A2) at (1,2){};
 \node[authorised] (B2) at (-1,2){};
 \node[authorised] (A1) at (0,1){};
 \node[authorised] (B1) at (1,1){};
 \node[authorised] (C1) at (2,1){};
 \node[authorised] (D1) at (-1,1){};
 \node[empty] (A0) at (0,0){};
 \node[empty] (B0) at (1,0){};
 \node[empty] (C0) at (2,0){};
 \node[empty] (D0) at (-1,0){};
 \draw (A2)--(A1)--(A0);
 \draw (A2)--(B1)--(B0);
 \draw (A2)--(C1)--(C0);
 \draw (B2)--(D1)--(D0);
 \draw (A3)--(A2);
 \draw (A3)--(B2);
 \draw (-1,-0.4) node {{\scriptsize $2$}};
\end{tikzpicture}
\end{center}
An irregular class of generic form has fission tree $\mathbf T^{\rm VI}$ if and only if it is of the form
\[
\bm{\breve{\Theta}}^{\rm VI}_{\rm gen}=2\cir{\lambda' z^2+\alpha' z}_\infty + \cir{\mu' z^2 + \beta' z}_\infty+ \cir{\mu' z^2 + \gamma' z}_\infty + \cir{\mu' z^2 + \delta' z}_\infty,
\]
with $\lambda'\neq \mu'$, and $\beta'$, $\gamma'$, $\delta'$ pairwise distinct.

The generic class of representations has rank $5$, its number of Fourier sphere coefficients is~${k=2}$, and the corresponding partition of the set of Stokes circles of $\bm{\breve{\Theta}}^{\rm VI}_{\rm gen}$ is given by
\begin{gather*}
 N_1=N_1^-=\big\{\cir{\lambda' z^2+\alpha' z}\big\},
\\
N_2=N_2^-=\big\{\cir{\mu' z^2 + \beta' z},\cir{\mu' z^2 + \gamma' z}, \cir{\mu' z^2+\delta' z}\big\}.
\end{gather*}

\subsubsection{Nongeneric representations}

Since the two principal subtrees are not isomorphic, there are two distinct classes of nongeneric nearby representations. The standard representation corresponds to the nongeneric class of representations for $i=2$.

The other nongeneric class of representations, for $i=1$, has the following fission forest $\mathbf F^{\rm VI}_{1}$:
$$
\begin{tikzpicture}
\tikzstyle{mandatory}=[circle,fill=black,minimum size=5pt,draw, inner sep=0pt]
\tikzstyle{authorised}=[circle,fill=white,minimum size=5pt,draw,inner sep=0pt]
\tikzstyle{empty}=[circle,fill=black,minimum size=0pt,inner sep=0pt]
\tikzstyle{root}=[fill=black,minimum size=5pt,draw,inner sep=0pt]
\tikzstyle{indeterminate}=[circle,densely dotted,fill=white,minimum size=5pt,draw, inner sep=0pt]
 \draw (-2,1) node {$1$};
 \draw (-2,2) node {$2$};
 \node[root] (R) at (1,2){};
 \node[root] (A1) at (-1,1){};
 \node[authorised] (B1) at (0,1){};
 \node[authorised] (C1) at (1,1){};
 \node[authorised] (D1) at (2,1){};
 \node[empty] (A0) at (-1,0){};
 \node[empty] (B0) at (0,0){};
 \node[empty] (C0) at (1,0){};
 \node[empty] (D0) at (2,0){};
 \draw (A1)--(A0);
 \draw (R)--(B1)--(B0);
 \draw (R)--(C1)--(C0);
 \draw (R)--(D1)--(D0);
 \draw (-1,-0.4) node {{\scriptsize $3$}};
\end{tikzpicture}
$$
Explicitly, elements of the orbit $\SL_2(\mathbb C)\cdot \bm{\breve{\Theta}}^{\rm VI}_{\rm gen}$ belonging to this class of representations are of the form
\[
\bm{\breve{\Theta}}^{\rm VI}_1=2\cir{0}_{a''} + \cir{\alpha'' z}_\infty+ \cir{\beta'' z}_\infty + \cir{\gamma'' z}_\infty,
\]
with $a''\in \mathbb P^1\smallsetminus\{\infty\}$, and $\alpha'', \beta'',\gamma''\in \mathbb C$ pairwise distinct.
The corresponding non-modified irregular class is thus
\[
\bm{\Theta}^{\rm VI}_1=3\cir{0}_{a''} + \cir{\alpha'' z}_\infty+ \cir{\beta'' z}_\infty + \cir{\gamma'' z}_\infty,
\]
and its fission forest is indeed $\mathbf F^{\rm VI}_{1}$.

Notice that this class of representations corresponds to the well-known Harnad dual \cite{harnad1994dual} of the standard Lax representation.

\subsubsection{Diagram}

The corresponding diagram $\Gamma^{\rm VI}$ is the affine $D_4$ Dynkin diagram famously associated to Pain\-lev\'e~VI
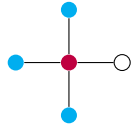
$$
\begin{tikzpicture}[scale=0.35]
\tikzstyle{vertex}=[circle,fill=cyan,minimum size=6pt,inner sep=0pt]
\tikzstyle{vertex_2}=[circle,fill=purple,minimum size=6pt,inner sep=0pt]
\tikzstyle{vertex_3}=[circle,fill=white,draw,minimum size=6pt,inner sep=0pt]
\node[vertex_2] (A) at (0,0){} ;
\node[vertex] (B) at (0,2){} ;
\node[vertex] (C) at (0,-2){};
\node[vertex_3] (D) at (2,0){};
\node[vertex] (E) at (-2,0){};
\draw (A)--(B);
\draw (A)--(C);
\draw (A)--(D);
\draw (A)--(E);
%\draw (A)node[below left]{$2$};
\end{tikzpicture}
$$
Here, the right node is not coloured, to indicate that it does not belong to the core diagram, but comes from the conjugacy class of the formal monodromy of the Stokes circle with multiplicity~2, corresponding to the central vertex.

The classes of representations correspond to three different readings of the diagram, as indicated on the figure below
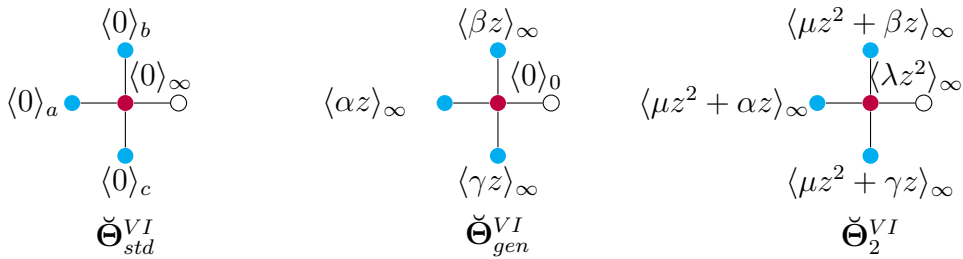
$$
\begin{tikzpicture}[scale=0.35]
\tikzstyle{vertex}=[circle,fill=cyan,minimum size=6pt,inner sep=0pt]
\tikzstyle{vertex_2}=[circle,fill=purple,minimum size=6pt,inner sep=0pt]
\tikzstyle{vertex_3}=[circle,fill=white,draw,minimum size=6pt,inner sep=0pt]
\begin{scope}
\node[vertex_2] (A) at (0,0){} ;
\node[vertex] (B) at (0,2){} ;
\node[vertex] (C) at (0,-2){};
\node[vertex_3] (D) at (2,0){};
\node[vertex] (E) at (-2,0){};
\draw (A)--(B);
\draw (A)--(C);
\draw (A)--(D);
\draw (A)--(E);
\draw (-3.5,0) node {$\cir{0}_a$};
\draw (0,3) node {$\cir{0}_b$};
\draw (0,-3) node {$\cir{0}_c$};
\draw (1.3,1) node {$\cir{0}_{\infty}$};
\draw (0,-5) node {$\bm{\breve{\Theta}}^{\rm VI}_{\rm std}$};
\end{scope}
\begin{scope}[xshift=14cm]
\node[vertex_2] (A) at (0,0){} ;
\node[vertex] (B) at (0,2){} ;
\node[vertex] (C) at (0,-2){};
\node[vertex_3] (D) at (2,0){};
\node[vertex] (E) at (-2,0){};
\draw (A)--(B);
\draw (A)--(C);
\draw (A)--(D);
\draw (A)--(E);
\draw (-5,0) node {$\cir{\alpha'' z}_\infty$};
\draw (0,3) node {$\cir{\beta'' z}_\infty$};
\draw (0,-3) node {$\cir{\gamma'' z}_\infty$};
\draw (1.4,1) node {$\cir{0}_{a''}$};
\draw (0,-5) node {$\bm{\breve{\Theta}}^{\rm VI}_{1}$};
\end{scope}
\begin{scope}[xshift=28cm]
\node[vertex_2] (A) at (0,0){} ;
\node[vertex] (B) at (0,2){} ;
\node[vertex] (C) at (0,-2){};
\node[vertex_3] (D) at (2,0){};
\node[vertex] (E) at (-2,0){};
\draw (A)--(B);
\draw (A)--(C);
\draw (A)--(D);
\draw (A)--(E);
\draw (-5.7,0) node {{\scriptsize $\cir{\mu' z^2+\beta' z}_\infty$}};
\draw (0,3) node {{\scriptsize$\cir{\mu' z^2+\gamma' z}_\infty$}};
\draw (0,-3) node {{\scriptsize$\cir{\mu' z^2+\delta' z}_\infty$}};
\draw (2.5,1) node {{\scriptsize $\cir{\lambda' z^2+\alpha' z}_{\infty}$}};
\draw (0,-5) node {$\bm{\breve{\Theta}}^{\rm VI}_{\rm gen}$};
\end{scope}
\end{tikzpicture}
$$

\begin{Remark}
Let us make a few comments about all these Painlev\'e examples. Notice that, in several cases, the short fission tree breaks part of the symmetry of the equation, for example, for Painlev\'e~VI it breaks the symmetry of the 4 regular singular points, and for Painlev\'e~III it breaks the symmetry between the two singularities. This is because the construction depends on the location of the point at infinity, and because to have a generic representation with minimal rank, we need to have one of the singularities at infinity (which coincides with the standard choice). On the other hand, in several cases, the short fission tree makes apparent some symmetry that was not directly visible from the fission forest: for example, for Painlev\'e~IV the generic tree makes it clear that there is a symmetry between the two Stokes circles at infinity and the two singularities at finite distance, and for degenerate Painlev\'e~III, it makes clear that the forests of two nongeneric classes of nearby representations coincide. The basic reason underlying these observations is that the different combinatorial invariants behave in a~nice way with respect to different natural operations on connections: the fission forest is invariant under M\"obius transformations and twisting by any rank one connection, but not invariant under~$\SL_2(\mathbb C)$, while the diagram is invariant under~$\SL_2(\mathbb C)$, but not under M\"obius transformations changing the location of infinity or twists with rank one connections having singularities at finite distance. The takeaway from this is that, for any given nonabelian Hodge space, it can be useful to combine all these invariants to understand as many symmetries of the situation as possible.
\end{Remark}

\subsection*{Acknowledgements} I am funded by the PNRR Grant CF 44/14.11.2022, ``Cohomological Hall Algebras of Smooth Surfaces and Applications''. I~also acknowledge the funding of FCiências.ID, the great working environment of the Grupo de F\'{\i}sica Matem\'atica, where most of this work was carried out, and Philip Boalch, Gabriele Rembado, Giordano Cotti, Giulio Ruzza, Davide Masoero, Gabriele Degano, Julian Barrag\'an Amado, for useful discussions. I am also grateful to the referees for their useful comments which helped improve the quality of the manuscript.

\pdfbookmark[1]{References}{ref}
\LastPageEnding

\end{document}